\documentclass[AMA,STIX1COL]{WileyNJD-v2}

\articletype{RESEARCH ARTICLE}%

\received{}
\revised{}
\accepted{}

\usepackage{subfigure}
\usepackage{caption}
\usepackage{array}
\newcommand{\tabincell}[2]{\begin{tabular}{@{}#1@{}}#2\end{tabular}}  

\begin{document}

\title{Meshless Fragile Points Methods Based on Petrov-Galerkin Weak-Forms for Transient Heat Conduction Problems in Complex Anisotropic Nonhomogeneous Media.}

\author[1]{Yue Guan*}

\author[1]{Satya N. Atluri}

\authormark{GUAN \textsc{et al}}

\address[1]{\orgdiv{Department of Mechanical Engineering}, \orgname{Texas Tech University}, \orgaddress{\state{TX}, \country{United States}}}

\corres{*Yue Guan, Department of Mechanical Engineering, Texas Tech University, Lubbock, TX 79415, United States. \email{yuguan@ttu.edu}}


\abstract[Summary]{Three kinds of Fragile Points Methods based on Petrov-Galerkin weak-forms (PG-FPMs) are proposed for analyzing heat conduction problems in nonhomogeneous anisotropic media. This is a follow-up of the previous study on the original FPM based on a symmetric Galerkin weak-form. The trial function is piecewise-continuous, written as local Taylor expansions at the Fragile Points. A modified Radial Basis Function-based Differential Quadrature (RBF-DQ) method is employed for establishing the local approximation. The Dirac delta function, Heaviside step function, and the local fundamental solution of the governing equation are alternatively used as test functions. Vanishing or pure contour integral formulation in subdomains or on local boundaries can be obtained. Extensive numerical examples in 2D and 3D are provided as validations. The collocation method (PG-FPM-1) is superior in transient analysis with arbitrary point distribution and domain partition. The finite volume method (PG-FPM-2) shows the best efficiency, saving 25\% to 50\% computational time comparing with the Galerkin FPM. The singular solution method (PG-FPM-3) is highly efficient in steady-state analysis. The anisotropy and nonhomogeneity give rise to no difficulties in all the methods. The proposed PG-FPM approaches represent an improvement to the original Galerkin FPM, as well as to other meshless methods in earlier literature.}


\maketitle

\section{Introduction}

With an increasing application of new materials in various engineering fields, the necessity to solve heat conduction problems in continuously nonhomogeneous and anisotropic media is growing \cite{Sladek2005, Quint2011, Zhang2019}. Despite the great success of the Finite Element Method (FEM) and Finite Volume Method (FVM) in commercial software packages, they still suffer from inherent drawbacks such as locking, mesh distortion, etc., and there is still an increasing interest in developing new advanced methods.

Among the most popular numerical tools in heat conduction analysis, the Finite Element Method (FEM) \cite{Zienkiewicz2005} based on a Galerkin weak-form can be applied in complex domains. Yet it usually demands a high-quality mesh, and the accuracy can be threatened by mesh distortion. The Finite Difference Method (FDM) \cite{Chen1999}, on the other hand, is efficient but only applicable in regular domains with uniform meshes. The Finite Volume Method (FVM) \cite{Chai1994, Gersborg-Hansen2006, Li2012, Rapp2017} combines the versatility of the FEM in unstructured grids and the efficiency of the FDM, hence is widely used in engineering practice. The conventional FVM strictly satisfies the conservation law, which is essential in fluid mechanics. However, the FVM based a classic two-point flux formula is inconsistent in anisotropic problems \cite{Droniou2006, Eymard2019}. Multiple-point flux approximations \cite{Liu2015, Prestini2017, Eymard2019}, e.g., generalized finite difference method, have to be employed. And the same as the FEM, the FVM also suffers from the mesh distortion problem and may have difficulties in analyzing problems with crack developments. The Boundary Element Method (BEM) \cite{Wrobel2003, Wen2009} is another mature technique in heat conduction analysis. However, the non-availability of infinite space fundamental solutions in continuously nonhomogeneous anisotropic media limits its application. Besides, the BEM leads to fully populated matrix which is computationally expensive, especially in transient problems.

Meshless methods, on the other hand, are increasingly popular in recent decades due to their high adaptivity and low cost in generating meshes. Representative meshless methods based on symmetric Galerkin weak-forms contain the Diffuse Element Method (DEM) \cite{Nayroles1992}, the Element-Free Galerkin (EFG) method \cite{Belytschko1994}, etc. Alternatively, the trial and test functions can be chosen from different function spaces. This is known as a Petrov-Galerkin scheme. An appropriate choice of test function could improve the performance of the approach. For example, the Petrov-Galerkin Diffuse Element Method (PG DEM) \cite{Krongauz1997} passes the patch test and exhibits better accuracy and convergence rate than the original DEM. The Local Boundary Integral Equation (LBIE) Method \cite{Zhu1999} using local fundamental solutions as test functions bypasses the problem of global fundamental solutions in nonhomogeneous anisotropic materials in the BEM. The Meshless Local Petrov-Galerkin (MLPG) method proposed by Atluri and Zhu \cite{Atluri1998} exploiting multiple test functions including Dirac delta function, Heaviside step function, and local fundamental solution of the governing equations simplifies the computational process by constructing pure contour or vanishing integrals on local boundaries. The MLPG method is truly meshless and has been employed in analyzing 2D and 3D heat conduction problems with anisotropy and nonhomogeneity \cite{Sladek2008, Shibahara2011}.

Nevertheless, all the meshless methods mentioned above are based on Moving Least Squares (MLS), Radial Basis Function (RBF) or other complicated local approximations. As a result, the numerical integrations in either Galerkin or Petrov-Galerkin weak-forms are tedious and extremely complicated. The trial function based on MLS approximation also lacks Delta-function property, hence the essential boundary conditions cannot be enforced directly. Dong et al. \cite{ LeitingDongTianYangKaileiWang2019} proposed a Fragile Points Method (FPM) based on symmetric Galerkin weak-form which avoids the complicated numerical integration using local, polynomial, piecewise-continuous trial functions. Guassian quadrature scheme with only one integration point in each subdomain is sufficient most of the time. The trial function also satisfies the Kronecker-delta property. Thus the essential boundary conditions can be applied strongly. The work is further extended to 2D elasticity problems \cite{Yang2021}, and 2D and 3D heat conduction problems \cite{Guan2020, Guan2020a}. The current paper is a follow-up of the previous study \cite{Guan2020, Guan2020a}. Several improved versions of the FPM are carried out in this paper, based on Petrov-Galerkin formulations. When the test functions are chosen as Dirac delta function, Heaviside step function, and local fundamental solutions (assuming locally homogeneous material properties), integrals in the subdomains or on the local boundaries may vanish or become pure contour integrals. Therefore, a higher computational efficiency can be expected. Besides, the final discretized algebraic equations in the PG-FPMs are easily ascribed a physical interpretation. A comparison of the most popular numerical methods in heat conduction analysis and the original and improved FPM approaches is presented in Table~\ref{table:Comp_1}.

The following discussion begins with a brief description of the boundary value problem under study. The local, polynomial, discontinuous trial functions as well as the local approximations are introduced in Section~\ref{sec:TF}. The formulation and implementation of three kinds of Fragile Points Methods based on Petrov-Galerkin weak-forms (PG-FPMs) are developed in Section~\ref{sec:formulation}, followed by a brief discussion about the time discretization method in Section~\ref{sec:time}. In Section~\ref{sec:2D} and \ref{sec:3D}, a number of 2D and 3D numerical examples are presented respectively to illustrate the accuracy and efficiency of the three proposed PG-FPM approaches. A discussion on the penalty parameters used in the PG-FPMs is given at last in Section~\ref{sec:PS}.

\begin{table}[htbp]
\caption{A comparison of the PG-FPMs and other numerical methods.}
\centering
{
\begin{tabular}{  m{70pt}  m{70pt}<{\centering}  m{70pt}<{\centering}  m{70pt}<{\centering}  m{70pt}<{\centering}  m{70pt}<{\centering} }
\toprule
\textbf{Method} & \textbf{FEM} & \textbf{FVM} & \textbf{EFG} & \textbf{MLPG} & \textbf{FPM \& PG-FPMs} \\
\midrule
Trial functions & $C^0$ polynomial & Polynomial & MLS, etc. & MLS, etc. &  Discontinuous polynomial \\
\midrule
Test functions & As above & Heaviside step function & As above & MLS, weight function, Dirac delta function, Heaviside step function, local fundamental solution, $\cdots$ & Discontinuous polynomial, Dirac delta function, Heaviside step function, local fundamental solution, $\cdots$\\
\midrule
Weak-form & Galerkin & Petrov-Galerkin & Galerkin & Petrov-Galerkin & Galerkin \& Petrov-Galerkin \\
\midrule
Weak-form integration & Simple (Gaussian quadrature) & Simple (Gaussian quadrature) & Very complicated & Very complicated & Simple (Gaussian quadrature) \\
\midrule
Delta function property & Yes & Yes & Yes/No & Yes/No & Yes \\
\bottomrule
\end{tabular}}
\label{table:Comp_1}
\end{table}

\section{The heat conduction problem and governing equation} \label{sec:governing_eqn}

In this study, we consider a boundary value problem for heat conduction in a continuously nonhomogeneous anisotropic medium, which is described by the following governing equation:
\begin{align}
\begin{split}
\rho(\mathbf{x}) c(\mathbf{x})  \frac{\partial u}{\partial t} (\mathbf{x}, t) = \nabla \cdot \left[ \mathbf{k} \nabla u (\mathbf{x}, t) \right] + Q(\mathbf{x}, t), \quad \mathbf{x} \in \Omega, t \in \left[ 0, T \right],
\label{eqn:governing}
\end{split}
\end{align}
where $\Omega$ is the entire 2D or 3D domain under study, the spatial coordinate $\mathbf{x} = \left[ x , y \right]^\mathrm{T}$ in 2D or $\left[ x , y , z \right]^\mathrm{T}$ in 3D, $u(\mathbf{x}, t)$ is the temperature field, and $Q(\mathbf{x}, t)$ is the density of heat sources. $\mathbf{k} (\mathbf{x})$, $\rho (\mathbf{x})$ and $c(\mathbf{x})$ are the thermal conductivity tensor, mass density and specific heat capacity of the medium.

The following boundary and initial conditions are assumed:
\begin{align}
 \text{1).} \;  & \text{Dirichlet bc}:  & u(\mathbf{x}, t) = \widetilde{u}_D (\mathbf{x}, t),& \quad \text{on} \; \Gamma_D ,&  \label{eqn:BC1} \\
 \text{2).} \;  & \text{Neumann bc}: & q(\mathbf{x}, t) := \mathbf{n}^\mathrm{T} \mathbf{k} \nabla u(\mathbf{x}, t)  =  \widetilde{q}_N (\mathbf{x}, t),& \quad \text{on} \; \Gamma_N ,& \label{eqn:BC2} \\
 \text{3).} \; & \text{Robin (convective) bc}: & q(\mathbf{x}, t) = h(\mathbf{x}) \left[ \widetilde{u}_R (\mathbf{x}) - u(\mathbf{x}, t) \right],& \quad \text{on} \; \Gamma_R ,&\label{eqn:BC3} \\
 \text{4).} \; & \text{Symmetric bc}: & \mathbf{n}^\mathrm{T} \nabla u(\mathbf{x}, t) = 0,& \quad \text{on} \; \Gamma_S ,& \label{eqn:BC4} \\
  \text{5).} \; & \text{Initial condition}: & \left. u(\mathbf{x}, t) \right|_{t=0} = u(\mathbf{x}, 0) ,& \quad \text{in} \; \Omega \cup \partial \Omega, &
 \label{eqn:IC}
\end{align}
where the global boundary $\partial \Omega = \Gamma_D \cup \Gamma_N \cup \Gamma_R \cup \Gamma_S$, $\mathbf{n}$ is the unit outward normal vector of $\partial \Omega$, $h(\mathbf{x})$ is the heat transfer coefficient, and $\widetilde{u}_R (\mathbf{x})$ is the temperature of the medium outside the convective boundary. Note that for an isotropic problem, the symmetric condition is equivalent to Neumann boundary condition with $\widetilde{q}_N = 0$.

\section{Trial functions and meshless approximations} \label{sec:TF}

In this section, we introduce the local, polynomial, point-based discontinuous trial functions used in the Petrov-Galerkin Fragile Points Methods (PG-FPMs). A meshless local Radial Basis Function-based Differential Quadrature (RBF-DQ) method is employed to approximate the derivatives at randomly scattered points. More details can be found in \cite{LeitingDongTianYangKaileiWang2019, Guan2020, Yang2021}.

\begin{figure}[htbp] 
  \centering 
    \subfigure[]{ 
    \label{fig:Schem_2D} 
    \includegraphics[width=0.48\textwidth]{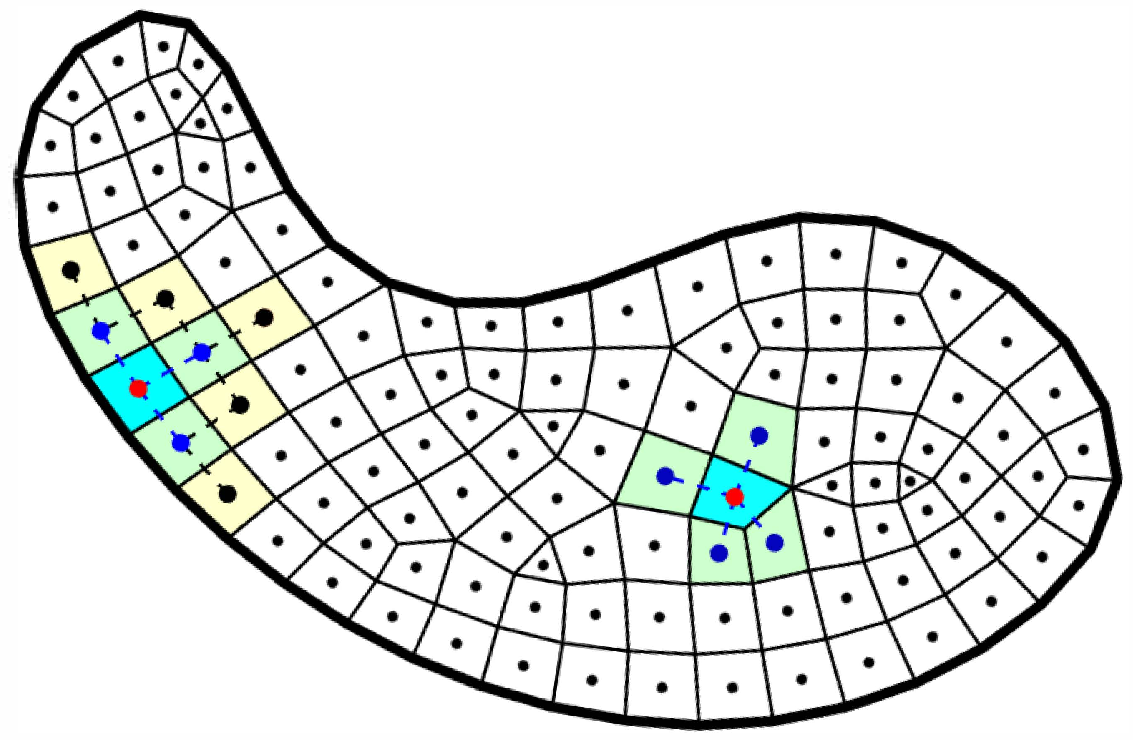}}  
    \subfigure[]{ 
    \label{fig:Schem_3D} 
    \includegraphics[width=0.48\textwidth]{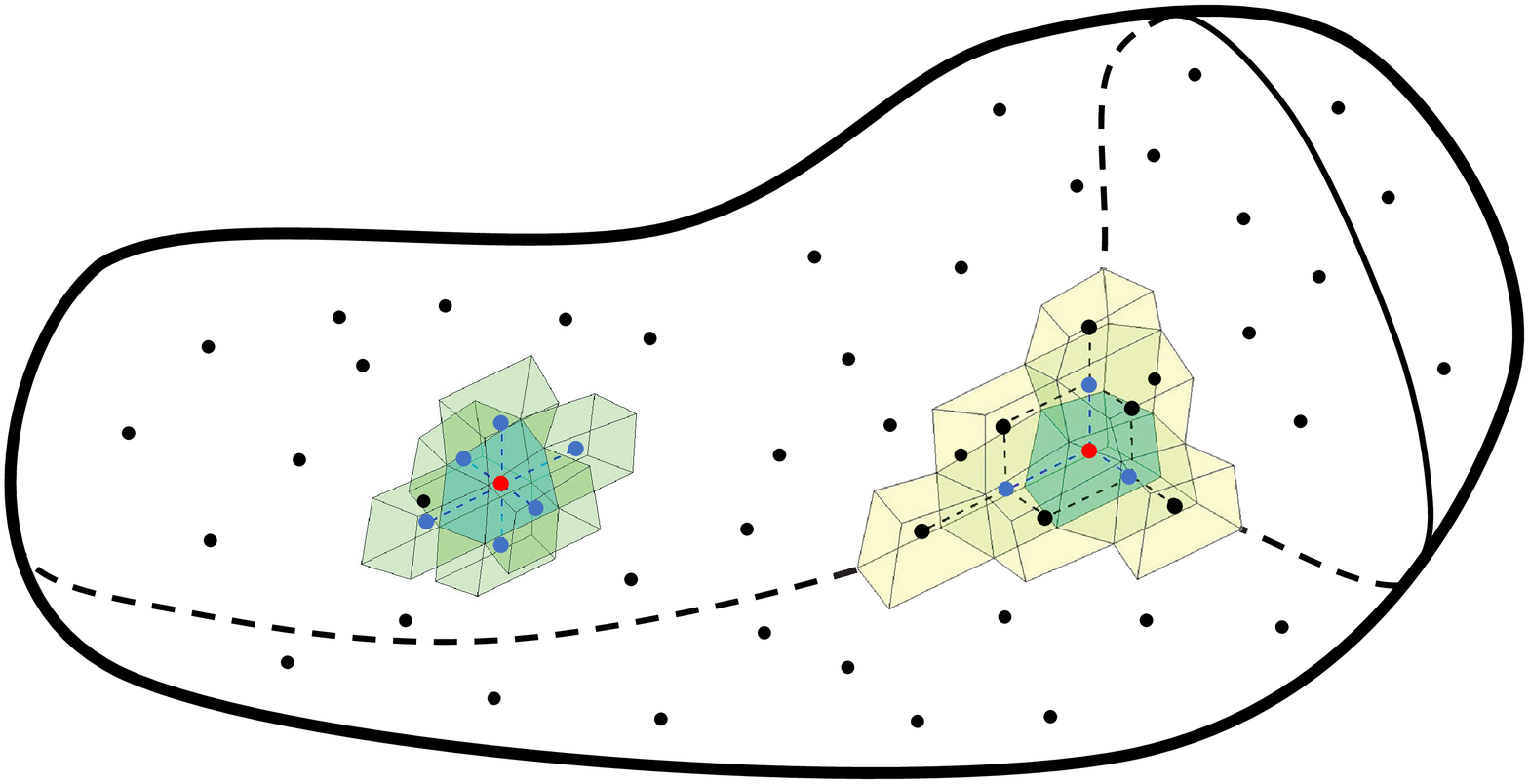}}  
  \caption{The domain $\Omega$ and its partitions. (a) 2D domain with a mixed quadrilateral and triangular partition (converted from FEA meshing). (b) 3D domain with a hexahedron partition (partially shown).} 
  \label{fig:Schem} 
\end{figure}

In the PG-FPMs, with a set of random Fragile Points, the global domain is partitioned into several confirming and nonoverlapping subdomains, within which only one Internal Point exists. A variety of partitioning schemes are available, e.g., the Voronoi Diagram partition \cite{Voronoi1908}, quadrilateral and triangular partition (in 2D), tetrahedron and hexahedron partition (in 3D), etc. In practice, the conventional FEM meshing can also be converted into FPM subdomains, while the internal Fragile Points are defined as the centroid of each FEM geometrical element. For example, Figure~\ref{fig:Schem_2D} shows a mixed quadrilateral and triangular partition and the corresponding Fragile Point distribution converted from ABAQUS preprocessing results.

In each subdomain, the {\it point-based} trial function ($u_h$) can be written as a Taylor expansion at the corresponding Internal Point. For instance, in subdomain $E_0$:
\begin{align}
\begin{split}
u_h (\mathbf{x}) = u_0 + \left(\mathbf{x} - \mathbf{x}_0\right)^\mathrm{T} \nabla u \Big| _{P_0} + \frac{1}{2}  \left(\mathbf{x} - \mathbf{x}_0 \right)^\mathrm{T} \mathbf{H} \left( u \right) \Big| _{P_0} \left(\mathbf{x} - \mathbf{x}_0 \right), \quad \mathbf{x} \in E_0,
\end{split} \label{eq:trial_0}
\end{align}
where $P_0$ is the Internal Point within $E_0$, $u_0$ is the value of $u_h$ at $P_0$, and $\mathbf{x}_0$ is the location vector of $P_0$. $\mathbf{H} \left( u \right)$ is the Hessian matrix of the temperature field, i.e., $\mathbf{H} \left( u \right) = \frac{\partial}{\partial \mathbf{x}} \left( \nabla u \right)$.

The gradient and Hessian matrix of temperature at point $P_0$ can then be approximated by the value of $u_h$ at several supporting Points. In the present work, as shown in Figure~\ref{fig:Schem}, the support of a given point $P_0 \in E_0$ is defined to involve all the nearest neighboring points in subdomains sharing boundaries with $E_0$. To ensure the accuracy of the second derivatives, at least three supporting points in each direction are required. Thus for Fragile Points close to the global boundary $\partial \Omega$, the second neighboring points (shown yellow in Figure~\ref{fig:Schem}) should also be taken into consideration. These supporting points are named as $P_1 \left( \mathbf{x}_1 \right), P_2 \left( \mathbf{x}_2 \right), \cdots ,P_m \left( \mathbf{x}_m \right)$.

\begin{figure}[htbp] 
  \centering 
    \subfigure[]{ 
    \label{fig:Support_01} 
    \includegraphics[width=0.48\textwidth]{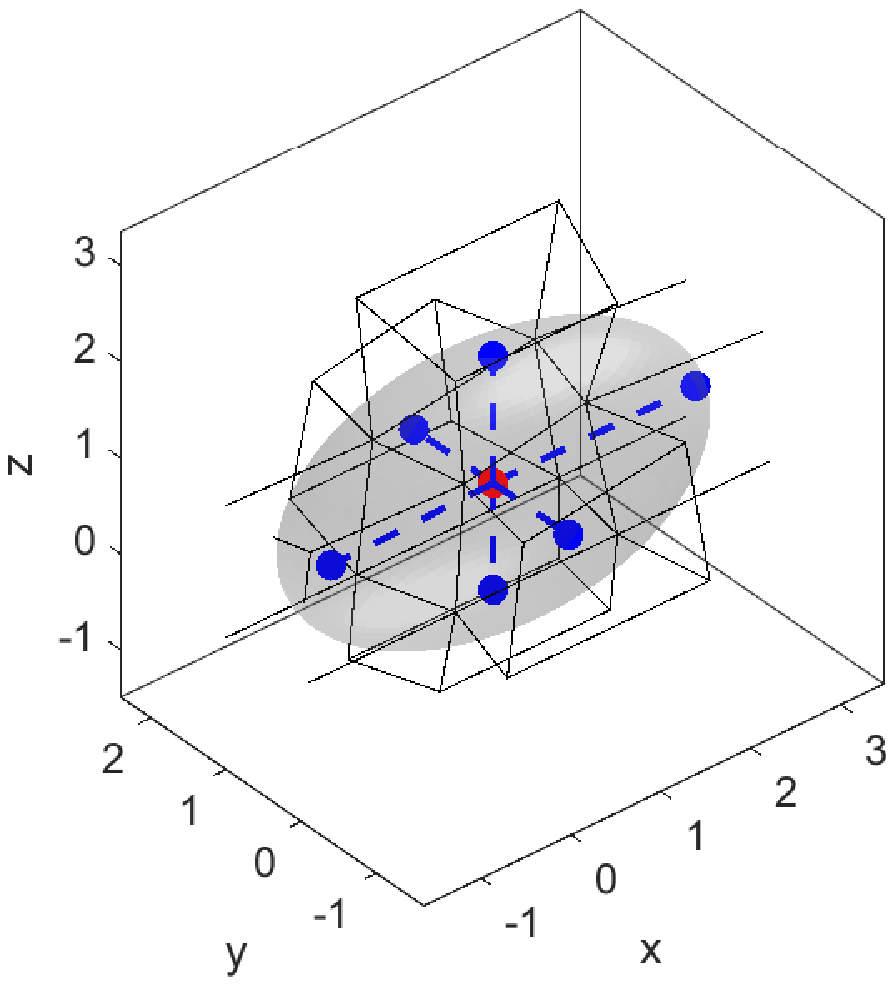}}  
    \subfigure[]{ 
    \label{fig:Support_02} 
    \includegraphics[width=0.48\textwidth]{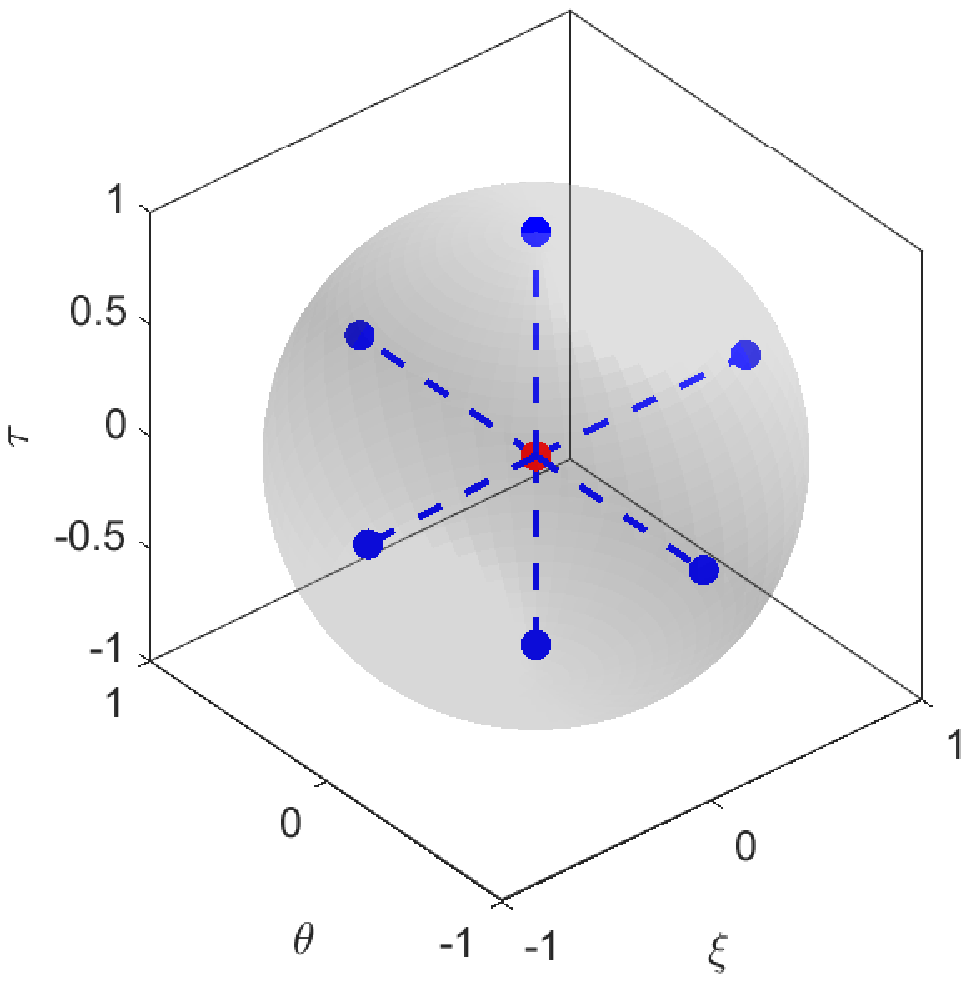}}  
  \caption{Support of a point in 3D. (a) Supporting points in the physical domain. (b) Geometric transformation into a standard domain.} 
  \label{fig:Support} 
\end{figure}

Here we introduce a local Radial Basis Function-based Differential Quadrature (RBF-DQ) method. The approach is first proposed in 2D by Shu et al. \cite{Shu2003}. In the current work, a modified linearly complete RBF-DQ method is presented for multi-dimensional and unevenly distributed points. First, for Fragile Point $P_0$, all its supporting points are transformed into a standard computational domain (see Figure~\ref{fig:Support}). That is:
\begin{align}
\begin{split}
& P_i \left( \mathbf{x}_i \right) \Rightarrow P_i \left( \mathbf{\boldsymbol{\xi}_i} \right), \quad \text{where } \boldsymbol{\xi} = \left\{ \begin{matrix} \left[ \xi, \theta \right]^\mathrm{T} \; \text{in 2D} \\ \left[ \xi, \theta, \tau \right]^\mathrm{T} \; \text{in 3D} \end{matrix} \right.\\
& \xi_i = \frac{1}{l_x} \left( x_i - x_0 \right), \quad \theta_i = \frac{1}{l_y} \left( y_i - y_0 \right), \quad \tau_i = \frac{1}{l_z} \left( z_i - z_0 \right), \\
& l_x = \mathrm{max} \left( \left| x_i - x_0 \right| \right), \quad l_y = \mathrm{max} \left( \left| y_i - y_0 \right| \right), \quad l_z = \mathrm{max} \left( \left| z_i - z_0 \right| \right).
\end{split}
\end{align}

Multiple RBFs, including the multiquadric (MQ), inverse-MQ and Gaussians can be used as basis function in the local RBF-DQ method. Here the MQ-RBF is shown as an example:
\begin{align}
\begin{split}
\psi \left( r \right) = \sqrt{r^2 + c^2},
\end{split}
\end{align}
where $r$ is the radial length from the conference point, $c$ is a constant parameter. 

The temperature field $u \left( \boldsymbol{\xi} \right) $ can be locally approximated by MQ-RBFs and additional polynomials as:
\begin{align}
\begin{split}
u \left( \boldsymbol{\xi} \right) = \sum_{i = 0}^m \lambda_i \psi \left( \left\| \boldsymbol{\xi} - \boldsymbol{\xi}_i \right\|_2 \right) + \zeta_{0} + \boldsymbol{\zeta} \cdot \boldsymbol{\xi},
\end{split} \label{eq:RBF_1}
\end{align}
where $\boldsymbol{\zeta} = \left[ \zeta_1, \zeta_2 \right]$ in 2D or $\boldsymbol{\zeta} = \left[ \zeta_1, \zeta_2,  \zeta_3 \right]$ in 3D. To make the problem well-determined, the following conditions are given in addition to the $\left( m+1\right)$ collocation equations at $P_i$:
\begin{align}
\begin{split}
\sum_{i = 0}^m \lambda_i = 0, \quad \sum_{i = 0}^m \lambda_i \boldsymbol{\xi}_i = \mathbf{0}.
\end{split}
\end{align}
Thus Eqn.~\ref{eq:RBF_1} can be rewritten as:
\begin{align}
\begin{split}
u \left( \boldsymbol{\xi} \right) = \sum_{i = 0}^m \zeta_{i} g_i \left( \boldsymbol{\xi} \right),
\end{split}
\end{align}
where
\begin{align}
\begin{split} \notag
& g_0 \left( \boldsymbol{\xi} \right) = 1; \quad \left[ \begin{matrix} g_1 \left( \boldsymbol{\xi} \right) & \cdots & g_N \left( \boldsymbol{\xi} \right) \end{matrix} \right] = \boldsymbol{\xi}^\mathrm{T} ; \\
& g_i \left( \boldsymbol{\xi} \right) = \psi_i \left( \boldsymbol{\xi} \right) - \left[ \begin{matrix} \psi_0 \left( \boldsymbol{\xi} \right) & \psi_1 \left( \boldsymbol{\xi} \right) & \cdots & \psi_N \left( \boldsymbol{\xi} \right) \end{matrix} \right]  \left[ \begin{matrix} 1 & 1 & \cdots & 1 \\ \boldsymbol{\xi}_0 & \boldsymbol{\xi}_1 & \cdots & \boldsymbol{\xi}_N \end{matrix} \right]^\mathrm{-1} \left[\begin{matrix} 1 \\  \boldsymbol{\xi}_i \end{matrix}\right], \\
& \psi_i \left( \boldsymbol{\xi} \right) = \psi \left( \left\| \boldsymbol{\xi} - \boldsymbol{\xi}_i \right\|_2 \right) ,  \quad i = N+1, \cdots, m,
\end{split}
\end{align}
$N$ is the dimension of the problem.

In the DQ method, any partial derivative of temperature field at point $P_0$ can be approximated by a weighted linear sum of the value of $u_h$ at all the supporting points, i.e.:
\begin{align}
\begin{split}
D u \left( \boldsymbol{\xi}_0 \right) = \sum_{i=0}^m W_i^{D} u \left( \boldsymbol{\xi}_i \right),
\end{split}
\end{align}
where $D$ is a linear differential operator, $W_i^{D}$ is the corresponding weighting coefficient at point $P_i$. These weighting coefficients can be determined by a set of base vectors $g_i \left( \boldsymbol{\xi} \right)$ ($i = 0, 1, \cdots, m$). Therefore, we arrive at the final approximation of the first and second derivatives at point $P_0$:
\begin{align}
\begin{split}
\mathbf{D} u \Big| _{P_0} = \mathbf{B} \mathbf{u}_E,
\end{split}
\end{align}
where
\begin{align}
\begin{split} 
& \mathbf{D} = \left\{ \begin{matrix} \left[ \begin{matrix} \frac{\partial}{\partial x} & \frac{\partial}{\partial y} & \frac{\partial^2}{\partial x^2} & \frac{\partial^2}{\partial y^2} & \frac{\partial^2}{\partial x \partial y} \end{matrix} \right]  ^ {\mathrm{T}} \; & \text{in 2D} \\ \left[ \begin{matrix} \frac{\partial}{\partial x} & \frac{\partial}{\partial y} & \frac{\partial}{\partial z} & \frac{\partial^2}{\partial x^2} & \frac{\partial^2}{\partial y^2} & \frac{\partial^2}{\partial z^2} & \frac{\partial^2}{\partial x \partial y} & \frac{\partial^2}{\partial y \partial z} & \frac{\partial^2}{\partial x \partial z} \end{matrix} \right]  ^ {\mathrm{T}} \; & \text{in 3D}\end{matrix} \right. ,\\
& \mathbf{u}_E = \left[ \begin{matrix} u_0 & u_1 & u_2 & \cdots & u_m \end{matrix} \right] ^ {\mathrm{T}}, \\
& \mathbf{B}  = \mathbf{J}^\mathrm{-1} \left( \mathbf{D}_{\boldsymbol{\xi}} \mathbf{G}_0^\mathrm{T} \right) \mathbf{G}^\mathrm{-T} , \\
& \mathbf{J} =  \left\{ \begin{matrix}  \mathrm{diag} \left( \left[ \begin{matrix} l_x & l_y & l_x^2 & l_y^2 & l_x l_y \end{matrix}  \right] \right) \; & \text{in 2D} \\ \mathrm{diag} \left( \left[ \begin{matrix} l_x & l_y & l_z & l_x^2 & l_y^2 & l_z^2 & l_x l_y & l_y l_z & l_x l_z \end{matrix}  \right] \right) \; & \text{in 3D}\end{matrix} \right. ,\\
& \mathbf{D}_{\boldsymbol{\xi}} = \left\{ \begin{matrix} \left[ \begin{matrix} \frac{\partial}{\partial \xi} & \frac{\partial}{\partial \theta} & \frac{\partial^2}{\partial \xi^2} & \frac{\partial^2}{\partial \theta^2} & \frac{\partial^2}{\partial \xi \partial \theta} \end{matrix} \right]  ^ {\mathrm{T}} \; & \text{in 2D} \\ \left[ \begin{matrix} \frac{\partial}{\partial \xi} & \frac{\partial}{\partial \theta} & \frac{\partial}{\partial \tau} & \frac{\partial^2}{\partial \xi^2} & \frac{\partial^2}{\partial \theta^2} & \frac{\partial^2}{\partial \tau^2} & \frac{\partial^2}{\partial \xi \partial \theta} & \frac{\partial^2}{\partial \theta \partial \tau} & \frac{\partial^2}{\partial \xi \partial \tau} \end{matrix} \right]  ^ {\mathrm{T}} \; & \text{in 3D}\end{matrix} \right. ,\\
& \mathbf{G} = \left[ \begin{matrix} \mathbf{G}_0 & \mathbf{G}_1 & \cdots & \mathbf{G}_m \end{matrix} \right], \quad  \mathbf{G}_i = \left[ \begin{matrix} g_0 \left( \boldsymbol{\xi}_i \right) & g_1 \left( \boldsymbol{\xi}_i \right) & \cdots & g_m \left( \boldsymbol{\xi}_i \right) \end{matrix} \right]^\mathrm{T}.
\end{split}
\end{align}

Substituting the approximation into Eqn.~\ref{eq:trial_0}, the relation between $u_h$ and $\mathbf{u}_E$ is obtained:
\begin{align}
\begin{split}
u_h (\mathbf{x}) = \mathbf{N} \left( \mathbf{x} \right) \mathbf{u}_E, \quad \mathbf{x} \in E_0,
\end{split}
\end{align}
where $\mathbf{N} \left( \mathbf{x} \right)$ is the shape function of $u_h$:
\begin{align}
\begin{split}
& \mathbf{N} \left( \mathbf{x} \right) = \overline{\mathbf{N}} \left( \mathbf{x} \right) \mathbf{B} + \left[ \begin{matrix} 1 & 0 & \cdots & 0 \end{matrix} \right] _ {1 \times (m+1)}, \\
& \overline{\mathbf{N}} \left( \mathbf{x} \right) = \left\{ \begin{matrix}  \left[ \begin{matrix} x-x_0 & y - y_0 & \frac{1}{2} \left( x - x_0 \right)^2 & \frac{1}{2} \left( y - y_0 \right)^2 & \left( x - x_0 \right) \left( y - y_0 \right) \end{matrix} \right] \; & \text{in 2D} \\ \left[ \begin{matrix} \begin{matrix} x-x_0 & y - y_0 & \frac{1}{2} \left( x - x_0 \right)^2 & \frac{1}{2} \left( y - y_0 \right)^2 & \frac{1}{2} \left( z - z_0 \right)^2 & \cdots \end{matrix} \\ \begin{matrix} \left( x - x_0 \right) \left( y - y_0 \right) & \left( y - y_0 \right) \left( z - z_0 \right) & \left( x - x_0 \right) \left( z - z_0 \right) \end{matrix} \end{matrix} \right]  \; & \text{in 3D}\end{matrix} \right. .
\end{split}
\end{align}

\begin{figure}[htbp] 
  \centering 
    \subfigure[]{ 
    \label{fig:Shape_Func_Linear} 
    \includegraphics[width=0.48\textwidth]{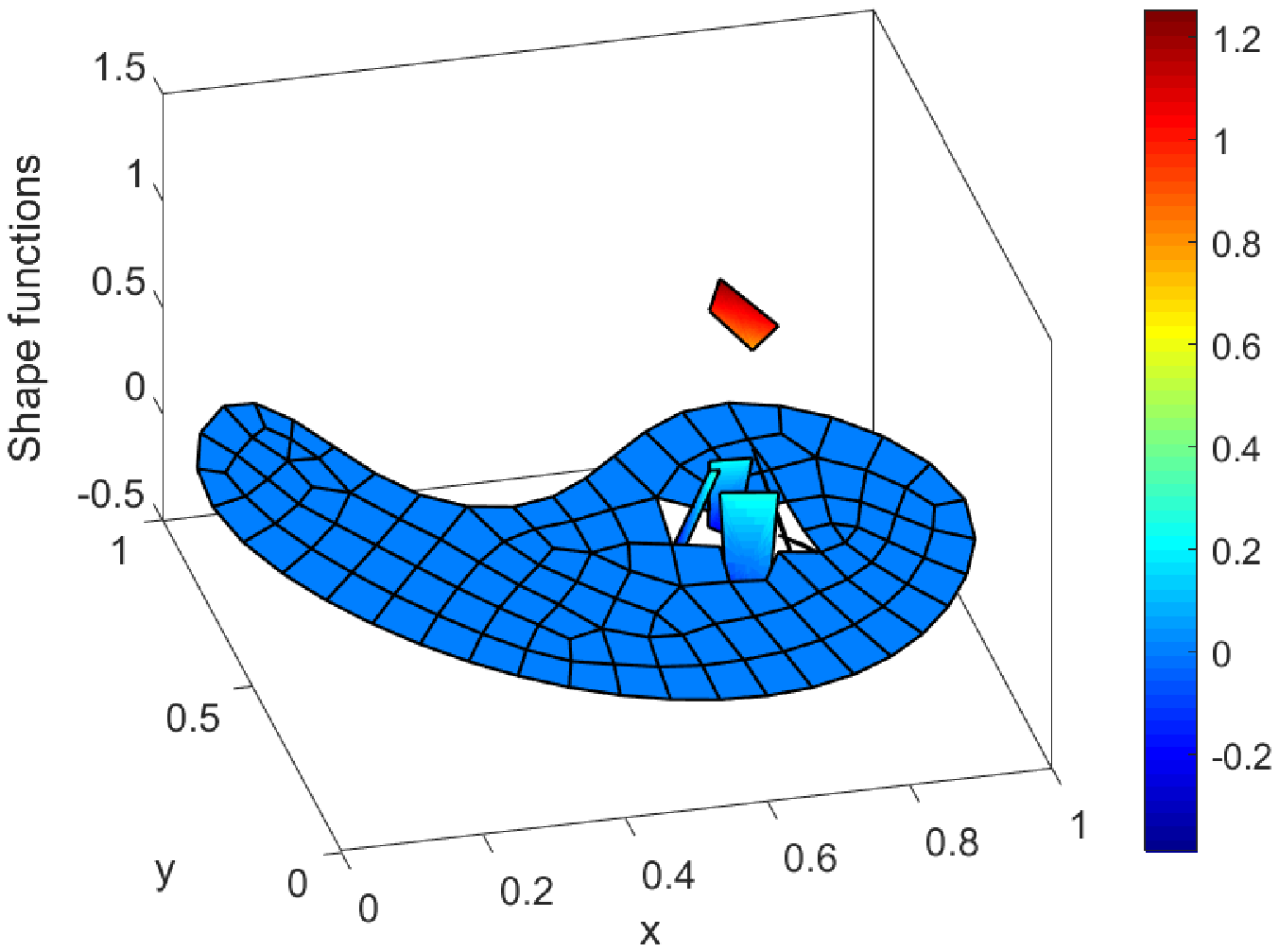}}  
    \subfigure[]{ 
    \label{fig:Shape_Func_Quad} 
    \includegraphics[width=0.48\textwidth]{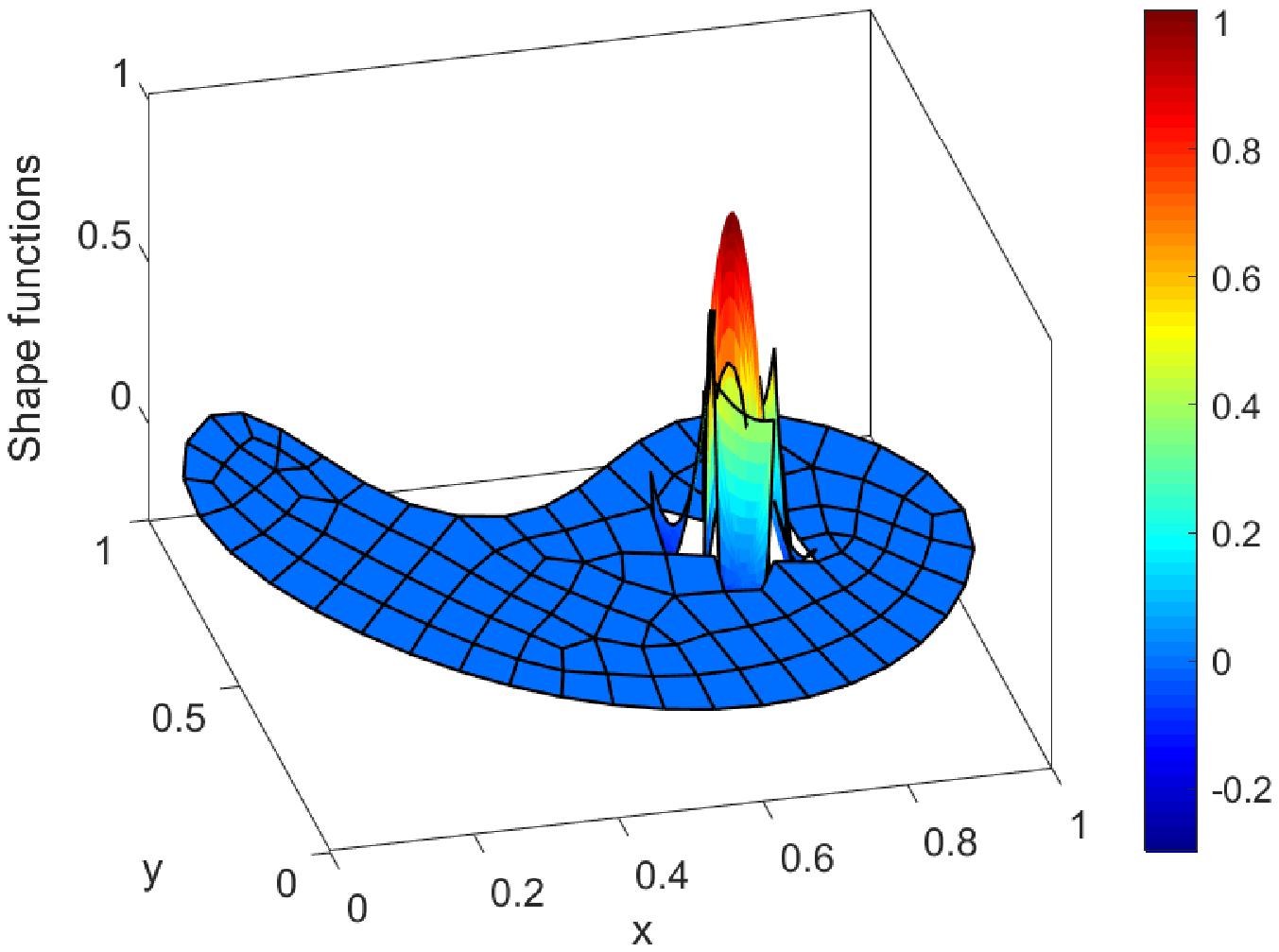}}  
  \caption{The shape functions in 2D. (a) Linear shape function. (b) Quadratic shape function.} 
  \label{fig:Shape_Func} 
\end{figure}

In some kinds of PG-FPMs, the quadratic terms in the Taylor expansion can be neglected. Thus the trial function $u_h$ is simplified to be linear. The gradient of temperature at the Internal Points $\nabla u \Big| _{P_0}$ can be approximated either by the proposed local RBF-DQ method or by a Generalized Finite Difference (GFD) method (see previous literature \cite{Guan2020, Liszka1980} for details). Both approximations utilizing in the PG-FPMs achieve similar accuracy and are not distinguished in the following discussion. Figure.~\ref{fig:Shape_Func} presents the graphs of linear and quadratic shape functions used in the PG-FPMs in a 2D domain respectively.

\section{Weak forms and test functions} \label{sec:formulation}

In this section, we present all the weak-forms of the governing equation and numerical implementations for three different kinds of PG-FPM approaches.

\subsection{The weak form 1 with Dirac delta function (collocation method)}

First, Eqn.~\ref{eqn:governing} can be written in a local weak-form with test function $v$ in subdomain $E$:
\begin{align}
\begin{split}
\int_E \rho c v \frac{\partial u}{\partial t} \mathrm{d} \Omega - \int_E v  \nabla \cdot \left( \mathbf{k} \nabla u \right) \mathrm{d} \Omega = \int_E v Q \mathrm{d} \Omega,
\label{eqn:weakform1_1}
\end{split}
\end{align}
where 
\begin{align}
\begin{split} \notag
& \nabla \cdot \left( \mathbf{k} \nabla u \right) = \overline{\mathbf{k}} \left( \mathbf{D} u \right), \\
& \overline{\mathbf{k}} = \left\{ \begin{matrix} \left[ \begin{matrix} \nabla \mathbf{k} & k_{11} & k_{22} & 2 k_{12}  \end{matrix} \right] \; & \text{in 2D} \\ \left[ \begin{matrix} \nabla \mathbf{k} & k_{11} & k_{22} & k_{33} & 2 k_{12} & 2 k_{23} & 2 k_{13} \end{matrix} \right]  \; & \text{in 3D}\end{matrix} \right. .
\end{split}
\end{align}

If the test function $v$ is taken to be the Dirac delta function, the collocation method will be derived. The gradient of thermal conductivity tensor $\nabla \mathbf{k}$ at each Fragile Point can be approximated by the value of $\mathbf{k}$ at its supporting points using the same local RBF-DQ method or GFD method presented above.

To ensure the consistency and stability of the method, the Numerical Flux Corrections are employed. Let $\Gamma$ denote the set of all internal and external boundaries, i.e., $\Gamma = \Gamma_h + \partial \Omega$, where $\Gamma_h$ is the set of all internal boundaries. Since the Delta test function vanishes on all these boundaries, another test function $\widetilde{v}$ is defined on $\Gamma$ to impose the boundary and continuity conditions. Here we define the jump operator $\left[ \! \left[ \cdot \right] \! \right]$ and average operator $\left\{ \cdot \right\}$ on $\Gamma$ as:
\begin{align}\label{eqn:operator}
\begin{split}
\left[ \! \left[ w \right] \! \right] = \begin{cases}
w \Big|_e^{E_1} -  w \Big|_e^{E_2} & e \in \Gamma_h\\
w \Big|_e & e \in \partial \Omega
\end{cases}, \quad
\left\{ w \right\} = \begin{cases}
\frac{1}{2} \left( w \Big|_e^{E_1} +  w \Big|_e^{E_2} \right) & e \in \Gamma_h\\
w \Big|_e & e \in \partial \Omega
\end{cases},
\end{split}
\end{align}
where $E_1$ and $E_2$ are the two neighboring subdomains of $e$, i.e., $e = \partial E_1 \cap \partial E_2$, for $e \in \Gamma_h$. The order of $E_1$ and $E_2$ does not affect the final formulation of the PG-FPMs. Test function $\widetilde{v} \left( \mathbf{x} \right)$ possesses the shape of Heaviside step function:
\begin{align}
\begin{split} 
& \widetilde{v} \left( \mathbf{x} \right) = \left[ \! \left[ \widetilde{v}_e \right] \! \right], \qquad \widetilde{v}_e \big|^{E_i}_e = v_i,
\end{split}
\end{align}
where $v_i$ is the value of test function $v_h$ at Internal Point $Pi \in E_i$. 

Multiplying the boundary condition equations (Eqn.~\ref{eqn:BC1} – \ref{eqn:BC4}) by $\widetilde{v}$ and integrating over the subdomain boundaries, we can obtain:
\begin{align}
& \int_e \widetilde{v} u \mathrm{d} \Gamma = \int_e \widetilde{v} \widetilde{u}_D \mathrm{d} \Gamma, & \text{for} \; e \in \Gamma_D \\
& \int_e \widetilde{v} \mathbf{n}^\mathrm{T} \mathbf{k} \nabla u \mathrm{d} \Gamma =  \int_e \widetilde{v} \widetilde{q}_N \mathrm{d} \Gamma, & \text{for} \;  e \in \Gamma_N, \\
& \int_e \widetilde{v} \mathbf{n}^\mathrm{T} \mathbf{k} \nabla u \mathrm{d} \Gamma =  \int_e \widetilde{v} h \left( \widetilde{u}_R - u \right) \mathrm{d} \Gamma, & \text{for} \;  e \in \Gamma_R,\\
& \int_e \widetilde{v} \mathbf{n}^\mathrm{T} \nabla u \mathrm{d} \Gamma = 0, & \text{for} \;  e \in \Gamma_S.
\end{align}
The continuity requirement across the internal boundaries can satisfied with Numerical Flux $\left[ \! \left[ u_h \right] \! \right]$:
\begin{align}
\begin{split} \label{eq:NF_CC}
&  \int_e \widetilde{v} \left[ \! \left[ u \right] \! \right] \mathrm{d} \Gamma = 0, \qquad  \text{for} \;  e \in \Gamma_h.
\end{split}
\end{align}
The above continuity and boundary conditions are imposed in Eqn.~\ref{eqn:weakform1_1} with penalty parameters $\eta_1$ and $\eta_2$ respectively. After summing over all subdomains, the final formula of the collocation method (PG-FPM-1) is achieved:
\begin{align}
\begin{split}
& \sum_{E \in \Omega} \int_E \rho c v \frac{\partial u}{\partial t} \mathrm{d} \Omega - \sum_{E \in \Omega} \int_E v  \overline{\mathbf{k}} \left( \mathbf{D} u \right) \mathrm{d} \Omega + \sum_{e \in \Gamma_h} \frac{1}{A_e} \int_e \frac{\eta_1}{h_e^2} \overline{k} \widetilde{v} \left[ \! \left[ u \right] \! \right] \mathrm{d} \Gamma \\
& \qquad  + \sum_{e \in \Gamma_D}  \frac{1}{A_e} \int_e \frac{\eta_2}{h_e^2} \overline{k} \widetilde{v} u \mathrm{d} \Gamma + \sum_{e \in \Gamma_N \cup \Gamma_R} \frac{1}{A_e} \int_e \frac{\eta_2}{h_e} \widetilde{v} \mathbf{n}^\mathrm{T} \mathbf{k} \nabla u \mathrm{d} \Gamma + \sum_{e \in \Gamma_R} \frac{1}{A_e} \int_e \frac{\eta_2}{h_e} h \widetilde{v} u \mathrm{d} \Gamma \\
& \qquad + \sum_{e \in \Gamma_S}  \frac{1}{A_e} \int_e \frac{\eta_2}{h_e} \overline{k} \widetilde{v} \mathbf{n}^\mathrm{T} \nabla u \mathrm{d} \Gamma = \sum_{E \in \Omega} \int_E v Q \mathrm{d} \Omega + \sum_{e \in \Gamma_D}  \frac{1}{A_e} \int_e \frac{\eta_2}{h_e^2} \overline{k} \widetilde{v} \widetilde{u}_D \mathrm{d} \Gamma \\
& \qquad + \sum_{e \in \Gamma_N} \frac{1}{A_e} \int_e \frac{\eta_2}{h_e} \widetilde{v} \widetilde{q}_N \mathrm{d} \Gamma + \sum_{e \in \Gamma_R} \frac{1}{A_e} \int_e \frac{\eta_2}{h_e} h \widetilde{v} \widetilde{u}_R \mathrm{d} \Gamma,
\label{eqn:weakform1_2}
\end{split}
\end{align}
where $A_e = \int_e \mathrm{d} \Gamma$ is the length (in 2D) or area (in 3D) of the subdomain boundary, $\overline{k}$ is an estimated average thermal conductivity in the domain. In transient analysis, $\overline{k}$ should also be larger than $\rho c \frac{h_e^2}{\Delta t}$, where $\Delta t$ is the smallest analyzing time step. $h_e$ is a boundary-dependent parameter with the unit of length. In the present work, $h_e$ is defined as the distance between the Internal Points in subdomains sharing the boundary for $e \in \Gamma_h$, and the smallest distance between the centroid of the subdomain and the external boundary for $e \in \partial \Omega$. The penalty parameters $\eta_1$ and $\eta_2$ are nondimensional positive numbers independent of the boundary size. The method is only stable when the penalty parameters are large enough. However, excessively large penalty parameters can be harmful for the accuracy. The recommended values of $\eta_1$ and $\eta_2$ will be discussed in Section~\ref{sec:PS}.

It should be pointed out that unlike the Discontinuous Galerkin (DG) Methods \cite{Arnold2001, Mozolevski2007} in which discontinuous trial functions and Numerical Flux Corrections are also employed, the trial functions defined in the FPM and PG-FPMs have an inherent ``weak continuity''. Therefore, the required penalty parameter $\eta_1$ in the present methods can be much smaller, which assures a better accuracy as compared to the DG methods. The Dirichlet boundary conditions enforced by Interior Penalty (IP) Numerical Flux terms with penalty parameter $\eta_2$ in Eqn.~\ref{eqn:weakform1_2} can also be imposed in an alternative way. When boundary Fragile Points are adopted, i.e., for $\partial E_i \cap \partial \Omega \not= \varnothing$, $P_i \in \left( \partial E_i \cap \partial \Omega \right)$, we can enforce $u = \widetilde{u}_D$ strongly at these Points and thus the corresponding IP terms in Eqn.~\ref{eqn:weakform1_2} vanish. The present IP Numerical Flux  Correction approach, on the other hand, is suitable for any point distributions, regardless of whether there are Points distributed on the external boundaries.

Substituting the test functions $v$ and $\widetilde{v}$ into Eq.~\ref{eqn:weakform1_2}, the formula can also be written in a collocation form:
\begin{align}
\begin{split}
& \rho_i c_i \frac{\partial u_i}{\partial t} - \overline{\mathbf{k}}_i \left( \mathbf{D} u \right)_i + \sum_{e \in \partial E_i \cap \Gamma_h} \frac{1}{A_e} \frac{\eta_1}{h_e^2} \int_e \overline{k} \left[ \! \left[ u \right] \! \right] \mathrm{d} \Gamma + \sum_{e \in \partial E_i \cap \Gamma_D}  \frac{1}{A_e} \frac{\eta_2}{h_e^2} \int_e \overline{k} u \mathrm{d} \Gamma \\
& \qquad + \sum_{e \in \partial E_i \cap \left( \Gamma_N \cup \Gamma_R \right)} \frac{1}{A_e} \frac{\eta_2}{h_e}  \int_e \mathbf{n}^\mathrm{T} \mathbf{k} \nabla u \mathrm{d} \Gamma + \sum_{e \in \partial E_i \cap \Gamma_R} \frac{1}{A_e} \frac{\eta_2}{h_e} \int_e  h u \mathrm{d} \Gamma \\
& \qquad + \sum_{e \in \partial E_i \cap \Gamma_S}  \frac{1}{A_e} \frac{\eta_2}{h_e}  \int_e \overline{k} \mathbf{n}^\mathrm{T} \nabla u \mathrm{d} \Gamma = Q_i + \sum_{e \in \partial E_i \cap \Gamma_D}  \frac{1}{A_e} \frac{\eta_2}{h_e^2}  \int_e \overline{k} \widetilde{u}_D \mathrm{d} \Gamma \\
& \qquad + \sum_{e \in \partial E_i \cap \Gamma_N} \frac{1}{A_e} \frac{\eta_2}{h_e}  \int_e  \widetilde{q}_N \mathrm{d} \Gamma + \sum_{e \in \partial E_i \cap \Gamma_R} \frac{1}{A_e} \frac{\eta_2}{h_e}  \int_e  h \widetilde{u}_R \mathrm{d} \Gamma.
\label{eqn:weakform1_3}
\end{split}
\end{align}
Or in a matrix form:
\begin{align}
\begin{split}
\mathbf{C} \dot{\mathbf{u}} + \mathbf{K} \mathbf{u} = \mathbf{q}, \quad \text{or} \quad  \dot{\mathbf{u}} = - \mathbf{C} ^\mathrm{T} \mathbf{K} \mathbf{u} + \mathbf{C} ^\mathrm{T} \mathbf{q}
\label{eqn:weakform1_Matrix}
\end{split}
\end{align}
where $\mathbf{C}$ and $\mathbf{K}$ are the global heat capacity and thermal conductivity matrices, $\mathbf{q}$ is the heat flux vector, and $\mathbf{u}$ is the unknown vector with nodal temperatures.

Note that in the first kind of PG-FPM, or the collocation PG-FPM, the thermal conductivity $\mathbf{K}$ is sparse and the global heat capacity matrix $\mathbf{C}$ is always diagonal. Therefore, the Jacobian matrix $\mathbf{J} = - \mathbf{C}^\mathrm{-1} \mathbf{K}$ for the ordinary differential equations (ODEs) in the time domain (Eq.~\ref{eqn:weakform1_Matrix}) is sparse and can be obtained easily. When comparing with the original FPM based on a symmetric Galerkin weak-form \cite{Guan2020}, the collocation method (PG-FPM-1) simplifies the assembling process and is more efficient in solving transient problems when boundary points ($P \in \partial \Omega$) are included, i.e., when the heat capacity matrix $\mathbf{C}$ in the original Galerkin FPM is no longer diagonal. However, the thermal conductivity matrix $\mathbf{K}$ in the current PG-FPM-1 is asymmetric. As the second derivatives of the trial function is included in Eq.~\ref{eqn:weakform1_2} – \ref{eqn:weakform1_3}, the quadratic terms in $u_h$ cannot be neglected, thus the PG-FPM-1 suffers from a complicated local approximation process. A mixed formulation using $\nabla u$ as independent variables and linear trial functions may help to remedy the problem. The mixed formulation for the FPM can be seen in the previous study \cite{LeitingDongTianYangKaileiWang2019} and is not further discussed in this paper.

\subsection{The weak form 2 with Heaviside step function (finite volume method)}

Alternatively, by integrating Eqn.~\ref{eqn:weakform1_1} by parts once, we obtain:
\begin{align}
\begin{split}
\int_E \rho c v \frac{\partial u}{\partial t} \mathrm{d} \Omega + \int_E \nabla v ^{\mathrm{T}} \mathbf{k} \nabla u \mathrm{d} \Omega = \int_E v Q \mathrm{d} \Omega + \int_{\partial E} v \mathbf{n}^{\mathrm{T}} \mathbf{k} \nabla u \mathrm{d} \Gamma.
\label{eqn:weakform2_1}
\end{split}
\end{align}
Summing the above equation over all subdomains:
\begin{align}
\begin{split}
& \sum_{E\in \Omega} \int_E \rho c v \frac{\partial u}{\partial t} \mathrm{d} \Omega + \sum_{E\in \Omega} \int_E \nabla v ^{\mathrm{T}} \mathbf{k} \nabla u \mathrm{d} \Omega = \int_\Omega v Q \mathrm{d} \Omega \\
& \qquad  + \sum_{e \in \Gamma_h} \left( \int_e \left\{ v \right\} \left[ \! \left[ \mathbf{n}^\mathrm{eT} \mathbf{k} \nabla u \right] \! \right] + \left[ \! \left[ v \right] \! \right] \left\{ \mathbf{n}^\mathrm{eT} \mathbf{k} \nabla u \right\} \right) \mathrm{d} \Gamma + \sum_{e \in \partial \Omega} \int_e  v  \mathbf{n}^\mathrm{eT} \mathbf{k} \nabla u \mathrm{d} \Gamma,
\label{eqn:weakform2_2}
\end{split}
\end{align}
where $\mathbf{n}^e$ is an unit vector normal to $e$. For $e \in \Gamma_h$, $\mathbf{n}^e = \mathbf{n}_1 = - \mathbf{n}_2$, with $\mathbf{n}_j$ being the unit normal vector pointing outward from $E_j$. The order of $E_1$ and $E_2$ is consistent with the definition of the jump operator $\left[ \!\left[ \cdot \right] \! \right]$. And for $e \in \partial \Omega$, $\mathbf{n}^e = \mathbf{n}_1$, where $E_1$ is the only one neighboring subdomain of $e$.

When $u$ is the exact solution, the continuity conditions lead to $\left[ \! \left[ \mathbf{n}^\mathrm{eT} \mathbf{k} \nabla u \right] \! \right] =0$ and $\left[ \! \left[ u \right] \! \right] =0$ on $e \in \Gamma_h$. Hence, we can replace the term $ \left\{ v \right\} \left[ \! \left[ \mathbf{n}^\mathrm{eT} \mathbf{k} \nabla u \right] \! \right]$ in Eq.~\ref{eqn:weakform2_2} by $\left\{ \nabla v^\mathrm{T} \mathbf{k} \mathbf{n}^\mathrm{e} \right\} \left[ \! \left[ u \right] \! \right]$ without influencing the consistency of the formula. The following relations are obtained from the boundary conditions:
\begin{align}
\begin{split}
\sum_{e \in \Gamma_D} \int_e  v  \mathbf{n}^\mathrm{eT} \mathbf{k} \nabla u \mathrm{d} \Gamma = & \sum_{e \in \Gamma_D} \int_e \left( \left[ \! \left[ v \right] \! \right] \left\{ \mathbf{n}^\mathrm{eT} \mathbf{k} \nabla u \right\} + \left\{ \nabla v^\mathrm{T} \mathbf{k} \mathbf{n}^\mathrm{e} \right\} \left[ \! \left[ u \right] \! \right]  \right) \mathrm{d} \Gamma  \\
& - \sum_{e \in \Gamma_D} \int_e \left( \nabla v^\mathrm{T} \mathbf{k} \mathbf{n}^\mathrm{e} \right) \widetilde{u}_D \mathrm{d} \Gamma, \\
\sum_{e \in \Gamma_N} \int_e  v  \mathbf{n}^\mathrm{eT} \mathbf{k} \nabla u \mathrm{d} \Gamma = & \sum_{e \in \Gamma_N} \int_e v  \widetilde{q}_N  \mathrm{d} \Gamma, \\
\sum_{e \in \Gamma_R} \int_e  v  \mathbf{n}^\mathrm{eT} \mathbf{k} \nabla u \mathrm{d} \Gamma = & \sum_{e \in \Gamma_R} \int_e h v  \left( \widetilde{u}_R - u \right) \mathrm{d} \Gamma, \\
\sum_{e \in \Gamma_S} \int_e  v  \mathbf{n}^\mathrm{eT} \mathbf{k} \nabla u \mathrm{d} \Gamma = & \sum_{e \in \Gamma_S} \int_e v \mathbf{n}^\mathrm{e T} \mathbf{k} \left( \mathbf{I} -  \mathbf{n}^\mathrm{e} \mathbf{n}^\mathrm{e T} \right) \nabla u \mathrm{d} \Gamma,
\end{split}
\end{align}
where $\mathbf{I}$ is the unit matrix. On $\Gamma_S$, we have $\nabla u = \left( \mathbf{I} -  \mathbf{n}^\mathrm{e} \mathbf{n}^\mathrm{e T} \right) \nabla u  + \mathbf{n}^\mathrm{e} \left( \mathbf{n}^\mathrm{e T} \nabla u \right)  = \left( \mathbf{I} -  \mathbf{n}^\mathrm{e} \mathbf{n}^\mathrm{e T} \right) \nabla u$. When the medium is isotropic, the term on $\Gamma_S$ vanishes, i.e., the symmetric boundary $\Gamma_S$ is equivalent to $\Gamma_N$ with $\widetilde{q}_N = 0$.

Finally, two IP Numerical Flux terms $\left[ \! \left[ u_h \right] \! \right]$ and $\widetilde{u}_D$ are applied to $\Gamma_h$ and $\Gamma_D$ respectively. The formula for the second kind of PG-FPM is obtained:
\begin{align}
\begin{split}
& \sum_{E \in \Omega} \int_E \rho c v \frac{\partial u}{\partial t} \mathrm{d} \Omega + \sum_{E \in \Omega} \int_E \nabla v^\mathrm{T} \mathbf{k} \nabla u \mathrm{d} \Omega - \sum_{e \in \Gamma_S} \int_e v \mathbf{n}^\mathrm{e T} \mathbf{k} \left( \mathbf{I} -  \mathbf{n}^\mathrm{e} \mathbf{n}^\mathrm{e T} \right) \nabla u \mathrm{d} \Gamma \\
& \qquad   + \sum_{e \in \Gamma_R} \int_e h v u \mathrm{d} \Gamma - \sum_{e \in \Gamma_h \cup \Gamma_D} \int_e \left( \left[ \! \left[ v \right] \! \right] \left\{ \mathbf{n}^\mathrm{eT} \mathbf{k} \nabla u \right\} + \left\{ \nabla v^\mathrm{T} \mathbf{k} \mathbf{n}^\mathrm{e} \right\} \left[ \! \left[ u \right] \! \right]  \right) \mathrm{d} \Gamma   \\
& \qquad + \sum_{e \in \Gamma_h} \frac{\eta_1}{h_e} \int_e \overline{k} \left[ \! \left[ v \right] \! \right] \left[ \! \left[ u \right] \! \right] \mathrm{d} \Gamma + \sum_{e \in \Gamma_D} \frac{\eta_2}{h_e} \int_e \overline{k} v u  \mathrm{d} \Gamma = \sum_{E \in \Omega} \int_E v Q \mathrm{d} \Omega \\
& \qquad + \sum_{e \in \Gamma_D} \int_e \left( \frac{\eta_2}{h_e} \overline{k} v - \nabla v^\mathrm{T} \mathbf{k} \mathbf{n}^\mathrm{e} \right) \widetilde{u}_D \mathrm{d} \Gamma  + \sum_{e \in \Gamma_N} \int_e v  \widetilde{q}_N  \mathrm{d} \Gamma + \sum_{e \in \Gamma_R} \int_e h v  \widetilde{u}_R  \mathrm{d} \Gamma.
\label{eqn:weakform2_4}
\end{split}
\end{align}
This is a symmetric weak form. If the test function $v$ is selected to have the same polynomial shape as the trial function $u_h$, the Galerkin FPM is achieved \cite{Guan2020}. For isotropic problems, the Galerkin FPM leads to symmetric and sparse matrices. Since there exist only the first derivatives of $u$ and $v$, linear trial and test functions are sufficient, and thus the local approximation calculations are heavily reduced.

In the present work, alternatively, we choose the Heaviside step function as the test function, i.e. $v \left( \mathbf{x} \right) = v_i$ for $\mathbf{x} \in E_i$, and obtain:
\begin{align}
\begin{split}
& \int_{E} \rho c \frac{\partial u}{\partial t} \mathrm{d} \Omega - \sum_{e \in \partial {E} \cap \Gamma_S} \int_e \mathbf{n}^\mathrm{e T} \mathbf{k} \left( \mathbf{I} -  \mathbf{n}^\mathrm{e} \mathbf{n}^\mathrm{e T} \right) \nabla u \mathrm{d} \Gamma + \sum_{e \in \partial {E} \cap  \Gamma_R} \int_e h u \mathrm{d} \Gamma \\
& \qquad    - \sum_{e \in \partial {E} \cap  \left( \Gamma_h \cup \Gamma_D \right) } \int_e  \left\{ \mathbf{n}^\mathrm{eT} \mathbf{k} \nabla u \right\} \mathrm{d} \Gamma  + \sum_{e \in \partial {E} \cap \Gamma_h} \frac{\eta_1}{h_e} \int_e \overline{k} \left[ \! \left[ u \right] \! \right] \mathrm{d} \Gamma + \sum_{e \in \partial {E} \cap \Gamma_D} \frac{\eta_2}{h_e} \int_e \overline{k} u  \mathrm{d} \Gamma \\
& \qquad  = \int_{E} Q \mathrm{d} \Omega + \sum_{e \in \partial {E} \cap \Gamma_D} \frac{\eta_2}{h_e} \int_e \overline{k} \widetilde{u}_D \mathrm{d} \Gamma  + \sum_{e   \in \partial {E} \cap \Gamma_N} \int_e  \widetilde{q}_N  \mathrm{d} \Gamma + \sum_{e \in  \partial {E} \cap  \Gamma_R} \int_e h \widetilde{u}_R  \mathrm{d} \Gamma ,
\label{eqn:weakform2_5}
\end{split}
\end{align}
where $E$ is employed as the first neighboring subdomain ($E_1$) in the jump operator $\left[ \! \left[ \cdot \right] \! \right]$, and $\mathbf{n}^\mathrm{e}$ is pointing outward from $E$.

Eqn.~\ref{eqn:weakform2_5} indicates the integration of the governing equation~\ref{eqn:governing} over subdomain $E$. Therefore, the second kind of PG-FPM is also known as a finite volume method (PG-FPM-2). Yet unlike the conventional FVMs, in which the gradient of temperature $\nabla u$ is interpolated on each subdomain boundary independently, in the present PG-FPM, $\nabla u$ is approximated on the Internal Points and is assumed constant in the entire subdomain. Consequently, the local discretization is consistent for both isotropic and anisotropic materials, whereas the two-point discretization scheme widely used in conventional FVMs is only consistent for isotropic problems \cite{Liu2015, Prestini2017}. The conservation laws are weakly satisfied by the IP Numerical Flux Corrections in the current PG-FPM-2. With appropriate penalty parameters, the PG-FPM-2 shows good consistency and stability in both isotropic and anisotropic media.

The heat capacity and thermal conductivity matrices ($\mathbf{C}$ and $\mathbf{K}$) in the finite volume method (PG-FPM-2) are still sparse but not symmetric. However, comparing with the Galerkin FPM, the assembling process of the present PG-FPM-2 is significantly simplified. With the test function $v_h$ defined in each subdomain and completely independent of the others, the bandwidths of matrices $\mathbf{C}$ and $\mathbf{K}$ in the current PG-FPM-2 are minimized. Therefore, the finite volume method (PG-FPM-2) has remarkably higher efficiency as compared to the Galerkin FPM with the same number of Fragile Points.

As a result of the linear trial functions and step test functions, the weak-form integration is extremely simple in the present PG-FPM-2 approach. Guassian quadrature scheme with only one integration point in each subdomain is sufficient. If the internal Fragile Points are placed at the centroid of the subdomains, they coincide with the integration points and thus lead to a diagonal heat capacity matrix. For transient problems, the same as the collocation method (PG-FPM-1), the corresponding Jacobian matrix $\mathbf{J}$ is sparse with a relatively small bandwidth. Therefore, the efficiency of the current finite volume method (PG-FPM-2) is the highest among the Galerkin FPM, the PG-FPMs and all the previous meshless methods based on MLS or other complicated local approximations. Note that when boundary Fragile Points are employed, the Points are no longer coincided with the integration points. However, these degrees-of-freedom (DoFs) can be eliminated by imposing Dirichlet boundary conditions directly and thus the final heat capacity matrix remains diagonal.

\subsection{The weak form 3 with a fundamental solution (singular solution method)}

At last, integrating Eq.~\ref{eqn:weakform1_1} by parts twice yields the following asymmetric local weak form:
\begin{align}
\begin{split}
\int_E \rho c v \frac{\partial u}{\partial t} \mathrm{d} \Omega - \int_E \nabla \cdot \left( \mathbf{k} \nabla v \right) u \mathrm{d} \Omega = \int_E v Q \mathrm{d} \Omega + \int_{\partial E} v \mathbf{n}^{\mathrm{T}} \mathbf{k} \nabla u \mathrm{d} \Gamma - \int_{\partial E} \nabla v^\mathrm{T} \mathbf{k} \mathbf{n} u \mathrm{d} \Gamma.
\label{eqn:weakform3_1}
\end{split}
\end{align}
Summing it over all subdomains, when $u$ is the exact solution, one obtains:
\begin{align}
\begin{split}
& \sum_{E \in \Omega} \int_E \rho c v \frac{\partial u}{\partial t} \mathrm{d} \Omega \\
& \qquad - \sum_{E \in \Omega} \int_E \nabla \cdot \left( \mathbf{k} \nabla v \right) u \mathrm{d} \Omega - \sum_{e \in \Gamma_h} \int_e \left( \left[ \! \left[ v \right] \! \right] \left\{ \mathbf{n}^\mathrm{eT} \mathbf{k} \nabla u \right\} + \left\{ \nabla v^\mathrm{T} \mathbf{k} \mathbf{n}^\mathrm{e} \right\} \left[ \! \left[ u \right] \! \right]  \right)  \mathrm{d} \Gamma \\
& \qquad + \sum_{e \in \Gamma_h} \int_e \left( \left\{ v \right\} \left[ \! \left[ \mathbf{n}^\mathrm{eT} \mathbf{k} \nabla u \right] \! \right] + \left[ \! \left[ \nabla v^\mathrm{T} \mathbf{k} \mathbf{n}^\mathrm{e} \right] \! \right]  \left\{ u \right\}  \right)  \mathrm{d} \Gamma \\
& \qquad = \sum_{E \in \Omega} \int_E v Q \mathrm{d} \Omega + \sum_{e \in \partial \Omega} \int_{e} \left( v \mathbf{n}^{\mathrm{e T}} \mathbf{k} \nabla u  - \nabla v^\mathrm{T} \mathbf{k} \mathbf{n}^\mathrm{e} u \right) \mathrm{d} \Gamma.
\label{eqn:weakform3_2}
\end{split}
\end{align}
Substituting all the boundary conditions into Eqn.~\ref{eqn:weakform3_2}. Additional Numerical Flux terms are applied on all the external boundaries ($\partial \Omega = \Gamma_D \cup \Gamma_N \cup \Gamma_R \cup \Gamma_S$) and internal boundaries $\Gamma_h$ to enforce the boundary conditions and continuity requirement. Hence, the formula of the third kind of weak-form is given as:
\begin{align}
\begin{split}
& \sum_{E \in \Omega} \int_E \rho c v \frac{\partial u}{\partial t} \mathrm{d} \Omega - \sum_{E \in \Omega} \int_E \nabla \cdot \left( \mathbf{k} \nabla v \right) u \mathrm{d} \Omega + \sum_{e \in \Gamma_h} \frac{\eta_1}{h_e} \int_e \overline{k} \left[ \! \left[ v \right] \! \right] \left[ \! \left[ u \right] \! \right] \mathrm{d} \Gamma + \sum_{e \in \Gamma_D} \frac{\eta_2}{h_e} \int_e \overline{k} v u  \mathrm{d} \Gamma \\
& \qquad - \sum_{e \in \Gamma_h \cup \Gamma_D} \int_e \left( \left[ \! \left[ v \right] \! \right] \left\{ \mathbf{n}^\mathrm{eT} \mathbf{k} \nabla u \right\} + \left\{ \nabla v^\mathrm{T} \mathbf{k} \mathbf{n}^\mathrm{e} \right\} \left[ \! \left[ u \right] \! \right]  \right)  \mathrm{d} \Gamma \\
& \qquad + \sum_{e \in \Gamma_h \cup \Gamma_N \cup \Gamma_R \cup \Gamma_S} \int_e \left( \left\{ v \right\} \left[ \! \left[ \mathbf{n}^\mathrm{eT} \mathbf{k} \nabla u \right] \! \right] + \left[ \! \left[ \nabla v^\mathrm{T} \mathbf{k} \mathbf{n}^\mathrm{e} \right] \! \right]  \left\{ u \right\}  \right)  \mathrm{d} \Gamma \\
& \qquad + \sum_{e \in \Gamma_R} \int_e h \left( 2 v + \frac{\eta_2 h_e}{\overline{k}} \nabla v^\mathrm{T} \mathbf{k} \mathbf{n}^\mathrm{e}  \right) u \mathrm{d} \Gamma + \sum_{e \in \Gamma_N \cup \Gamma_R} \frac{\eta_2 h_e}{\overline{k}} \int_e \nabla v^\mathrm{T} \mathbf{k} \mathbf{n}^\mathrm{e} \mathbf{n}^\mathrm{e T} \mathbf{k} \nabla u \mathrm{d} \Gamma\\
& \qquad - \sum_{e \in \Gamma_S} 2 \int_e v \mathbf{n}^\mathrm{e T} \mathbf{k} \nabla u \mathrm{d} \Gamma +  \sum_{e \in \Gamma_S} \int_e \left( 2 v +  \frac{\eta_2 h_e}{\overline{k}} \nabla v^\mathrm{T} \mathbf{k} \mathbf{n}^\mathrm{e} \right) \mathbf{n}^\mathrm{e T} \mathbf{k} \mathbf{n}^\mathrm{e} \mathbf{n}^\mathrm{e T} \nabla u \mathrm{d} \Gamma \\
& \qquad = \sum_{E \in \Omega} \int_E v Q \mathrm{d} \Omega + \sum_{e \in \Gamma_D} \int_e \left( \frac{\eta_2}{h_e}  \overline{k} v - 2 \nabla v^\mathrm{T} \mathbf{k} \mathbf{n}^\mathrm{e}  \right) \widetilde{u}_D  \mathrm{d} \Gamma   \\
& \qquad + \sum_{e \in \Gamma_N} \int_e \left( 2 v + \frac{\eta_2 h_e}{\overline{k}} \nabla v^\mathrm{T} \mathbf{k} \mathbf{n}^\mathrm{e} \right)  \widetilde{q}_N  \mathrm{d} \Gamma + \sum_{e \in \Gamma_R} \int_e h \left( 2 v + \frac{\eta_2 h_e}{\overline{k}} \nabla v^\mathrm{T} \mathbf{k} \mathbf{n}^\mathrm{e}  \right) \widetilde{u}_R  \mathrm{d} \Gamma.
\label{eqn:weakform3_4}
\end{split}
\end{align}

In order to simplify the above equation, the test function $v$ is selected to be the fundamental solution of the steady-state heat equation, i.e., $\nabla \cdot \left( \mathbf{k} \nabla v \right) = 0$. We assume the material to be homogenous and orthotropic in each subdomain, i.e., $\mathbf{k} \left( \mathbf{x} \right)  = \mathbf{k} = \text{diag} \left( \left[ k_{11}, k_{22} \right] \right)$ or $\text{diag} \left( \left[ k_{11}, k_{22}, k_{33} \right] \right)$. The corresponding fundamental solutions and test function $v$ can be written as:
\begin{align}
\begin{split}
& \psi \left( \mathbf{x} \right) = \left\{ \begin{matrix} \text{ln} \left[ \left( \frac{1}{k_{11}} \left( x - \overline{x}_0 \right)^2 + \frac{1}{k_{22}} \left( y - \overline{y}_0 \right)^2 \right)^{1/2} \right] \; & \text{in 2D} \\ \left[ \frac{1}{k_{11}} \left( x - \overline{x}_0 \right)^2 + \frac{1}{k_{22}} \left( y - \overline{y}_0 \right)^2 + \frac{1}{k_{33}} \left( z - \overline{z}_0 \right)^2 \right]^{-1/2}  \; & \text{in 3D}\end{matrix} \right. , \\
& v \left( \mathbf{x} \right) = v_i \overline{\psi} \left( \mathbf{x} \right) = v_i \psi \left( \mathbf{x} \right) / \psi \left( \mathbf{x}_i \right), \qquad \text{for} \; \mathbf{x} \in E_i.
\end{split}
\end{align}
where $\overline{P}_0 \left( \overline{x}_0, \overline{y}_0 \right)$ or $\overline{P}_0 \left( \overline{x}_0, \overline{y}_0, \overline{z}_0 \right)$ is a fixed point outside the studied domain $\Omega$, $P_i \left( \mathbf{x}_i \right)$ is the internal Point in subdomain $E_i$, and $v_i$ is the corresponding value of $v$ at $P_i$. The test function is discontinuous and defined independently in each subdomain.

The second term in Eqn.~\ref{eqn:weakform3_4} vanishes and the equation can be rewritten at a point level:
\begin{align}
\begin{split}
& \int_E \rho c \overline{\psi} \frac{\partial u}{\partial t} \mathrm{d} \Omega + \sum_{e \in \partial E \cap \Gamma_h} \frac{\eta_1}{h_e} \int_e \overline{k} \overline{\psi} \left[ \! \left[ u \right] \! \right] \mathrm{d} \Gamma + \sum_{e \in \partial E \cap \Gamma_D} \frac{\eta_2}{h_e} \int_e \overline{k} \overline{\psi} u  \mathrm{d} \Gamma \\
& \qquad - \sum_{e \in \partial E \cap \left( \Gamma_h \cup \Gamma_D \right)} \int_e \left( \overline{\psi} \left\{ \mathbf{n}^\mathrm{eT} \mathbf{k} \nabla u \right\} + \frac{1}{2} \nabla \overline{\psi}^\mathrm{T} \mathbf{k} \mathbf{n}^\mathrm{e} \left[ \! \left[ u \right] \! \right]  \right)  \mathrm{d} \Gamma \\
& \qquad + \sum_{e \in \partial E \cap \left( \Gamma_h \cup \Gamma_N \cup \Gamma_R \cup \Gamma_S \right)} \int_e \left( \frac{1}{2} \overline{\psi} \left[ \! \left[ \mathbf{n}^\mathrm{eT} \mathbf{k} \nabla u \right] \! \right] + \nabla \overline{\psi}^\mathrm{T} \mathbf{k} \mathbf{n}^\mathrm{e} \left\{ u \right\}  \right)  \mathrm{d} \Gamma \\
& \qquad + \sum_{e \in \partial E \cap \Gamma_R} \int_e h \left( 2 \overline{\psi} + \frac{\eta_2 h_e}{\overline{k}} \nabla \overline{\psi}^\mathrm{T} \mathbf{k} \mathbf{n}^\mathrm{e}  \right) u \mathrm{d} \Gamma - \sum_{e \in \partial E \cap \Gamma_S} 2 \int_e \overline{\psi} \mathbf{n}^\mathrm{e T} \mathbf{k} \nabla u \mathrm{d} \Gamma  \\
& \qquad +  \sum_{e \in \partial E \cap \Gamma_S} \int_e \left( 2 \overline{\psi} +  \frac{\eta_2 h_e}{\overline{k}} \nabla \overline{\psi}^\mathrm{T} \mathbf{k} \mathbf{n}^\mathrm{e} \right) \mathbf{n}^\mathrm{e T} \mathbf{k} \mathbf{n}^\mathrm{e} \mathbf{n}^\mathrm{e T} \nabla u \mathrm{d} \Gamma \\
& \qquad + \sum_{e \in \partial E \cap \left( \Gamma_N \cup \Gamma_R \right)} \frac{\eta_2 h_e}{\overline{k}} \int_e \nabla \overline{\psi}^\mathrm{T} \mathbf{k} \mathbf{n}^\mathrm{e} \mathbf{n}^\mathrm{e T} \mathbf{k} \nabla u \mathrm{d} \Gamma = \int_E \overline{\psi} Q \mathrm{d} \Omega \\
& \qquad + \sum_{e \in \partial E \cap \Gamma_D} \int_e \left( \frac{\eta_2}{h_e}  \overline{k} \overline{\psi} - 2 \nabla \overline{\psi}^\mathrm{T} \mathbf{k} \mathbf{n}^\mathrm{e}  \right) \widetilde{u}_D  \mathrm{d} \Gamma \\
& \qquad + \sum_{e \in \partial E \cap \Gamma_N} \int_e \left( 2 \overline{\psi} + \frac{\eta_2 h_e}{\overline{k}} \nabla \overline{\psi}^\mathrm{T} \mathbf{k} \mathbf{n}^\mathrm{e} \right)  \widetilde{q}_N  \mathrm{d} \Gamma  \\
& \qquad + \sum_{e \in \partial E \cap  \Gamma_R} \int_e h \left( 2 \overline{\psi} + \frac{\eta_2 h_e}{\overline{k}} \nabla \overline{\psi}^\mathrm{T} \mathbf{k} \mathbf{n}^\mathrm{e}  \right) \widetilde{u}_R  \mathrm{d} \Gamma,
\label{eqn:weakform3_5}
\end{split}
\end{align}
where $E$ is employed as the first neighboring subdomain ($E_1$) in the jump operator $\left[ \! \left[ \cdot \right] \! \right]$ and $\mathbf{n}^\mathrm{e}$.

The current approach is also named as a singular solution method (PG-FPM-3). The same as the finite volume method, the second derivatives of the temperature are absent, hence a linear trial function $u_h$ is sufficient in the singular solution method. As the test function $v$ is also independent in each subdomain and the weak-form can be written at a point level, the bandwidths of matrices $\mathbf{C}$ and $\mathbf{K}$ are minimized. Note that there are no point thermal conductivity matrices in the assembling of the global matrix $\mathbf{K}$, i.e., all the terms relating to $\mathbf{K}$ are defined on subdomain boundaries, thus the assembling process is significantly simplified. However, the test functions in the current PG-FPM-3 are not polynomial. As a result, more than one integration points may be required to assure the accuracy of the approach, especially for transient problems. The global heat capacity matrix $\mathbf{C}$ is no longer diagonal, which leads to a full Jacobian matrix and is harmful for the efficiency of the method. Numerical results have shown that one integration point is good enough for triangular or tetrahedron partitions. Hence partitions with more complicated subdomain geometries should be avoided in transient analysis with the present singular solution method (PG-FPM-3).

\section{Time discretization methods} \label{sec:time}

A set of semi-discrete equations are obtained after spatial discretization via the PG-FPMs:
\begin{align}
\begin{split}
\dot{\mathbf{u}} = - \mathbf{C} ^\mathrm{T} \mathbf{K} \mathbf{u} + \mathbf{C} ^\mathrm{T} \mathbf{q} = \mathbf{J} \mathbf{u} + \mathbf{C} ^\mathrm{T} \mathbf{q}, \quad t \in \left[ 0,  T \right].
\end{split}
\end{align}
This is a system of standard first-order ODEs, which can be further solved by Finite Difference methods or weak-form methods \cite{He1999, Elgohary2014} in the time domain. In the current work, we employ the Backward Euler (BE) scheme and a Local Variational Iteration Method (LVIM). The formulation and numerical implementation of the LVIM are presented in our previous works \cite{Guan2020, Wang2019} and thus omitted here. Note that the LVIM can be remarkably more efficient than the classic Backward Euler scheme, especially for problems with a sparse Jacobian matrix. The enhanced matrix $\widetilde{\mathbf{J}}$ remains sparse in the LVIM, hence the approach is sensitive to the bandwidths of $- \mathbf{C}^{-1} \mathbf{K}$ achieved by spatial discretization. A smaller bandwidth achieved by the finite volume or singular solution PG-FPMs can help to improve the efficiency of the transient analysis significantly.

\section{2D examples} \label{sec:2D}

In the following sections, a number of numerical results of the PG-FPM approaches are presented, as well as a comparison with the Galerkin FPM and the FEM results (achieved by ABAQUS). For simplicity, the three approaches are named as PG-FPM-1 (the collocation method), PG-FPM-2 (the finite volume method) and PG-FPM-3 (the singular solution method) respectively. The relative errors $e_0$ and $e_1$ used in the following sections are defined as:
\begin{align}
\begin{split}
e_0 = \frac{\left\| u_h - u \right\| _{2}}{\left\| u \right\| _{2}},  \qquad e_1 = \frac{\left\| \nabla u_h - \nabla u \right\| _{2}}{\left\| \nabla u \right\| _{2}}
\end{split}
\end{align}
where $u$ and $u_h$ are the exact and computed solutions, respectively, and 
\begin{align}
\begin{split}
\left\| u \right\| _{2} = \left( \int_{\Omega} u^2 \mathrm{d} \Omega \right)^{1/2}, \qquad \left\| \nabla u \right\| _{2} = \left( \int_{\Omega} \left\| \nabla u \right\|^2 \mathrm{d} \Omega \right)^{1/2}. \notag
\end{split}
\end{align}

\subsection{Isotropic homogeneous examples}

\begin{figure}[htbp] 
  \centering 
    \subfigure[]{ 
    \label{fig:Ex11_01} 
    \includegraphics[width=0.48\textwidth]{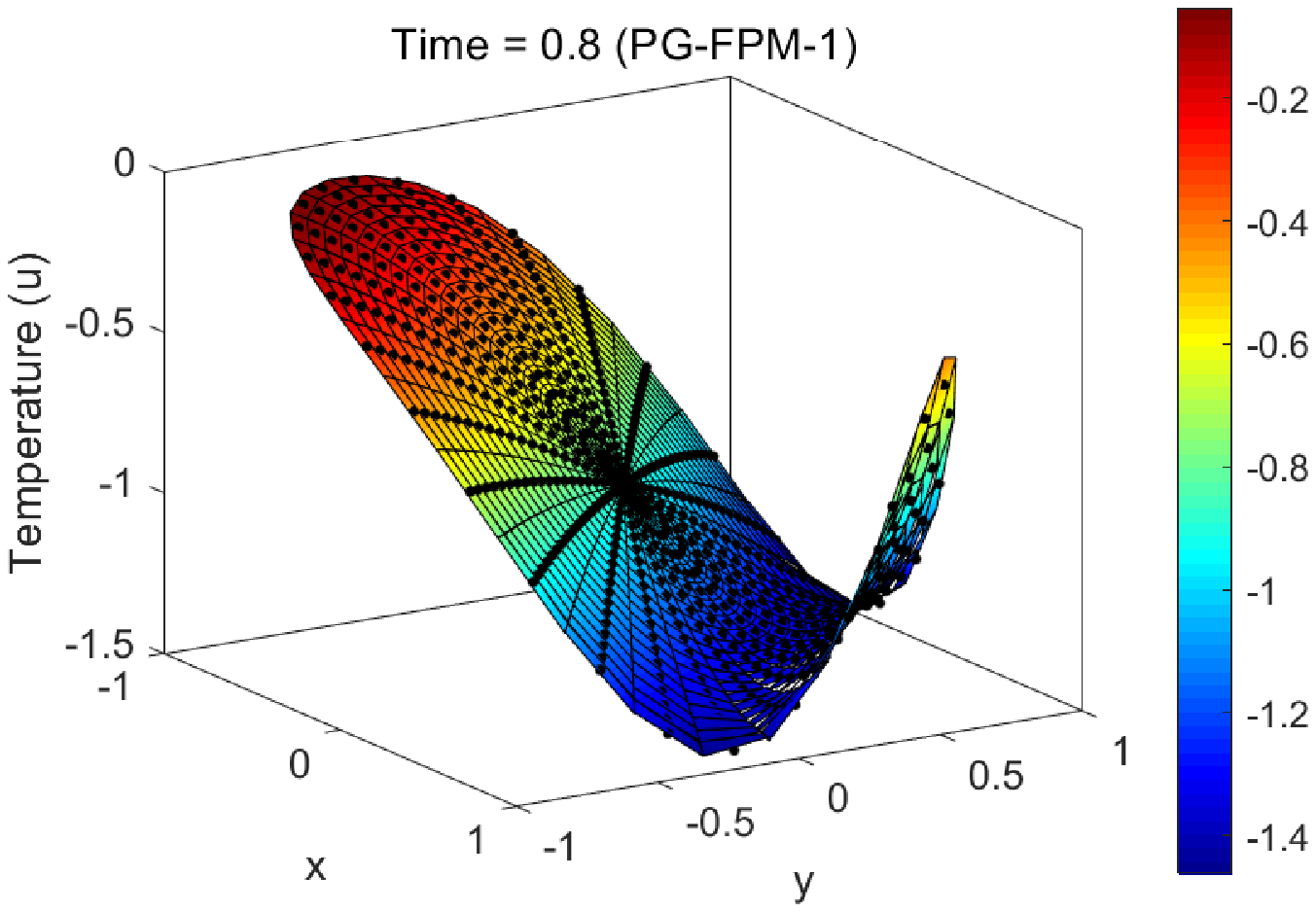}}  
    \subfigure[]{ 
    \label{fig:Ex11_02} 
    \includegraphics[width=0.48\textwidth]{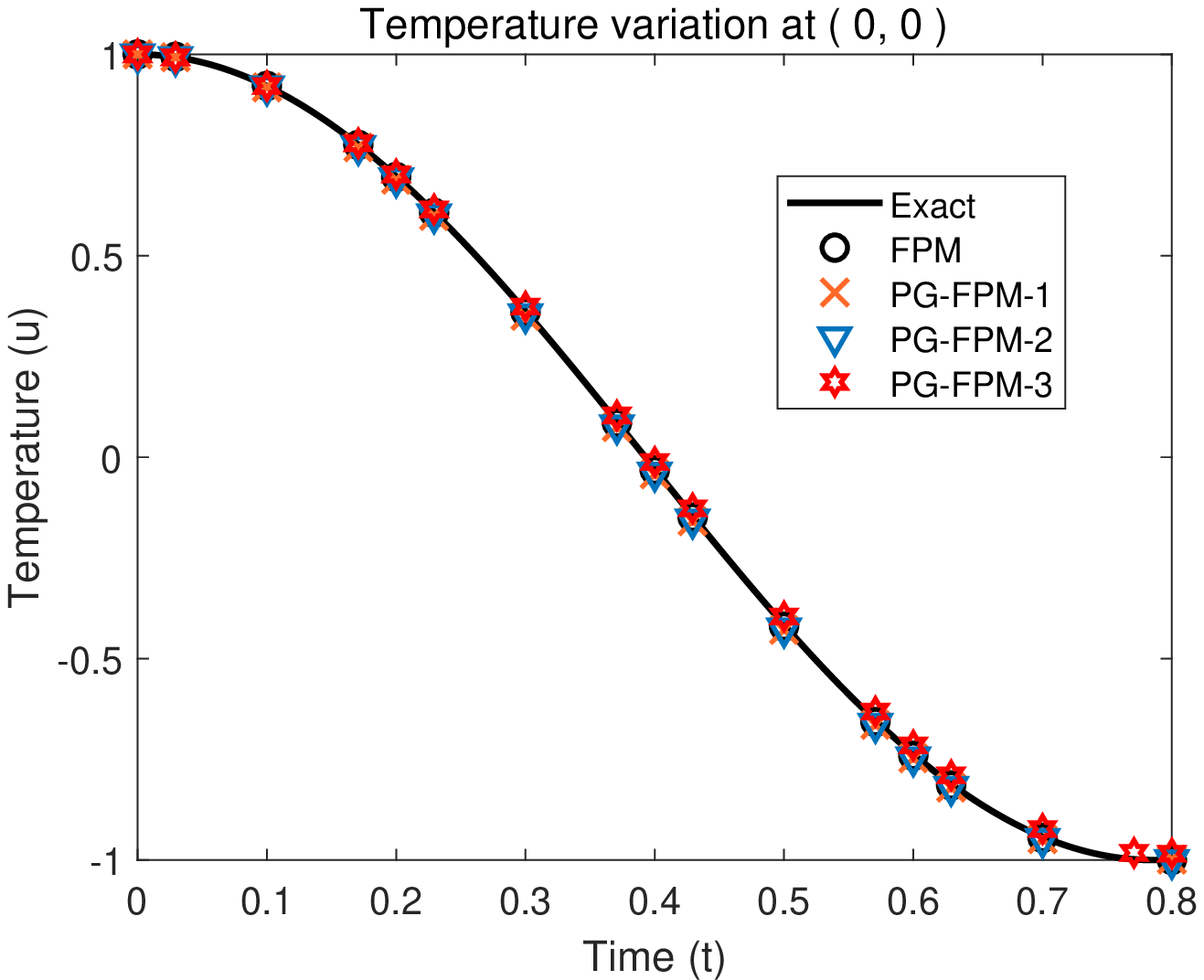}}  
  \caption{Ex. (1.1) - The computed solution. (a) spatial temperature distribution achieved by PG-FPM- 1 when $t = 0.8$. (b) transient temperature solution achieved by the FPM and PG-FPMs.}
  \label{fig:Ex11} 
\end{figure}

First, a benchmark 2D isotropic homogenous problem in a circular domain is under study. Without loss of generality, we assume the material density $\rho = 1$, specific heat capacity $c=1$, thermal conductivity tensor components $k_{11} = k_{22}=1$, $k_{12} = k_{21}=0$. The body source density $Q$ is absent. We consider a postulated analytical solution \cite{Johansson2011, Reeve2013}:
\begin{align}
\begin{split}
u(x,y,t) = e^{x+y} \mathrm{cos} (x+y+4t), \quad (x,y) \in \left\{ (x,y) \mid x^2 + y^2 \leq 1 \right\}.
\end{split}
\end{align}
Essential boundary conditions are prescribed on all external boundaries. The original Galerkin FPM and the three PG-FPMs are employed in the spatial discretization, while in the time domain, the LVIM is applied. A total of 600 Fragile Points are distributed uniformly in the domain. A mixed quadrilateral and triangular partition is exploited in all the approaches. In this example, all the Fragile Points are placed at the centroids of the subdomains. There is no Point exists on the external boundaries, hence the essential boundary conditions are applied by Interior Penalty (IP) Numerical Flux Corrections. The computational parameters used in the LVIM (see our previous study \cite{Guan2020} for details) are given as: the number of collocation points in each time interval $M = 5$, and the error tolerance in stopping criteria $tol = 1 \times 10^{-6}$. In the meshless Differential Quadrature (DQ) approximation in the PG-FPM-1 approach (the collocation method), the multiquadric Radial Basis Function (MQ-RBS) is employed, in which the nondimensional constant $c = 4$. These values maintain the same in all the following examples, unless otherwise stated. Here the time step $\Delta t = 0.2$. Table~\ref{table:Ex1} and Fig.~\ref{fig:Ex11} exhibit the solutions based on the FPM and PG-FPMs. As can be seen, the result of all the methods shows great consistency with the exact solution, as well as the numerical results achieved in previous literatures \cite{Johansson2011, Reeve2013, Guan2020}.

In Table~\ref{table:Ex1}, $\eta_1$ and $\eta_2$ denote the nondimensional penalty parameters used in the FPM and PG-FPMs. The recommended ranges of $\eta_1$ and $\eta_2$ will be given in section~\ref{sec:PS}. For the original FPM and PG-FPM-2 / 3 approaches, a large enough $\eta_1$ is required to ensure stability. Yet in the PG-FPM-1 (the collocation method), the required penalty parameter $\eta_1$ may be as small as zero, since the quadratic trial functions used in the PG-FPM-1 result in a better continuity across the subdomain boundaries. The PG-FPM-1 also presents a better estimate of the gradient of temperature comparing with the other methods. Moreover, its efficiency remains the same for any point distribution and domain partitions, whereas the efficiency of the Galerkin FPM and PG-FPM-2 / 3 decreases significantly when the Fragile Points do not coincide with the centroids of subdomains, as a result of the non-diagonal heat capacity matrix $\mathbf{C}$ and full Jacobian matrix ($\mathbf{J} = - \mathbf{C}^{-1} \mathbf{K}$). Neverthless, as compared with the finite volume method (PG-FPM-2) and the singular solution method (PG-FPM-3), the collocation method suffers from a more complicated local discretization process and leads to a larger bandwidth of the thermal conductivity matrix $\mathbf{K}$. 

The finite volume method (PG-FPM-2) shows the highest efficiency among the four approaches, saving approximately 25\% computational time and achieving similar accuracy as compared to the Galerkin FPM. With only one integration point adopted in each subdomain, the singular solution method (PG-FPM-3) presents similar computing efficiency as the finite volume method in this example. However, in some other transient problems, as a result of the complicated local fundamental solution, more than one integration points are required, and the global heat capacity matrix $\mathbf{C}$ in the PG-FPM-3 may have a larger bandwidth than all the other methods. This drawback is also seen in most of the previous meshless methods (EFG \cite{Belytschko1994}, MLPG \cite{Atluri1998, Sladek2008}, etc.), but is omitted in the FPM and PG-FPM-1 / 2 approaches. Besides, the PG-FPM-3 approach does not suffer this complicated integration problem in steady-state analysis, since there is no volume integration term in assembling the thermal conductivity matrix. Therefore, the singular solution method (PG-FPM-3) is recommended for steady-state analysis.

\begin{table}[htbp]
\caption{Relative errors and computational times of the FPM and PG-FPMs in solving Ex.~(1.1).}
\centering
{
\begin{tabular*}{500pt}{@{\extracolsep\fill}lcccc@{\extracolsep\fill}}
\toprule
\textbf{Method} & \tabincell{c}{\textbf{Computational} \\ \textbf{parameters}} &  \textbf{Relative errors} & \tabincell{c}{\textbf{Number of nonzero} \\ \textbf{elements in} $\mathbf{K}$ \& $\mathbf{C}$ \\  \textbf{in each column} \\ $N_{band} (\mathbf{K})$ \& $N_{band} (\mathbf{C})$} & \tabincell{c}{\textbf{Computational} \\ \textbf{time (s)}} \\
\midrule
FPM & $\eta_1 = 1$, $\eta_2 = 1 \times 10^{5}$ & \tabincell{c}{$e_0 = 5.2 \times 10^{-3}$\\$e_1 = 1.5 \times 10^{-1}$} & 24 \& 1 & 2.8 \\
PG-FPM-1 & $\eta_1 = 0$, $\eta_2 = 1 \times 10^{5}$ & \tabincell{c}{$e_0 = 3.7 \times 10^{-3}$\\$e_1 = 3.5 \times 10^{-2}$} & 17 \& 1 & 2.7 \\
PG-FPM-2 & $\eta_1 = 1$, $\eta_2 = 1 \times 10^{5}$ & \tabincell{c}{$e_0 = 8.6 \times 10^{-3}$\\$e_1 = 1.7 \times 10^{-1}$} & 12 \& 1 & 2.1 \\
PG-FPM-3 & $\eta_1 = 1$, $\eta_2 = 1 \times 10^{5}$ & \tabincell{c}{$e_0 = 9.0 \times 10^{-3}$\\$e_1 = 1.4 \times 10^{-1}$} & 12 \& 1 & 2.2 \\
\bottomrule
\end{tabular*}}
\label{table:Ex1}
\end{table}

\begin{figure}[htbp] 
  \centering 
    \subfigure[]{ 
    \label{fig:Ex12_01} 
    \includegraphics[width=0.48\textwidth]{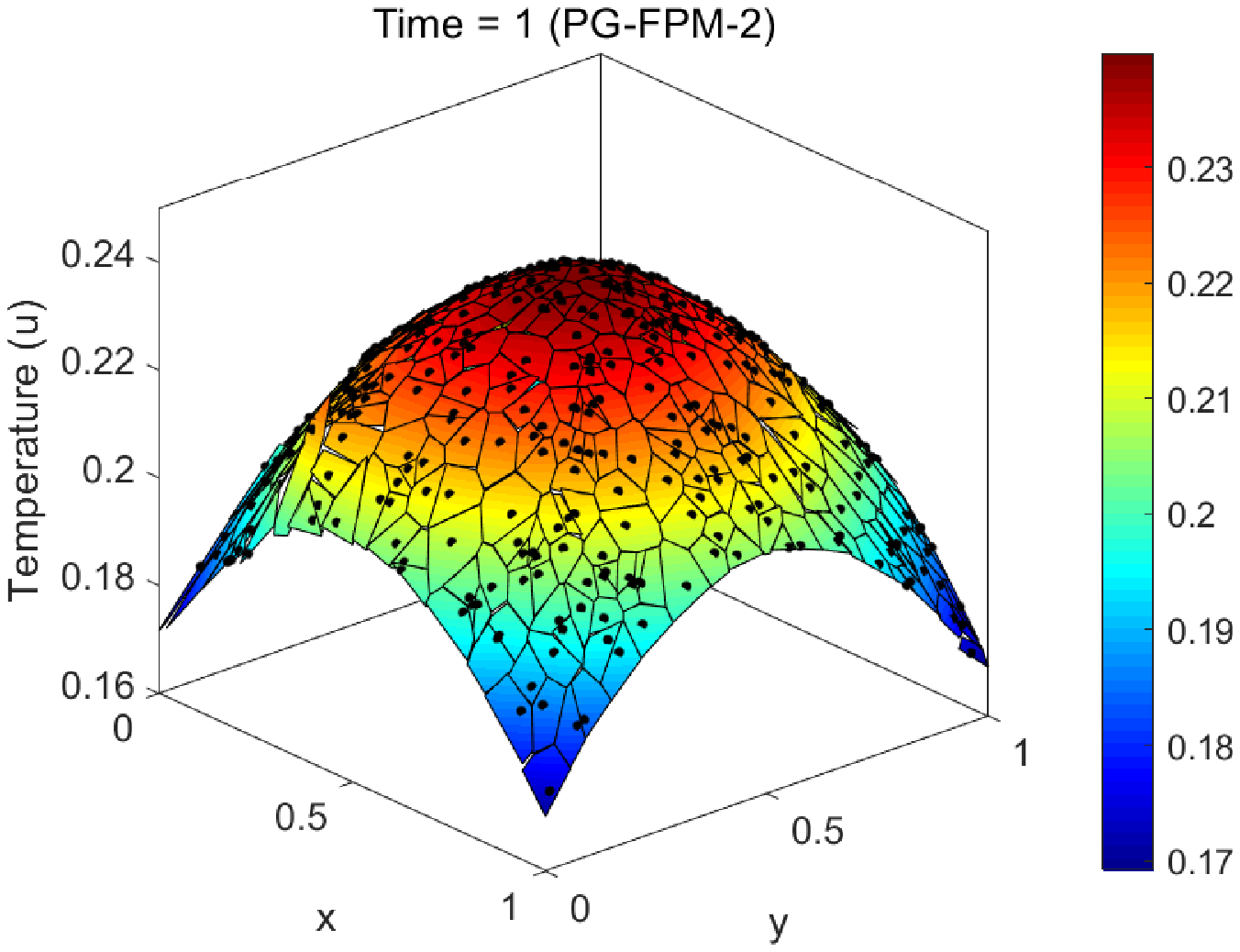}}  
    \subfigure[]{ 
    \label{fig:Ex12_02} 
    \includegraphics[width=0.48\textwidth]{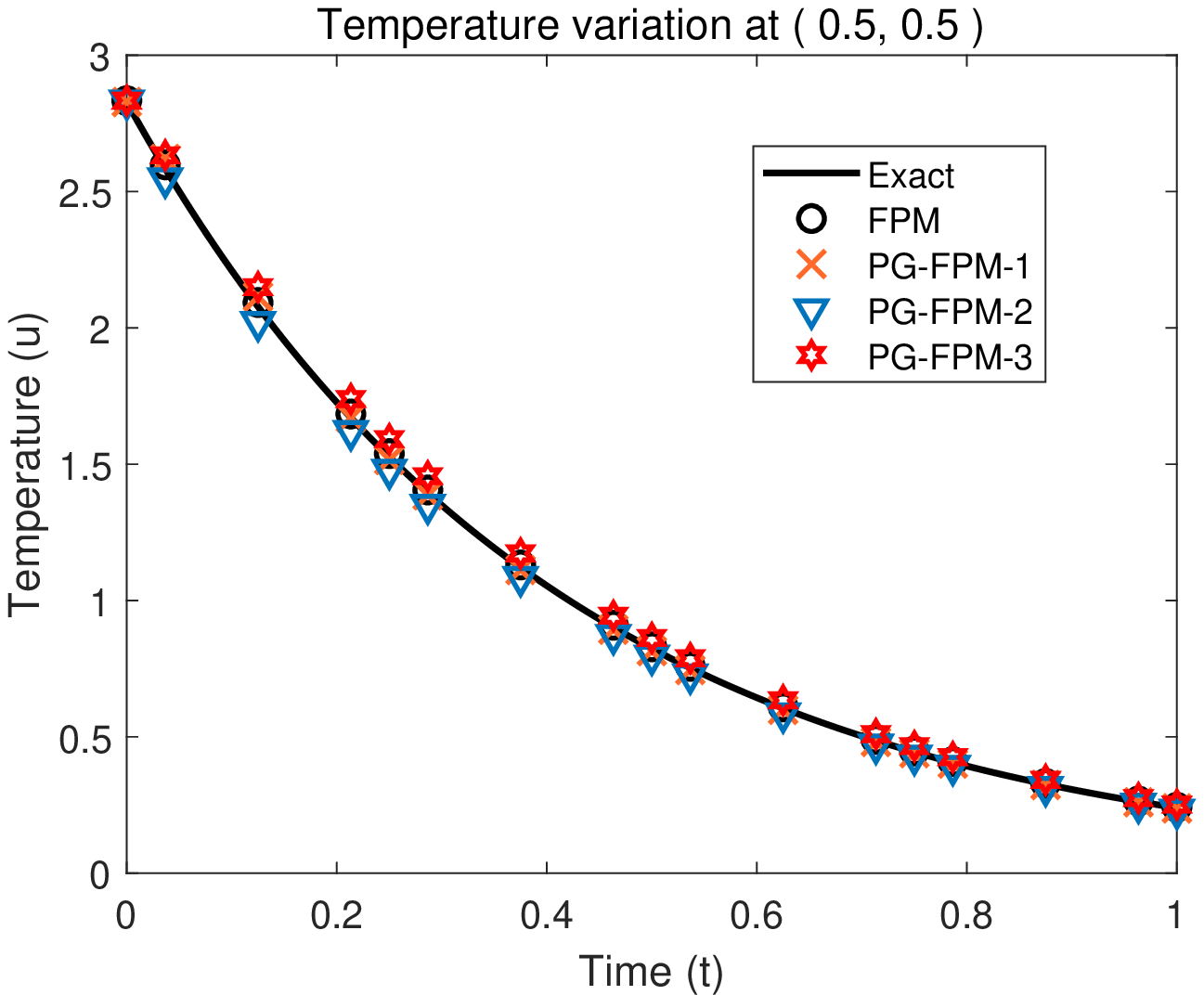}}  
  \caption{Ex. (1.2) - The computed solution with 400 random Points. (a) spatial temperature distribution achieved by PG-FPM- 2 when $t = 1$. (b) transient temperature solution achieved by the FPM and PG-FPMs.}
  \label{fig:Ex12} 
\end{figure}

The second numerical example (Ex.~(1.2)) is in a square domain. We consider the same material properties used in Ex.~(1.1) and the following postulated analytical solution \cite{Johansson2011}:
\begin{align}
\begin{split}
u(x,y,t) = \sqrt{2} e^{-\pi^2 t /4} \left[ \mathrm{cos} (\frac{\pi x}{2} - \frac{\pi}{4}) + \mathrm{cos} (\frac{\pi y}{2} - \frac{\pi}{4}) \right], \\
\quad (x,y) \in \left\{ (x,y) \mid x \in \left[ 0, 1 \right], \; y \in \left[ 0, 1 \right] \right\}.
\end{split}
\end{align}
Neumann boundary condition is applied on $x = 1$, while the other boundaries are under Dirichlet boundary conditions. 400 Fragile Points are distributed uniformly or randomly in the domain. Voronoi Diagram partition is utilized. The computed solutions achieved by the original Galerkin FPM and PG-FPMs are presented in Table~\ref{table:Ex12_1}, \ref{table:Ex12_2} and Fig.~\ref{fig:Ex12}. All the methods show excellent accuracy comparing with the exact solution. The spatial temperature distribution achieved by different approaches are generally identical, hence only the PG-FPM-2 result is shown in Fig.~\ref{fig:Ex12_01} as an example. Under uniform point distribution (see Table~\ref{table:Ex12_1}), the same as Ex.~(1.1), the finite volume method (PG-FPM-2) achieves the highest efficiency, followed by the singular solution method (PG-FPM-3). Whereas when random distributed points are employed (see Table~\ref{table:Ex12_2}), the collocation method (PG-FPM-1) shows the highest efficiency. The computational time is cut by a half by the PG-FPM-1, as a consequence of the diagonal heat capacity matrix. All the proposed three PG-FPM approaches improve the performance of the original Galerkin FPM in some aspects.

\begin{table}[htbp]
\caption{Relative errors and computational times of the FPM and PG-FPMs in solving Ex.~(1.2) with 400 uniform Points.}
\centering
{
\begin{tabular*}{500pt}{@{\extracolsep\fill}lcccc@{\extracolsep\fill}}
\toprule
\textbf{Method} & \tabincell{c}{\textbf{Computational} \\ \textbf{parameters}} &  \textbf{Relative errors} & \tabincell{c}{$N_{band} (\mathbf{K})$ \& $N_{band} (\mathbf{C})$}  & \tabincell{c}{\textbf{Computational} \\ \textbf{time (s)}} \\
\midrule
FPM & $\eta_1 = 1$, $\eta_2 = 1 \times 10^{5}$ & \tabincell{c}{$e_0 = 8.6 \times 10^{-4}$\\$e_1 = 1.3 \times 10^{-2}$} & 14 \& 1 & 1.5 \\
PG-FPM-1 & $\eta_1 = 0$, $\eta_2 = 1 \times 10^{5}$ & \tabincell{c}{$e_0 = 5.6 \times 10^{-3}$\\$e_1 = 4.1 \times 10^{-2}$} & 13 \& 1 & 1.4 \\
PG-FPM-2 & $\eta_1 = 1$, $\eta_2 = 1 \times 10^{5}$ & \tabincell{c}{$e_0 = 5.1 \times 10^{-4}$\\$e_1 = 1.9 \times 10^{-2}$} & 9 \& 1 & 1.0 \\
PG-FPM-3 & $\eta_1 = 1$, $\eta_2 = 1 \times 10^{5}$  & \tabincell{c}{$e_0 = 9.9 \times 10^{-4}$\\$e_1 = 1.4 \times 10^{-2}$} & 9 \& 1 & 1.1 \\
\bottomrule
\end{tabular*}}
\label{table:Ex12_1}
\end{table}

\begin{table}[htbp]
\caption{Relative errors and computational times of the FPM and PG-FPMs in solving Ex.~(1.2) with 400 random Points.}
\centering
{
\begin{tabular*}{500pt}{@{\extracolsep\fill}lcccc@{\extracolsep\fill}}
\toprule
\textbf{Method} & \tabincell{c}{\textbf{Computational} \\ \textbf{parameters}} &  \textbf{Relative errors} & \tabincell{c}{$N_{band} (\mathbf{K})$ \& $N_{band} (\mathbf{C})$}  & \tabincell{c}{\textbf{Computational} \\ \textbf{time (s)}} \\
\midrule
FPM & $\eta_1 = 1$, $\eta_2 = 1 \times 10^{5}$ & \tabincell{c}{$e_0 = 2.2 \times 10^{-3}$\\$e_1 = 9.9 \times 10^{-2}$} & 37 \& 19 & 3.1 \\
PG-FPM-1 & $\eta_1 = 0.1$, $\eta_2 = 1 \times 10^{5}$ & \tabincell{c}{$e_0 = 6.6 \times 10^{-3}$\\$e_1 = 7.6 \times 10^{-2}$} & 22 \& 1 & 1.5 \\
PG-FPM-2 & $\eta_1 = 1$, $\eta_2 = 1 \times 10^{5}$ & \tabincell{c}{$e_0 = 8.6 \times 10^{-4}$\\$e_1 = 9.9 \times 10^{-2}$} & 19 \& 7 & 2.8 \\
PG-FPM-3 & $\eta_1 = 1$, $\eta_2 = 1 \times 10^{5}$  & \tabincell{c}{$e_0 = 1.0 \times 10^{-2}$\\$e_1 = 1.1 \times 10^{-1}$} & 19 \& 7 & 3.1 \\
\bottomrule
\end{tabular*}}
\label{table:Ex12_2}
\end{table}

\begin{figure}[htbp] 
  \centering 
    \subfigure[]{ 
    \label{fig:Ex12_spy_00} 
    \includegraphics[width=0.48\textwidth]{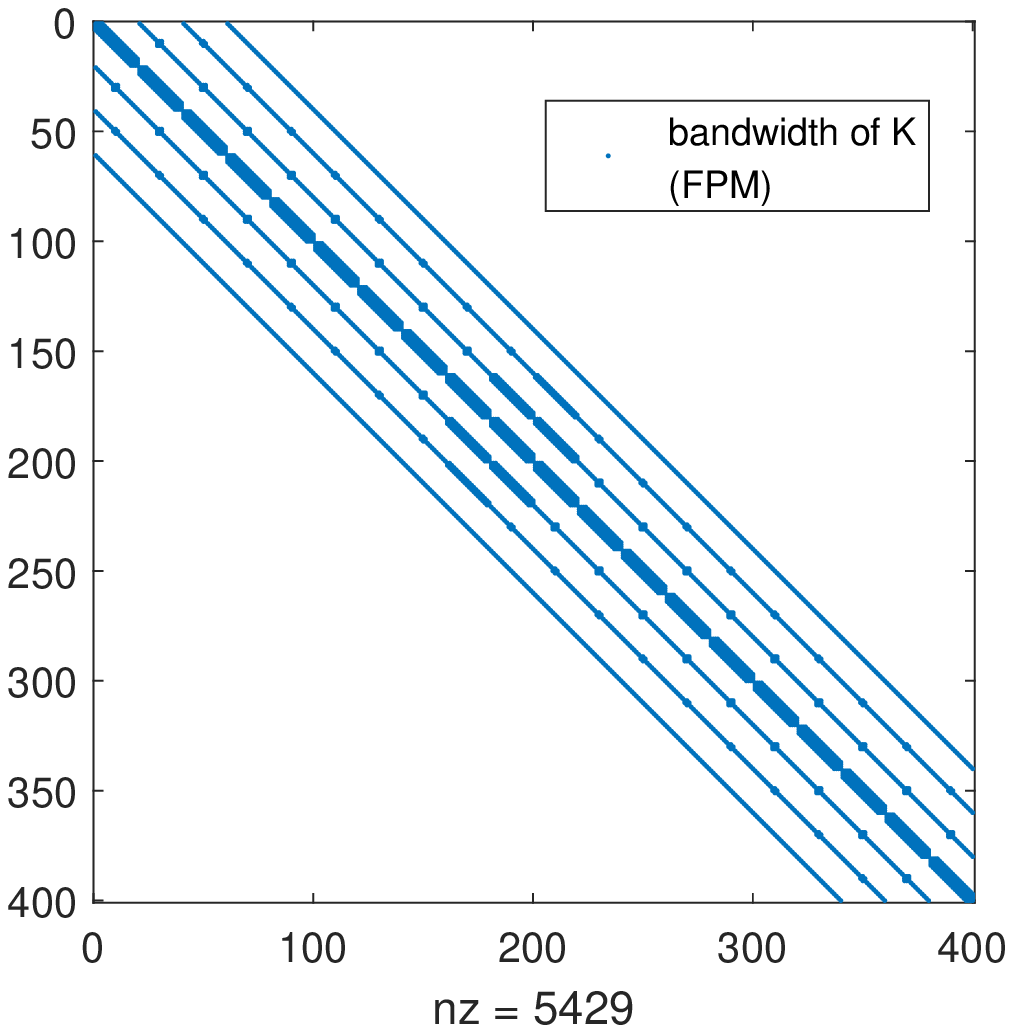}}  
    \subfigure[]{ 
    \label{fig:Ex12_spy_01} 
    \includegraphics[width=0.48\textwidth]{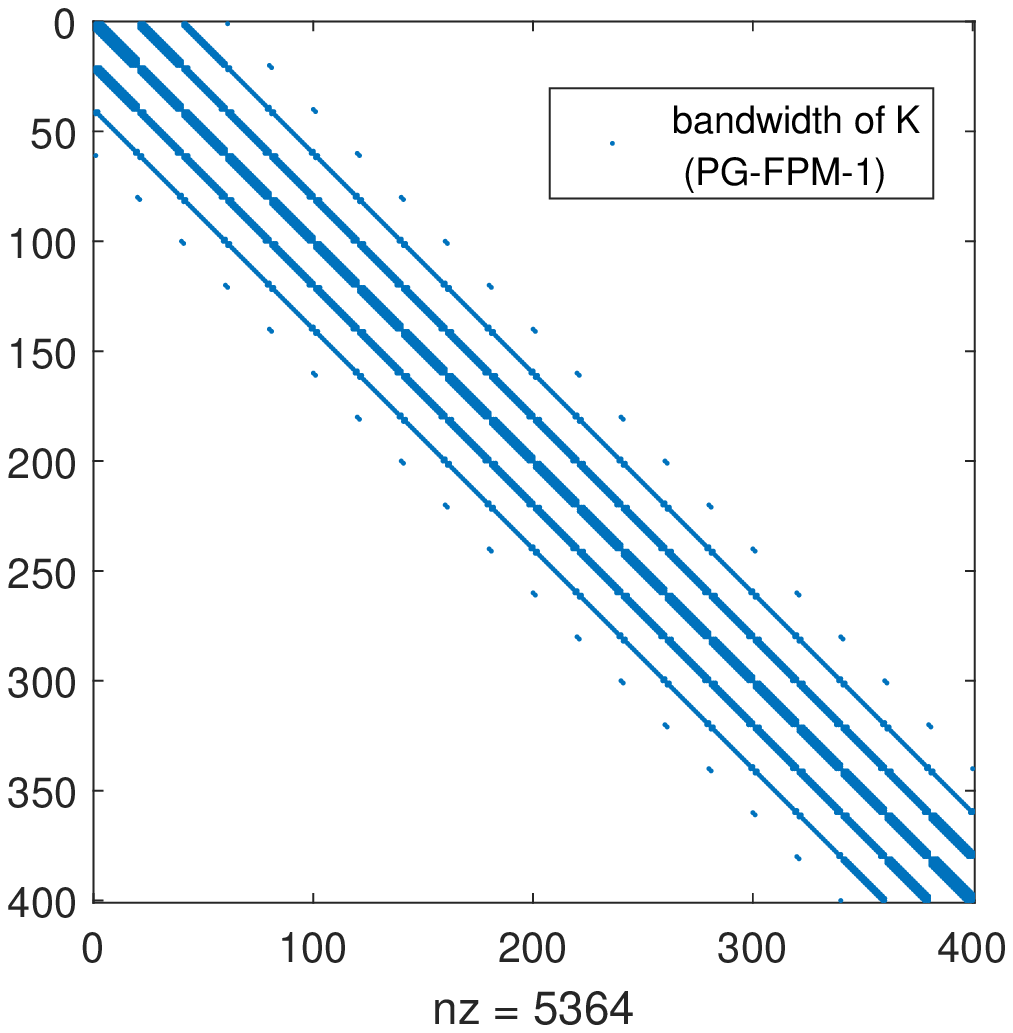}}  
    \subfigure[]{ 
    \label{fig:Ex12_spy_02} 
    \includegraphics[width=0.48\textwidth]{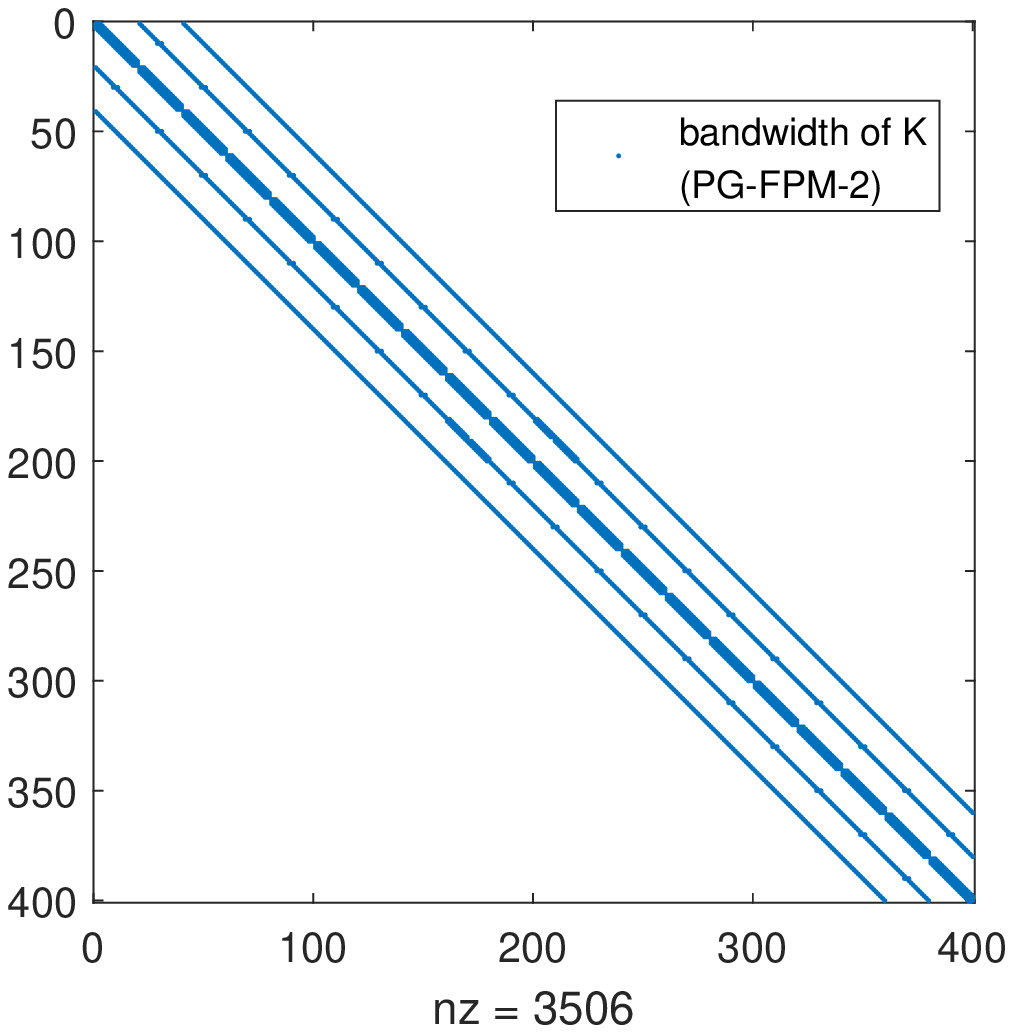}}  
    \subfigure[]{ 
    \label{fig:Ex12_spy_03} 
    \includegraphics[width=0.48\textwidth]{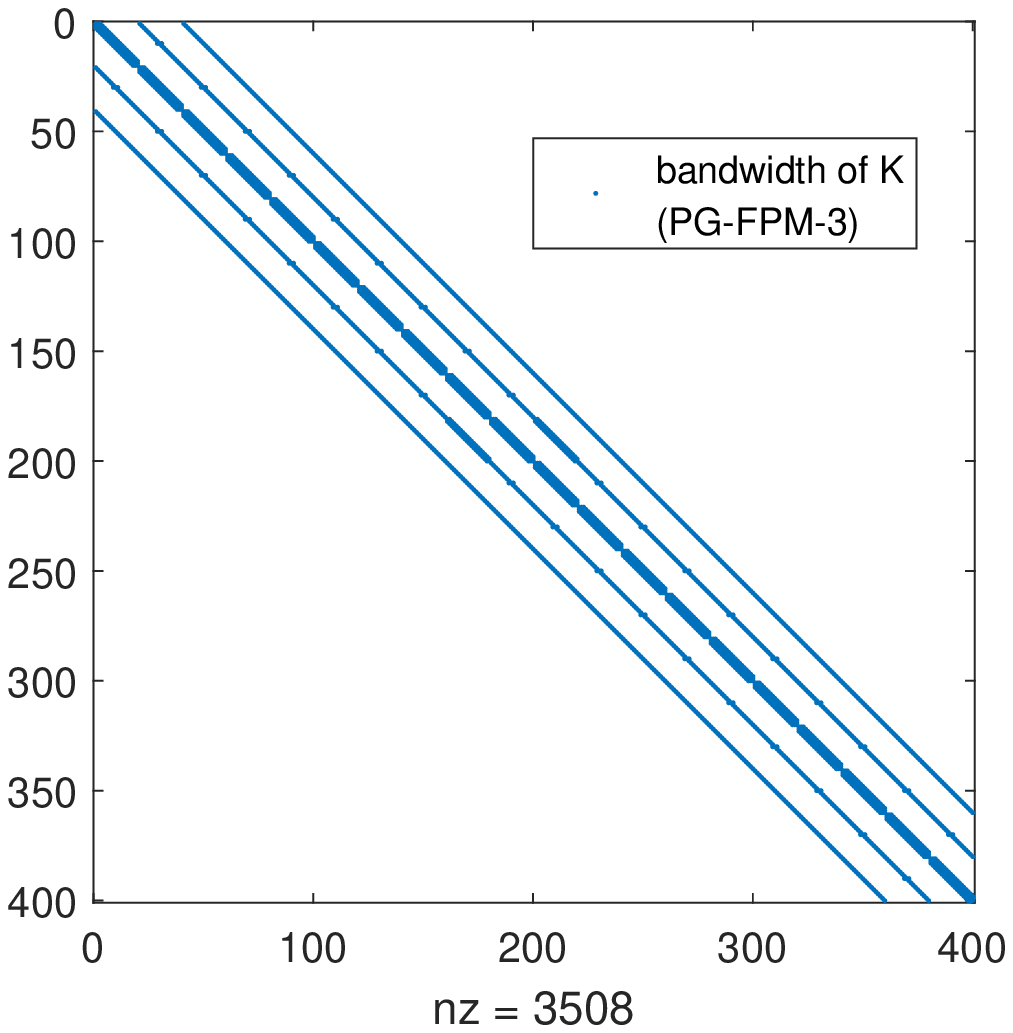}}  
  \caption{Ex. (1.2) - Visualized sparsity patterns of thermal conductivity matrix $\mathbf{K}$ in the original FPM and PGFPMs. (a) original FPM. (b) PG-FPM-1. (c) PG-FPM-2. (d) PG-FPM-3.}
  \label{fig:Ex11_spy} 
\end{figure}

Figure~\ref{fig:Ex11_spy} exhibits the visualized sparsity patterns and the numbers of non-zero elements in the thermal conductivity matrix $\mathbf{K}$ in the FPM and PG-FPMs. As can be seen intuitively, the PG-FPM-2 and PG-FPM-3 approaches decrease the bandwidth of $\mathbf{K}$ remarkably. The superiority of these highly sparse matrices can be more significant in problems with more Fragile Points and higher nonlinearity in the time domain (when more collocation points in each time interval in LVIM is applied).

\subsection{Anisotropic nonhomogeneous examples in a square domain}

In the following four examples, a mixed boundary value problem in anisotropic nonhomogeneous media is considered. The tested domain is a $L \times L$ square. The boundary and initial conditions are:
\begin{align}
\begin{split}
& u(x, 0, t) = u_0, \qquad  u(x, L, t) =  u_L, \qquad u_{,x}(0, y, t) = 0, \qquad  u_{,x} (L, y, t) =  0, \\
& u(x, y ,0) =  u_0,
\end{split}
\end{align}
where $u_{,x}$ is the partial derivative to $x$ of the temperature field $u$. Dirichlet boundary conditions are given on $y = 0$ and $y = L$. And symmetric boundary conditions are applied on the lateral sides. In isotropic media, the symmetry is equivalent to Neumann boundary conditions with $\widetilde{q}_N = 0$.

The material is functionally graded (FG) \cite{Miyamoto2013}, with the following material properties:
\begin{align}
\begin{split}
&\rho (x, y) =  1, \qquad c(x, y) = f(y), \qquad \mathbf{k} (x, y) = f(y) \left[ \begin{matrix} \hat{k}_{11} & \hat{k}_{12} \\ \hat{k}_{21} & \hat{k}_{22} \end{matrix} \right],
\end{split}
\end{align}
where $f(y)$ is the gradation function. In isotropic case, $\hat{k}_{ij} = \delta_{ij}$. Whereas in anisotropic case, $\hat{k}_{11} = \hat{k}_{22} = 2, \hat{k}_{12} = \hat{k}_{21} = 1$. The body source density $Q$ is absent. It is obvious that the resulting temperature field is independent of $x$, i.e., the example is equivalent to a 1D heat conduction problem.

The material gradation function $f(y)$ and boundary values used in Ex.~(1.3) – (1.6) are:
\begin{align}
 \text{Ex.~(1.3):} \quad  & \text{exponential}:  & f(y) = \mathrm{exp}(\delta y / L) &, \delta = 3, u_0 = 1, u_L =20; & \notag \\
 \text{Ex.~(1.4):} \quad  & \text{exponential}:  & f(y) = \left[ \mathrm{exp} (\delta y / L) + 5 \mathrm{exp} (-\delta y / L) \right] ^2 &, \delta = 2, u_0 = 1, u_L =20; & \notag \\
 \text{Ex.~(1.5):} \quad  & \text{trigonometric}: & f(y) = \left[ \mathrm{cos} (\delta y / L) + 5 \mathrm{sin} (\delta y /L) \right]^2 &, \delta = 2, u_0 = 0, u_L =100; & \notag \\
 \text{Ex.~(1.6):} \quad & \text{power-law}: & f(y) = \left( 1 + \delta y / L \right)^2 &, \delta = 3, u_0 = 1, u_L =20; & \notag
 \label{eqn:BCs}
\end{align}

These examples have also been studied using the Local Boundary Integral Equation (LBIE) Method \cite{Sladek2003}, the MLPG Method \cite{Sladek2004, Mirzaei2011, Mirzaei2014}, the meshless point interpolation method (PIM) \cite{Sladek2005} and the Galerkin FPM \cite{Guan2020}. These previous studies present consistent solutions. Therefore, in the current paper, the numerical results of the Galerkin FPM and the exact solutions \cite{Sladek2005} are exploited as benchmarks.

\begin{figure}[htbp] 
  \centering 
    \subfigure[]{ 
    \label{fig:Ex13_Trans} 
    \includegraphics[width=0.48\textwidth]{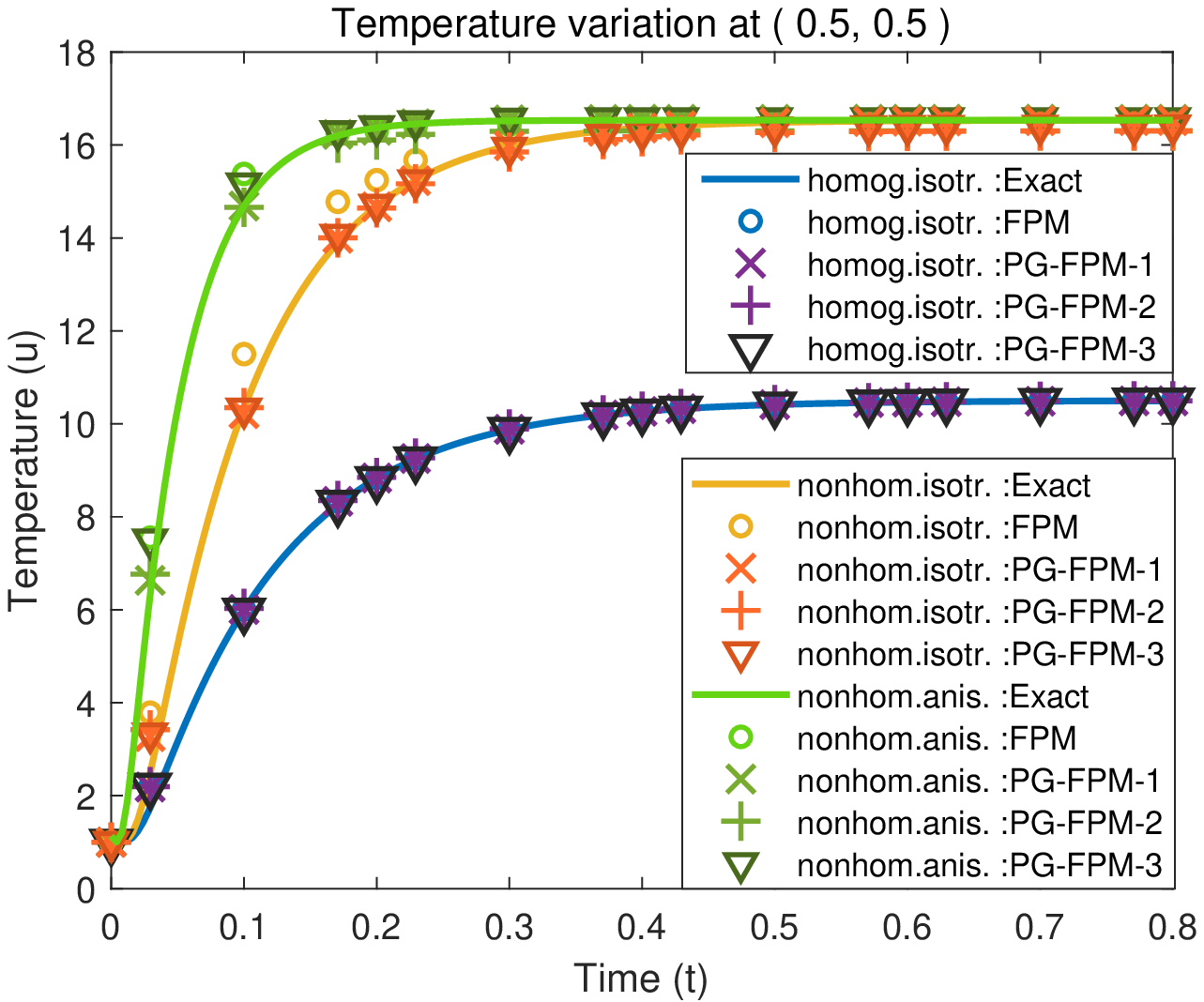}}  
    \subfigure[]{ 
    \label{fig:Ex13_Conf} 
    \includegraphics[width=0.48\textwidth]{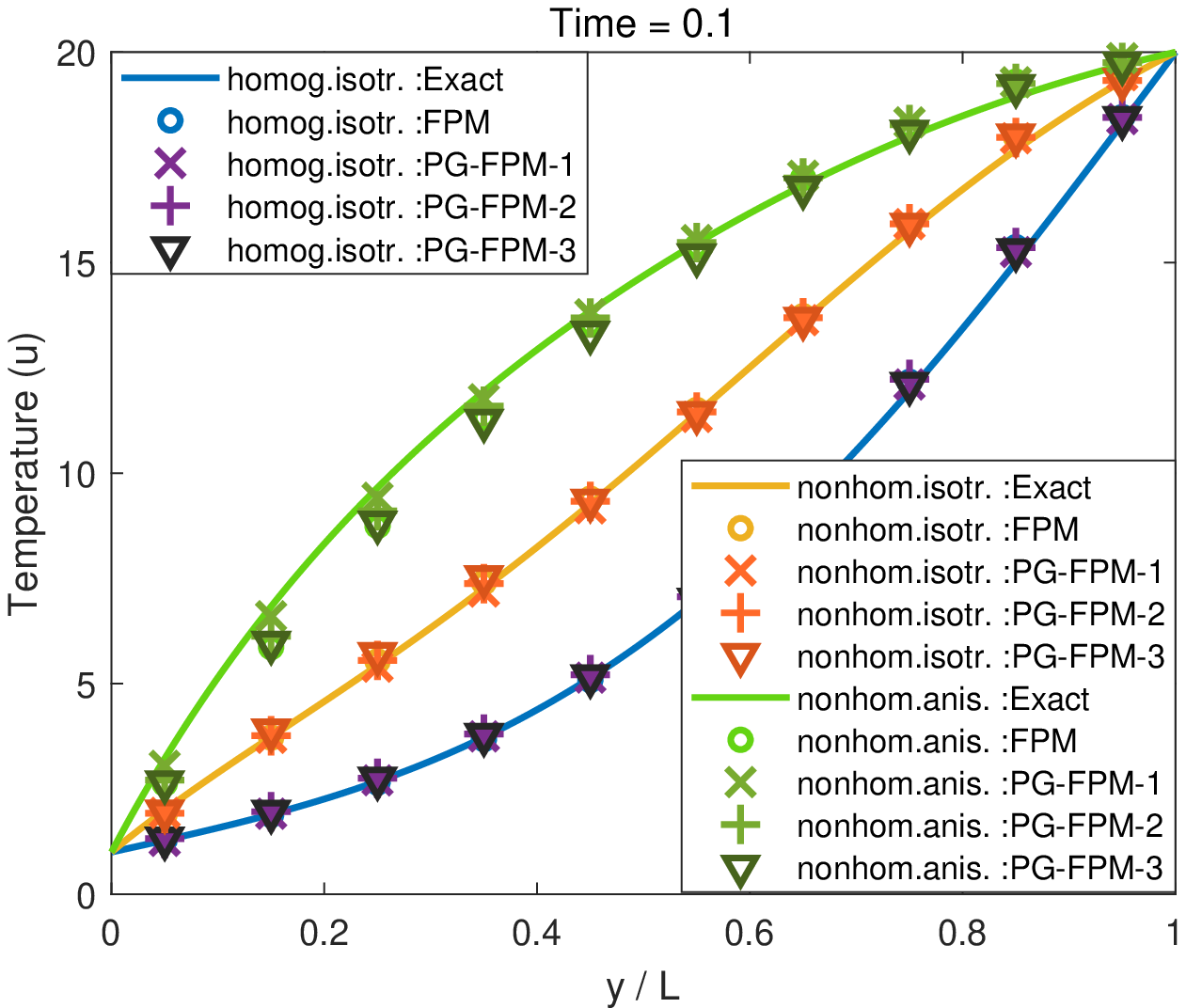}}  
  \caption{Ex. (1.3) - The computed solution achieved by the FPM and PG-FPMs with different material properties. (a) transient temperature solution at the midpoint of the domain in time scope $[0, 0.8]$. (b) vertical temperature distribution when $t=0.1$.} 
  \label{fig:Ex13} 
\end{figure}

Figure~\ref{fig:Ex13} shows the computed solution for Ex.~(1.3) in isotropic homogenous, isotropic nonhomogeneous, and anisotropic nonhomogeneous materials obtained by the Galerkin FPM and PG-FPMs. With 100 uniformly distributed Points, the results exhibit great agreement with the exact solution. Table~\ref{table:Ex13} shows the computational parameters, computational times, and average relative errors $\overline{e}_0$ in time interval $[0, 0.8]$ for all the methods with 400 Points. In the time domain, the LIVM approach is employed, with a time step $\Delta t = 0.1$. The collocation method (PG-FPM-1) acquires the best accuracy under various material properties using the same computational time as the original Galerkin FPM. Whereas the finite volume method (PG-FPM-2) and singular solution method (PG-FPM-3) improves the efficiency of the original Galerkin FPM by 30\%. While the nonhomogeneity and anisotropy of the material have a significant influence on the temperature distribution, they have little effects on the performance of proposed methods.

\begin{figure}[htbp] 
  \centering 
    \subfigure[]{ 
    \label{fig:Ex13_Sym} 
    \includegraphics[width=0.48\textwidth]{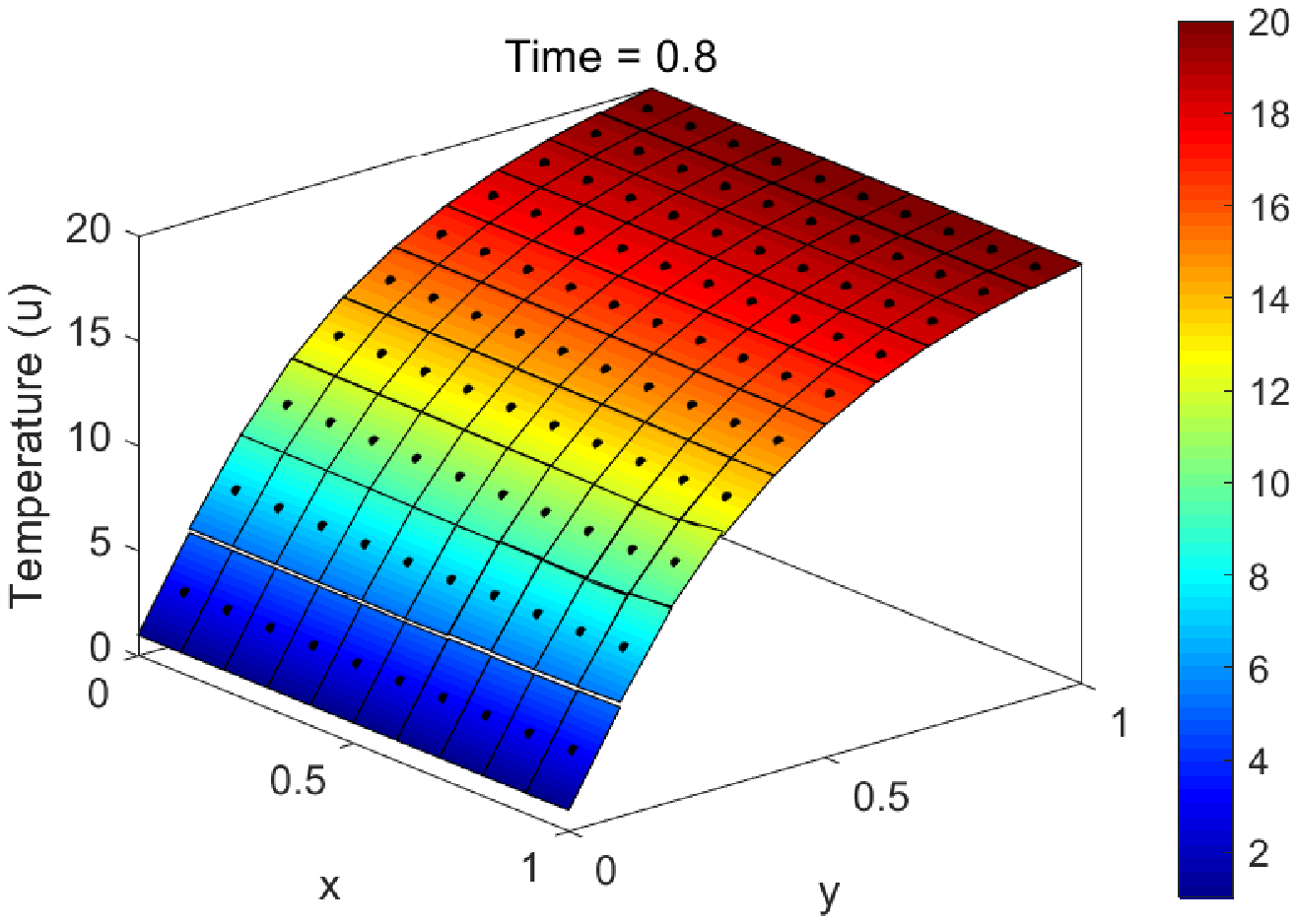}}  
    \subfigure[]{ 
    \label{fig:Ex13_Free} 
    \includegraphics[width=0.48\textwidth]{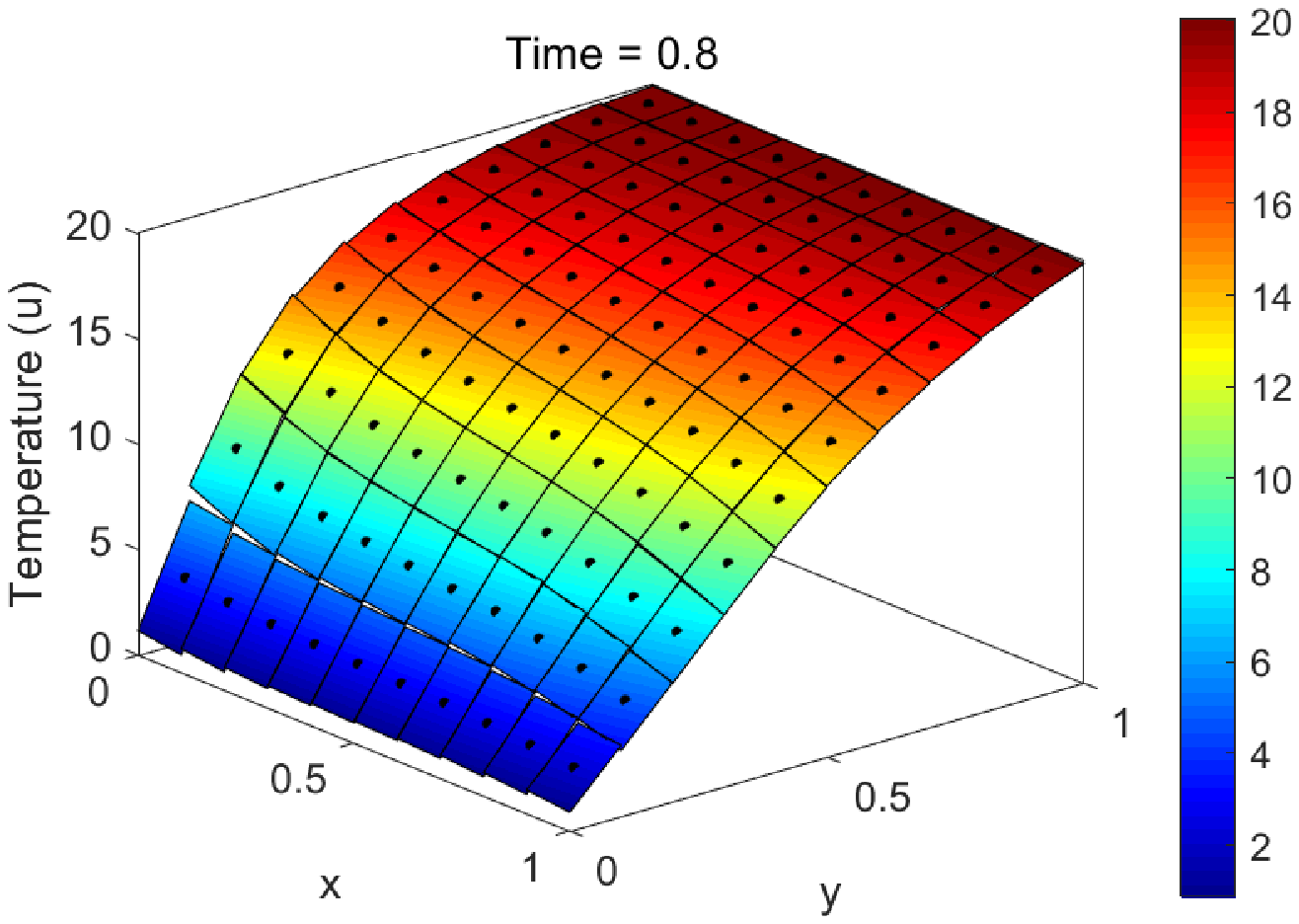}}  
  \caption{Ex. (1.3) - The computed solutions achieved by PG-FPM-3: a comparison of symmetric and free boundary conditions on the lateral sides in nonhomogeneous anisotropic material. (a) symmetric boundary conditions. (b) free  boundary conditions ($\widetilde{q}_N = 0$).} 
  \label{fig:Ex13_sym} 
\end{figure}

As have been stated, the symmetric boundary conditions are equivalent to free boundary conditions in isotropic materials. In anisotropic problems, however, their resulting temperature fields can be distinct. In Ex.~(1.3), when the lateral sides are under free boundary conditions ($\widetilde{q}_N = 0$), the spatial temperature distribution at $t = 0.8$ achieved by PG-FPM-3 is presented in Fig.~\ref{fig:Ex13_Free}. As can be seen, the temperature field is no longer independent of $x$. The result is also verified by the Galerkin FPM and other approaches.

\begin{table}[htbp]
\caption{Relative errors and computational time of the FPM and PG-FPMs in solving Ex.~(1.3).}
\centering
{
\begin{tabular*}{500pt}{@{\extracolsep\fill}lcccc@{\extracolsep\fill}}
\toprule
\textbf{Method} & \tabincell{c}{\textbf{Computational} \\ \textbf{parameters}} &  \textbf{Relative errors} & \tabincell{c}{$N_{band} (\mathbf{K})$ \& $N_{band} (\mathbf{C})$}  & \tabincell{c}{\textbf{Computational} \\ \textbf{time (s)}} \\
\midrule
\multicolumn{5}{c}{Homogenous isotropic ($\delta = 0$; $\hat{k}_{11} = \hat{k}_{22} = 1, \hat{k}_{12} = \hat{k}_{21} = 0$)} \\
\midrule
FPM & $\eta_1 = 1$, $\eta_2 = 1 \times 10^{5}$  & $\overline{e}_0 = 5.9 \times 10^{-3}$ & 12 \& 1 & 2.0 \\
PG-FPM-1 & $\eta_1 = 1$, $\eta_2 = 1 \times 10^{5}$  & $\overline{e}_0 = 6.3 \times 10^{-3}$ & 13 \& 1 & 2.0 \\
PG-FPM-2 & $\eta_1 = 1$, $\eta_2 = 1 \times 10^{5}$  & $\overline{e}_0 = 5.0 \times 10^{-3}$ & 8 \& 1 & 1.4 \\
PG-FPM-3 & $\eta_1 = 1$, $\eta_2 = 1 \times 10^{5}$  & $\overline{e}_0 = 5.6 \times 10^{-3}$ & 8 \& 1 & 1.5 \\
\midrule
\multicolumn{5}{c}{Nonhomogenous isotropic ($\delta = 3$; $\hat{k}_{11} = \hat{k}_{22} = 1, \hat{k}_{12} = \hat{k}_{21} = 0$)} \\
\midrule
FPM & $\eta_1 = 1$, $\eta_2 = 1 \times 10^{5}$  & $\overline{e}_0 = 2.8 \times 10^{-2}$ & 12 \& 1 & 2.0 \\
PG-FPM-1 & $\eta_1 = 1$, $\eta_2 = 1 \times 10^{5}$  & $\overline{e}_0 = 6.6 \times 10^{-3}$ & 13 \& 1 & 2.0 \\
PG-FPM-2 & $\eta_1 = 1$, $\eta_2 = 1 \times 10^{5}$  & $\overline{e}_0 = 1.2 \times 10^{-2}$ & 8 \& 1 & 1.4 \\
PG-FPM-3 & $\eta_1 = 1$, $\eta_2 = 1 \times 10^{5}$  & $\overline{e}_0 = 9.5 \times 10^{-3}$ & 8 \& 1 & 1.6 \\
\midrule
\multicolumn{5}{c}{Nonhomogenous anisotropic ($\delta = 3$; $\hat{k}_{11} = \hat{k}_{22} = 2, \hat{k}_{12} = \hat{k}_{21} = 1$)} \\
\midrule
FPM & $\eta_1 = 1$, $\eta_2 = 1 \times 10^{5}$  & $\overline{e}_0 = 3.1 \times 10^{-2}$ & 20 \& 1 & 2.2 \\
PG-FPM-1 & $\eta_1 = 1$, $\eta_2 = 1 \times 10^{5}$  & $\overline{e}_0 = 8.2 \times 10^{-3}$ & 13 \& 1 & 2.2 \\
PG-FPM-2 & $\eta_1 = 1$, $\eta_2 = 1 \times 10^{5}$  & $\overline{e}_0 = 1.4 \times 10^{-2}$ & 11 \& 1 & 1.4 \\
PG-FPM-3 & $\eta_1 = 1$, $\eta_2 = 1 \times 10^{5}$  & $\overline{e}_0 = 3.3 \times 10^{-2}$ & 11 \& 1 & 1.4 \\
\bottomrule
\end{tabular*}}
\label{table:Ex13}
\end{table}

\begin{figure}[htbp] 
  \centering 
    \subfigure[]{ 
    \label{fig:Ex14_Trans} 
    \includegraphics[width=0.48\textwidth]{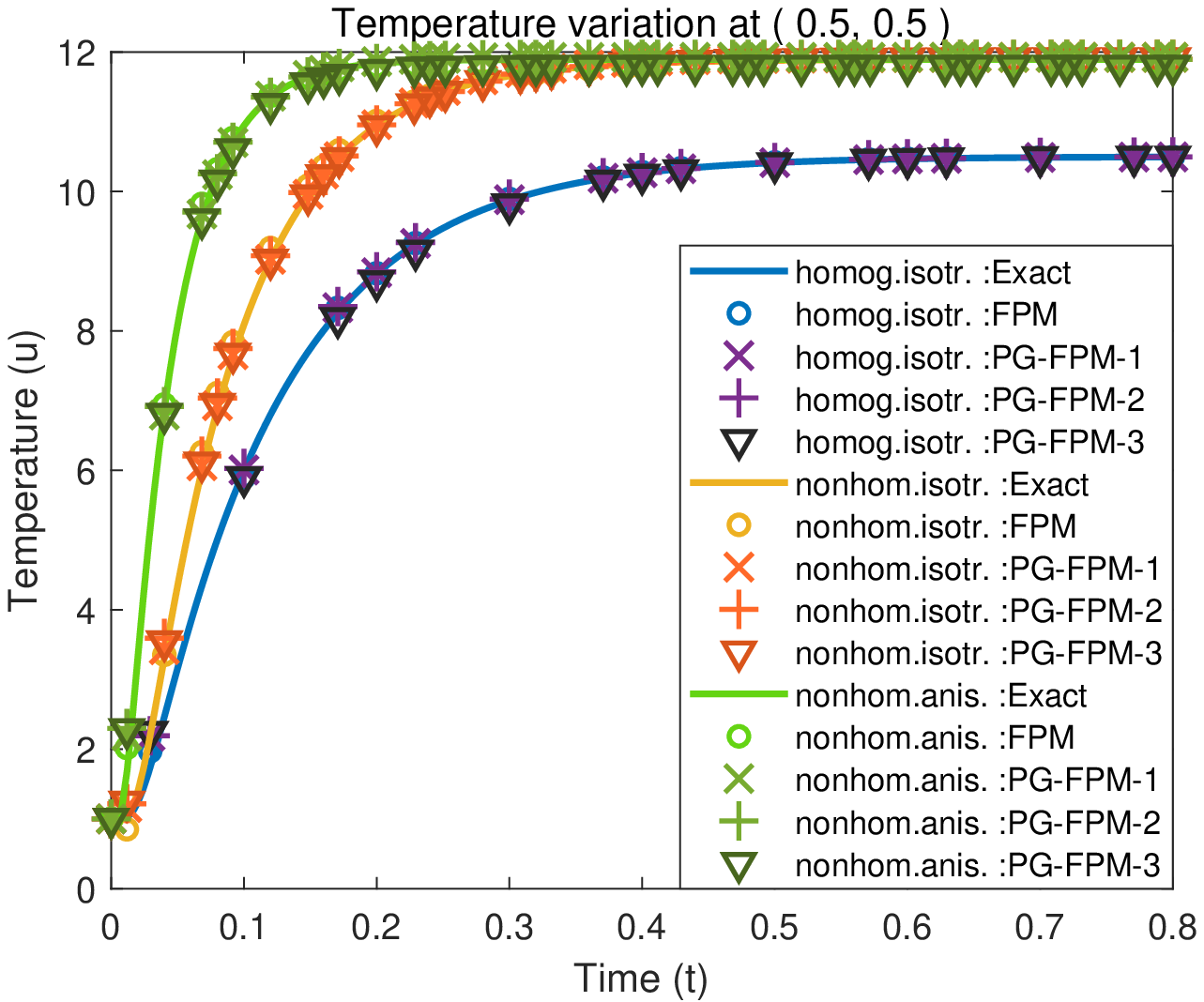}}  
    \subfigure[]{ 
    \label{fig:Ex14_Conf} 
    \includegraphics[width=0.48\textwidth]{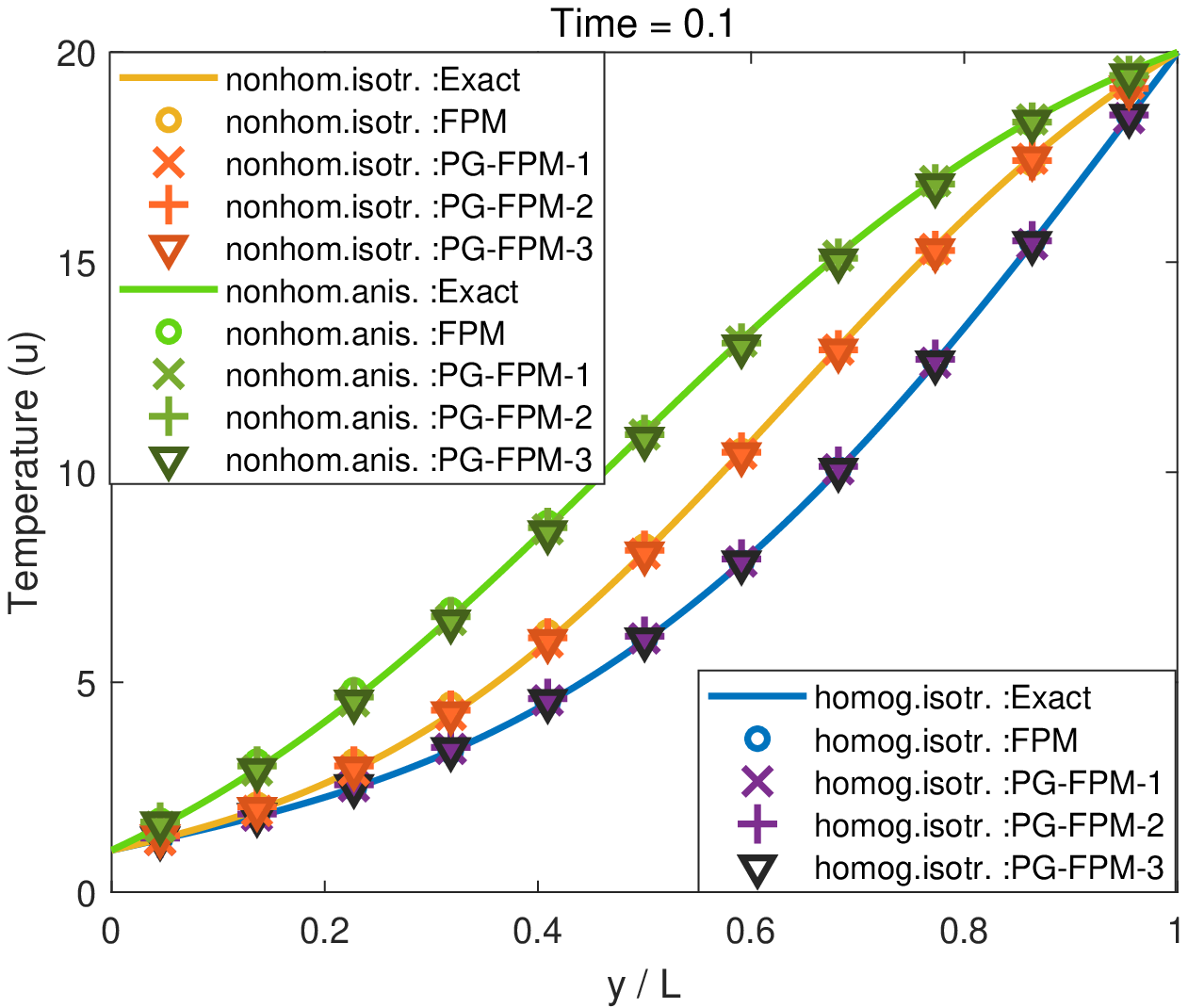}}  
  \caption{Ex. (1.4) - The computed solution achieved by the FPM and PG-FPMs with different material properties. (a) transient temperature solution at the midpoint of the domain in time scope $[0, 0.8]$. (b) vertical temperature distribution when $t=0.1$.} 
  \label{fig:Ex14} 
\end{figure}

\begin{figure}[htbp] 
  \centering 
    \subfigure[]{ 
    \label{fig:Ex15_Trans} 
    \includegraphics[width=0.48\textwidth]{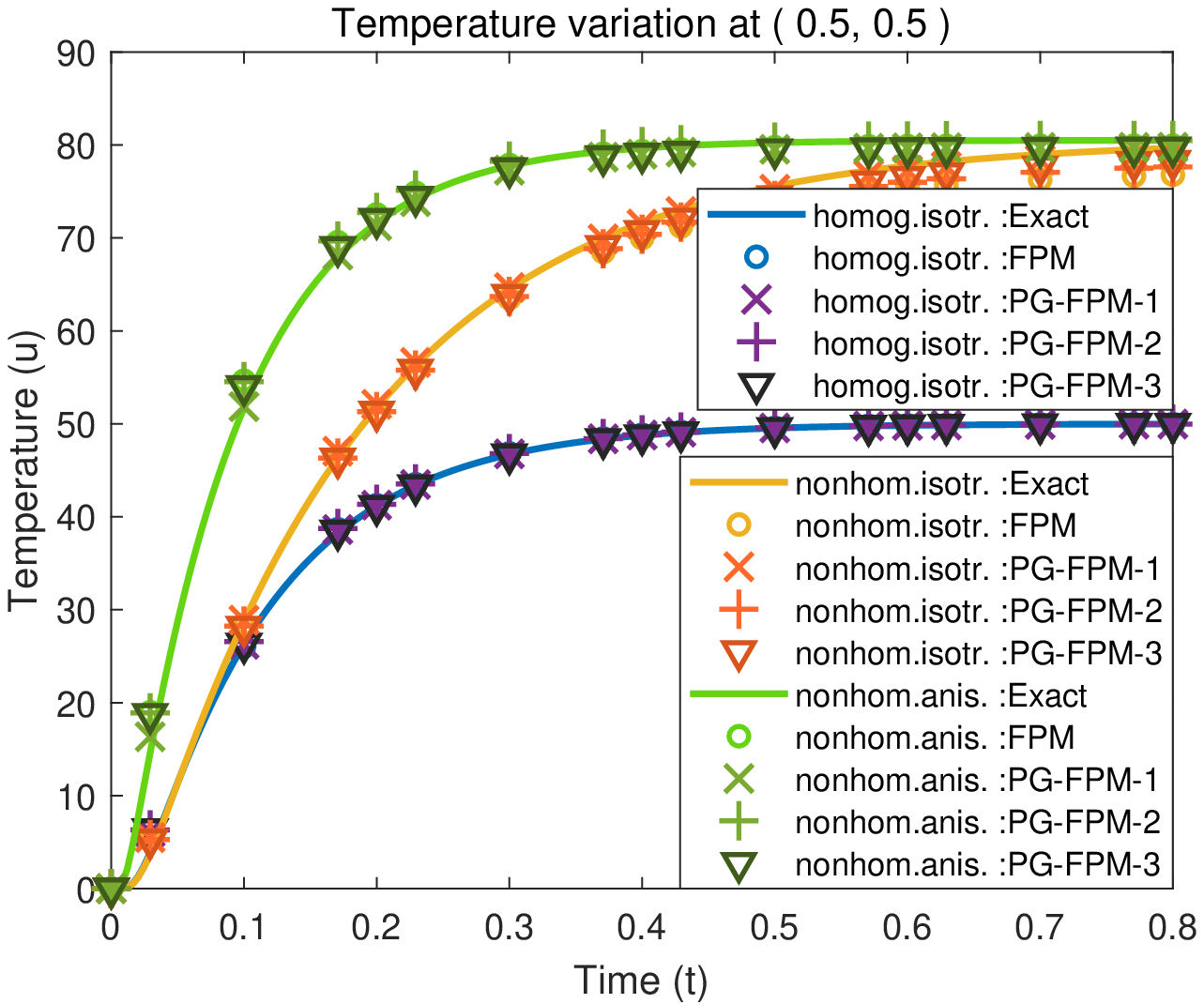}}  
    \subfigure[]{ 
    \label{fig:Ex15_Conf} 
    \includegraphics[width=0.48\textwidth]{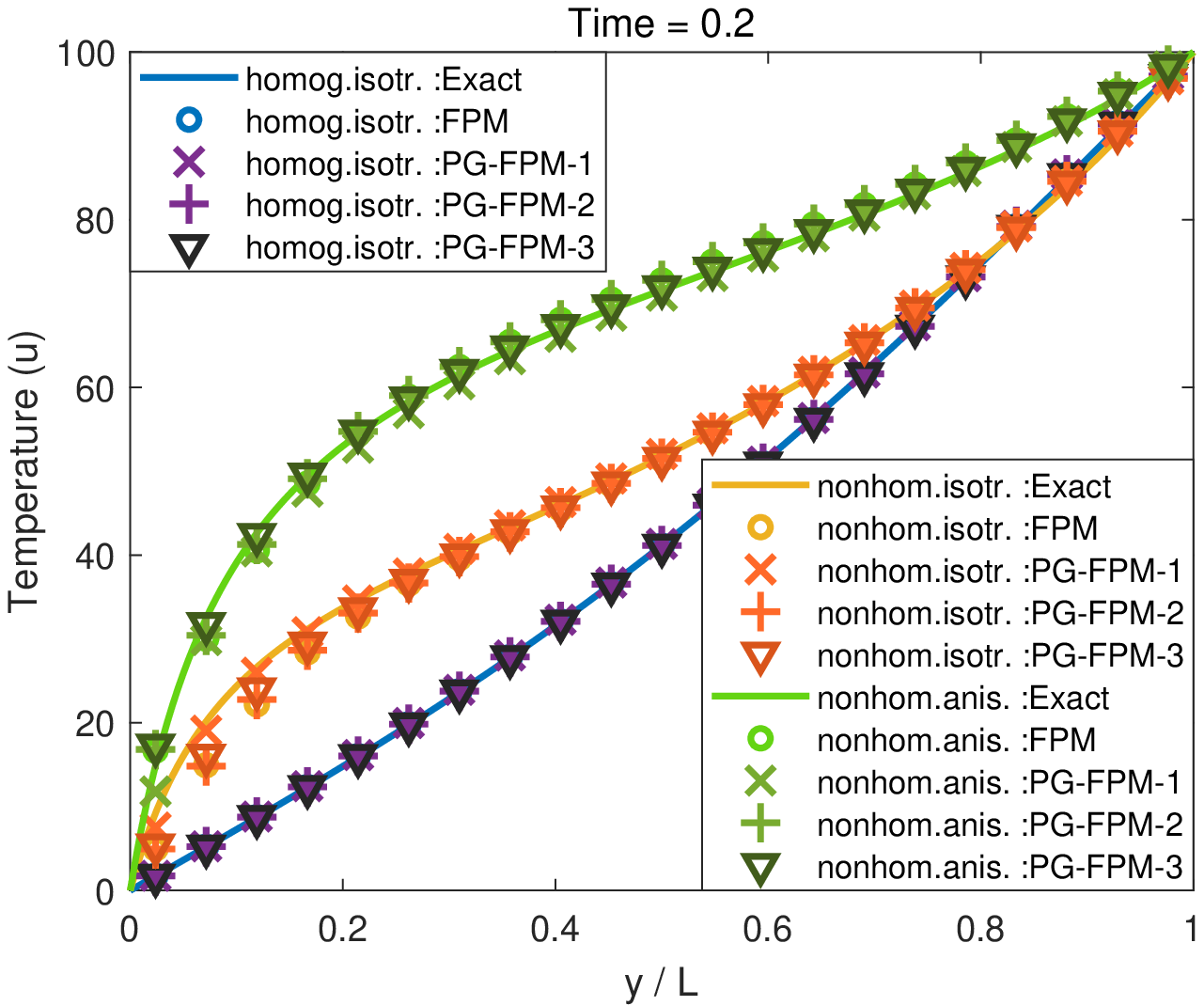}}  
  \caption{Ex. (1.5) - The computed solution achieved by the FPM and PG-FPMs with different material properties. (a) transient temperature solution at the midpoint of the domain in time scope $[0, 0.8]$. (b) vertical temperature distribution when $t=0.2$.} 
  \label{fig:Ex15} 
\end{figure}

\begin{figure}[htbp] 
  \centering 
    \subfigure[]{ 
    \label{fig:Ex16_Trans} 
    \includegraphics[width=0.48\textwidth]{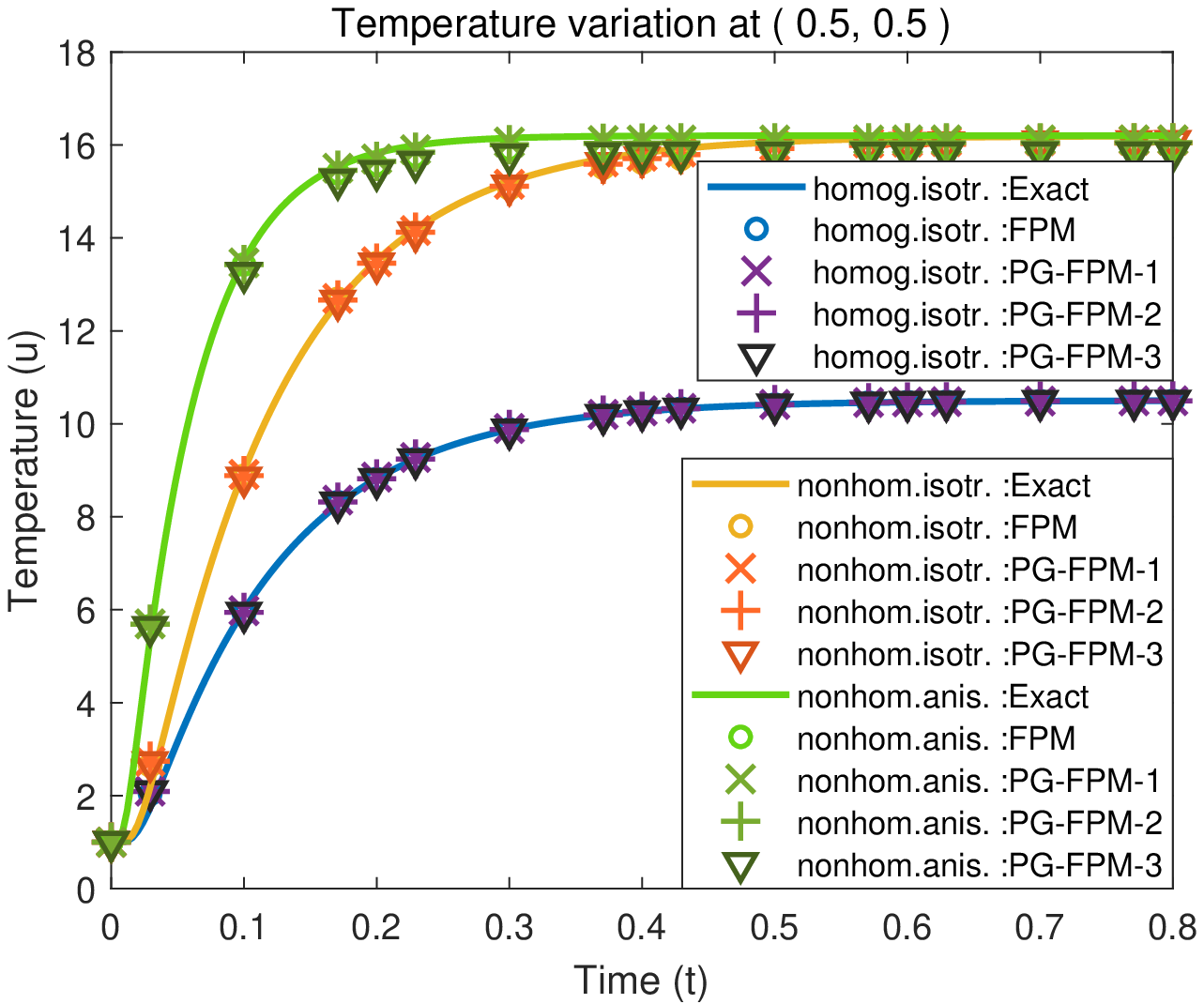}}  
    \subfigure[]{ 
    \label{fig:Ex16_Conf} 
    \includegraphics[width=0.48\textwidth]{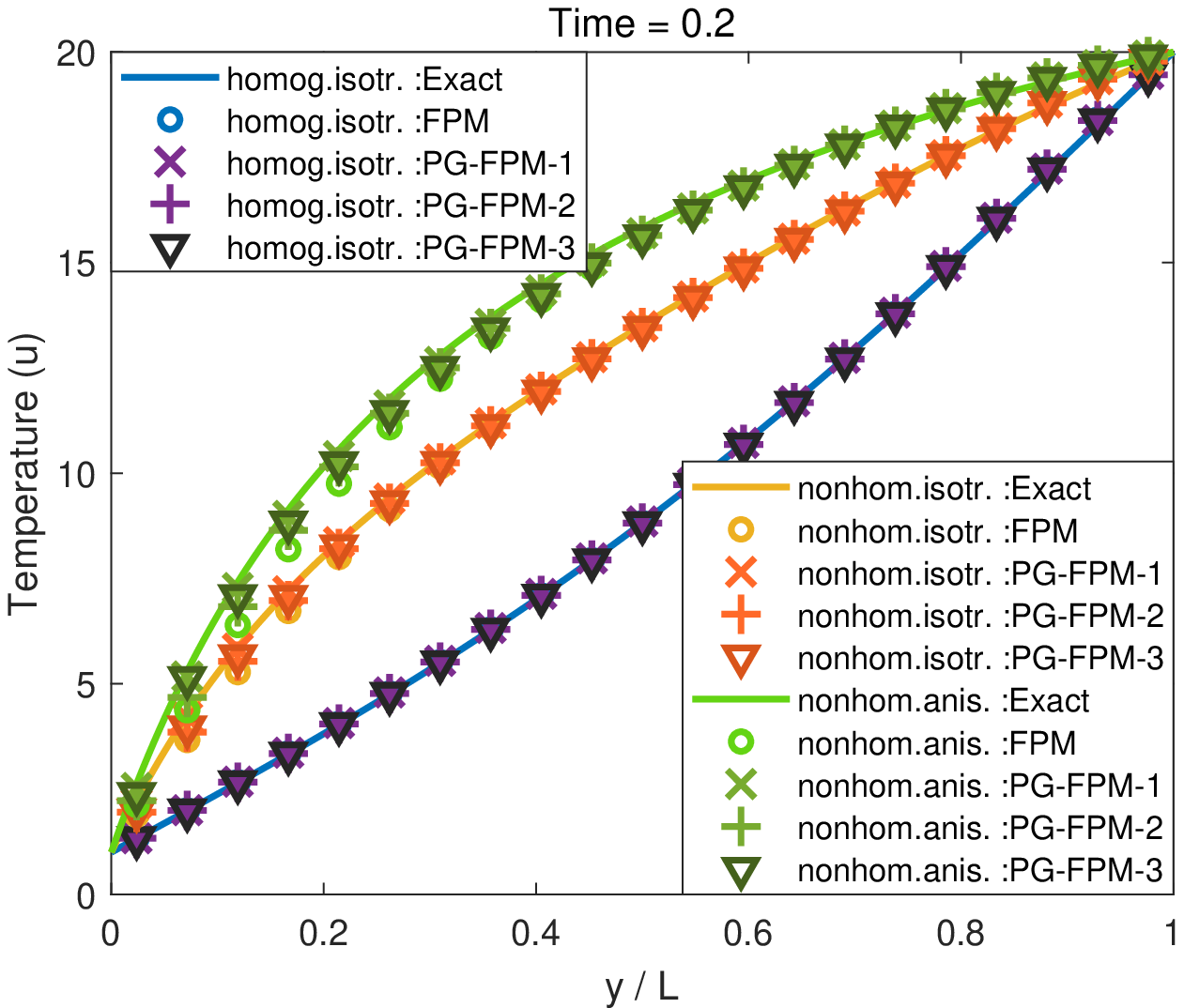}}  
  \caption{Ex. (1.6) - The computed solution achieved by the FPM and PG-FPMs with different material properties. (a) transient temperature solution at the midpoint of the domain in time scope $[0, 0.8]$. (b) vertical temperature distribution when $t=0.2$.} 
  \label{fig:Ex16} 
\end{figure}

The computed solutions of Ex.~(1.4) – (1.6) obtained by the Galerkin FPM and PG-FPMs are presented in Fig.~\ref{fig:Ex14} – \ref{fig:Ex16}. In all these examples, both the transient and spatial temperature distributions under different material properties show excellent agreement with the exact solution. The corresponding average errors and time costs are listed in Table~\ref{table:Ex14} – \ref{table:Ex16}. The time step $\Delta t = 0.1$. In Ex.~(1.4), with 121 uniformly distributed Points, the average errors for all the proposed PG-FPM approaches are less than 1\%. The collocation method (PG-FPM-1) shows the best accuracy in all the examples. For 2D problems with a small number of Fragile Points, the collocation method may cost a slightly longer time than the Galerkin FPM, as a result of the complicated local approximations. However, when the number of Points increases, the efficiency of the collocation method exceeds the Galerkin FPM, especially when the penalty parameter $\eta_1 = 0$. For instance, in Ex.~(1.5) and (1.6), with 441 Fragile Points employed, the collocation method presents better efficiency and accuracy than the Galerkin FPM. The finite volume method (PG-FPM-2), on the other hand, saves approximately 40\% of the computational time and maintains similar accuracy as compared to the Galerkin FPM. At last, for the singular solution method (PG-FPM-3), in Ex.~(1.4) and (1.5), when only one integration point is adopted in each subdomain, the method shows similar efficiency as the finite volume method. However, in Ex.~(1.6), with two integration points in each subdomain, the efficiency of the PG-FPM-3 decreases dramatically. In practical transient problems, it can be a challenge to determine the appropriate number of integration points. Therefore, the singular solution method is recommended for steady-state analysis only.

\begin{table}[htbp]
\caption{Relative errors and computational times of the FPM and PG-FPMs in solving Ex.~(1.4).}
\centering
{
\begin{tabular*}{500pt}{@{\extracolsep\fill}lcccc@{\extracolsep\fill}}
\toprule
\textbf{Method} & \tabincell{c}{\textbf{Computational} \\ \textbf{parameters}} &  \textbf{Relative errors} & \tabincell{c}{$N_{band} (\mathbf{K})$ \& $N_{band} (\mathbf{C})$}  & \tabincell{c}{\textbf{Computational} \\ \textbf{time (s)}} \\
\midrule
\multicolumn{5}{c}{Homogenous isotropic ($\delta = 0$; $\hat{k}_{11} = \hat{k}_{22} = 1, \hat{k}_{12} = \hat{k}_{21} = 0$)} \\
\midrule
FPM & $\eta_1 = 1$, $\eta_2 = 1 \times 10^{5}$  & $\overline{e}_0 = 5.0 \times 10^{-3}$ & 11 \& 1 & 0.42 \\
PG-FPM-1 & $\eta_1 = 1$, $\eta_2 = 1 \times 10^{5}$  & $\overline{e}_0 = 5.4 \times 10^{-3}$ & 13 \& 1 & 0.46 \\
PG-FPM-2 & $\eta_1 = 1$, $\eta_2 = 1 \times 10^{5}$  & $\overline{e}_0 = 7.1 \times 10^{-3}$ & 8 \& 1 & 0.24 \\
PG-FPM-3 & $\eta_1 = 1$, $\eta_2 = 1 \times 10^{5}$  & $\overline{e}_0 = 5.5 \times 10^{-3}$ & 8 \& 1 & 0.25 \\
\midrule
\multicolumn{5}{c}{Nonhomogenous isotropic ($\delta = 2$; $\hat{k}_{11} = \hat{k}_{22} = 1, \hat{k}_{12} = \hat{k}_{21} = 0$)} \\
\midrule
FPM & $\eta_1 = 1$, $\eta_2 = 1 \times 10^{5}$  & $\overline{e}_0 = 1.3 \times 10^{-2}$ & 11 \& 1 & 0.40 \\
PG-FPM-1 & $\eta_1 = 1$, $\eta_2 = 1 \times 10^{5}$  & $\overline{e}_0 = 6.9 \times 10^{-3}$ & 13 \& 1 & 0.47 \\
PG-FPM-2 & $\eta_1 = 1$, $\eta_2 = 1 \times 10^{5}$  & $\overline{e}_0 = 8.7 \times 10^{-3}$ & 8 \& 1 & 0.24 \\
PG-FPM-3 & $\eta_1 = 1$, $\eta_2 = 1 \times 10^{5}$  & $\overline{e}_0 = 7.5 \times 10^{-3}$ & 8 \& 1 & 0.28 \\
\midrule
\multicolumn{5}{c}{Nonhomogenous anisotropic ($\delta = 2$; $\hat{k}_{11} = \hat{k}_{22} = 2, \hat{k}_{12} = \hat{k}_{21} = 1$)} \\
\midrule
FPM & $\eta_1 = 1$, $\eta_2 = 1 \times 10^{5}$  & $\overline{e}_0 = 1.4 \times 10^{-2}$ & 20 \& 1 & 0.42 \\
PG-FPM-1 & $\eta_1 = 1$, $\eta_2 = 1 \times 10^{5}$  & $\overline{e}_0 = 8.6 \times 10^{-3}$ & 13 \& 1 & 0.44 \\
PG-FPM-2 & $\eta_1 = 1$, $\eta_2 = 1 \times 10^{5}$  & $\overline{e}_0 = 9.3 \times 10^{-3}$ & 11 \& 1 & 0.26 \\
PG-FPM-3 & $\eta_1 = 1$, $\eta_2 = 1 \times 10^{5}$  & $\overline{e}_0 = 9.8 \times 10^{-3}$ & 11 \& 5 & 0.27 \\
\bottomrule
\end{tabular*}}
\label{table:Ex14}
\end{table}

\begin{table}[htbp]
\caption{Relative errors and computational times of the FPM and PG-FPMs in solving Ex.~(1.5).}
\centering
{
\begin{tabular*}{500pt}{@{\extracolsep\fill}lcccc@{\extracolsep\fill}}
\toprule
\textbf{Method} & \tabincell{c}{\textbf{Computational} \\ \textbf{parameters}} &  \textbf{Relative errors} & \tabincell{c}{$N_{band} (\mathbf{K})$ \& $N_{band} (\mathbf{C})$}  & \tabincell{c}{\textbf{Computational} \\ \textbf{time (s)}} \\
\midrule
\multicolumn{5}{c}{Homogenous isotropic ($\delta = 0$; $\hat{k}_{11} = \hat{k}_{22} = 1, \hat{k}_{12} = \hat{k}_{21} = 0$)} \\
\midrule
FPM & $\eta_1 = 1$, $\eta_2 = 1 \times 10^{5}$  & $\overline{e}_0 = 6.4 \times 10^{-3}$ & 12 \& 1 & 2.6 \\
PG-FPM-1 & $\eta_1 = 0$, $\eta_2 = 1 \times 10^{5}$  & $\overline{e}_0 = 6.8 \times 10^{-3}$ & 13 \& 1 & 2.3 \\
PG-FPM-2 & $\eta_1 = 1$, $\eta_2 = 1 \times 10^{5}$  & $\overline{e}_0 = 5.4 \times 10^{-3}$ & 8 \& 1 & 1.6 \\
PG-FPM-3 & $\eta_1 = 1$, $\eta_2 = 1 \times 10^{5}$  & $\overline{e}_0 = 6.0 \times 10^{-3}$ & 8 \& 1 & 1.6 \\
\midrule
\multicolumn{5}{c}{Nonhomogenous isotropic ($\delta = 2$; $\hat{k}_{11} = \hat{k}_{22} = 1, \hat{k}_{12} = \hat{k}_{21} = 0$)} \\
\midrule
FPM & $\eta_1 = 1$, $\eta_2 = 1 \times 10^{5}$  & $\overline{e}_0 = 1.2 \times 10^{-2}$ & 12 \& 1 & 2.5 \\
PG-FPM-1 & $\eta_1 = 0$, $\eta_2 = 1 \times 10^{5}$  & $\overline{e}_0 = 4.0 \times 10^{-3}$ & 6 \& 1 & 2.2 \\
PG-FPM-2 & $\eta_1 = 1$, $\eta_2 = 1 \times 10^{5}$  & $\overline{e}_0 = 1.1 \times 10^{-2}$ & 8 \& 1 & 1.6 \\
PG-FPM-3 & $\eta_1 = 1$, $\eta_2 = 1 \times 10^{5}$  & $\overline{e}_0 = 8.9 \times 10^{-3}$ & 8 \& 1 & 1.6 \\
\midrule
\multicolumn{5}{c}{Nonhomogenous anisotropic ($\delta = 2$; $\hat{k}_{11} = \hat{k}_{22} = 2, \hat{k}_{12} = \hat{k}_{21} = 1$)} \\
\midrule
FPM & $\eta_1 = 1$, $\eta_2 = 1 \times 10^{5}$  & $\overline{e}_0 = 1.4 \times 10^{-2}$ & 22 \& 1 & 2.5 \\
PG-FPM-1 & $\eta_1 = 0$, $\eta_2 = 1 \times 10^{5}$  & $\overline{e}_0 = 4.7 \times 10^{-3}$ & 6 \& 1 & 2.1 \\
PG-FPM-2 & $\eta_1 = 1$, $\eta_2 = 1 \times 10^{5}$  & $\overline{e}_0 = 1.3 \times 10^{-2}$ & 12 \& 1 & 1.6 \\
PG-FPM-3 & $\eta_1 = 1$, $\eta_2 = 1 \times 10^{5}$  & $\overline{e}_0 = 1.3 \times 10^{-2}$ & 12 \& 1 & 1.6 \\
\bottomrule
\end{tabular*}}
\label{table:Ex15}
\end{table}

\begin{table}[htbp]
\caption{Relative errors and computational times of the FPM and PG-FPMs in solving Ex.~(1.6).}
\centering
{
\begin{tabular*}{500pt}{@{\extracolsep\fill}lcccc@{\extracolsep\fill}}
\toprule
\textbf{Method} & \tabincell{c}{\textbf{Computational} \\ \textbf{parameters}} &  \textbf{Relative errors} & \tabincell{c}{$N_{band} (\mathbf{K})$ \& $N_{band} (\mathbf{C})$}  & \tabincell{c}{\textbf{Computational} \\ \textbf{time (s)}} \\
\midrule
\multicolumn{5}{c}{Homogenous isotropic ($\delta = 0$; $\hat{k}_{11} = \hat{k}_{22} = 1, \hat{k}_{12} = \hat{k}_{21} = 0$)} \\
\midrule
FPM & $\eta_1 = 1$, $\eta_2 = 1 \times 10^{5}$  & $\overline{e}_0 = 5.9 \times 10^{-3}$ & 12 \& 1 & 2.6 \\
PG-FPM-1 & $\eta_1 = 0$, $\eta_2 = 1 \times 10^{5}$  & $\overline{e}_0 = 6.3 \times 10^{-3}$ & 6 \& 1 & 2.2 \\
PG-FPM-2 & $\eta_1 = 1$, $\eta_2 = 1 \times 10^{5}$  & $\overline{e}_0 = 5.0 \times 10^{-3}$ & 8 \& 1 & 1.6 \\
PG-FPM-3 & $\eta_1 = 1$, $\eta_2 = 1 \times 10^{5}$  & $\overline{e}_0 = 5.6 \times 10^{-3}$ & 8 \& 5 & 5.8 \\
\midrule
\multicolumn{5}{c}{Nonhomogenous isotropic ($\delta = 3$; $\hat{k}_{11} = \hat{k}_{22} = 1, \hat{k}_{12} = \hat{k}_{21} = 0$)} \\
\midrule
FPM & $\eta_1 = 1$, $\eta_2 = 1 \times 10^{5}$  & $\overline{e}_0 = 3.0 \times 10^{-2}$ & 12 \& 1 & 2.5 \\
PG-FPM-1 & $\eta_1 = 0$, $\eta_2 = 1 \times 10^{5}$  & $\overline{e}_0 = 7.5 \times 10^{-3}$ & 6 \& 1 & 2.1 \\
PG-FPM-2 & $\eta_1 = 1$, $\eta_2 = 1 \times 10^{5}$  & $\overline{e}_0 = 1.7 \times 10^{-2}$ & 8 \& 1 & 1.6 \\
PG-FPM-3 & $\eta_1 = 1$, $\eta_2 = 1 \times 10^{5}$  & $\overline{e}_0 = 1.3 \times 10^{-2}$ & 8 \& 5 & 5.9 \\
\midrule
\multicolumn{5}{c}{Nonhomogenous anisotropic ($\delta = 3$; $\hat{k}_{11} = \hat{k}_{22} = 2, \hat{k}_{12} = \hat{k}_{21} = 1$)} \\
\midrule
FPM & $\eta_1 = 1$, $\eta_2 = 1 \times 10^{5}$  & $\overline{e}_0 = 3.3 \times 10^{-2}$ & 22 \& 1 & 2.6 \\
PG-FPM-1 & $\eta_1 = 0$, $\eta_2 = 1 \times 10^{5}$  & $\overline{e}_0 = 8.7 \times 10^{-3}$ & 6 \& 1 & 2.1 \\
PG-FPM-2 & $\eta_1 = 1$, $\eta_2 = 1 \times 10^{5}$  & $\overline{e}_0 = 1.9 \times 10^{-2}$ & 12 \& 1 & 1.6 \\
PG-FPM-3 & $\eta_1 = 1$, $\eta_2 = 1 \times 10^{5}$  & $\overline{e}_0 = 3.3 \times 10^{-2}$ & 12 \& 5 & 6.0 \\
\bottomrule
\end{tabular*}}
\label{table:Ex16}
\end{table}

\subsection{Some practical examples}

\begin{figure}[htbp] 
  \centering 
    \subfigure[]{ 
    \label{fig:Ex17_Partition} 
    \includegraphics[width=0.48\textwidth]{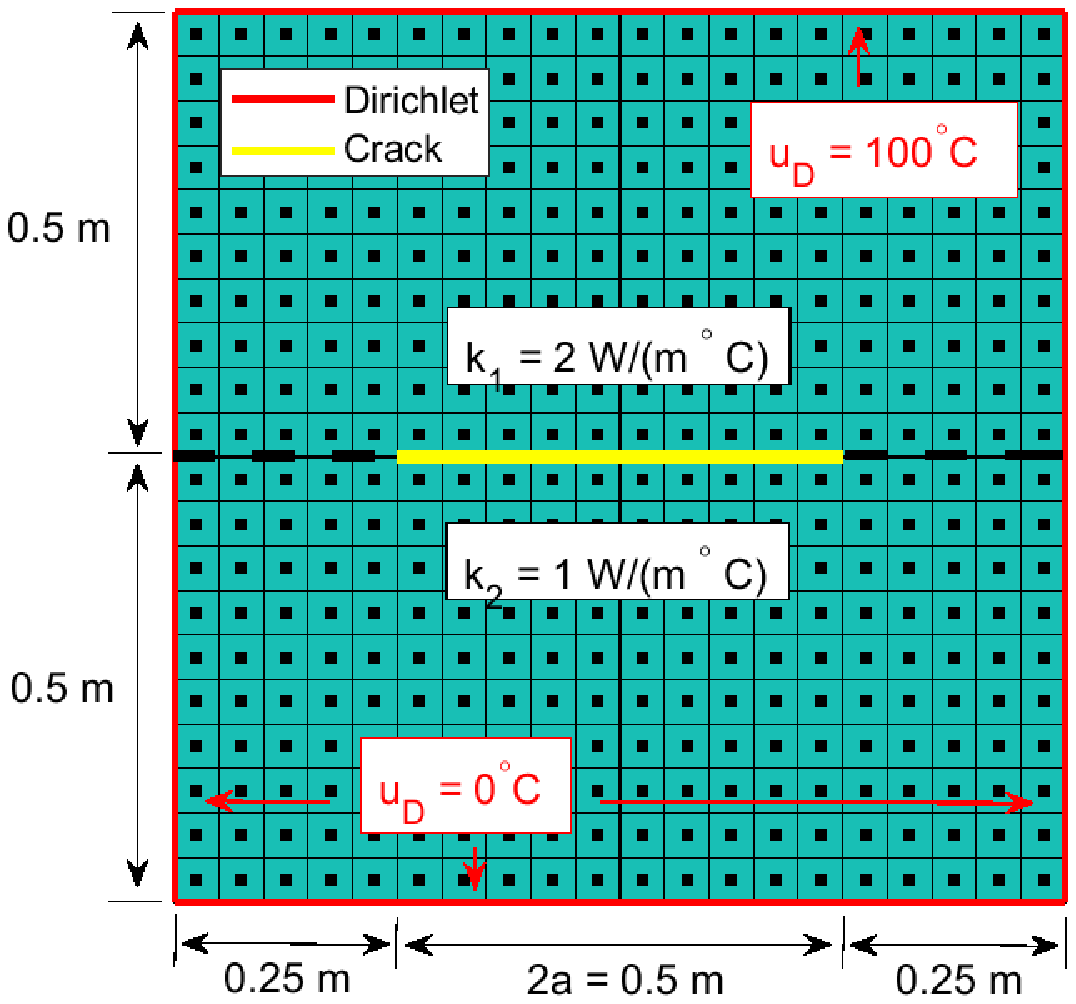}}  
    \subfigure[]{ 
    \label{fig:Ex17_Crack} 
    \includegraphics[width=0.48\textwidth]{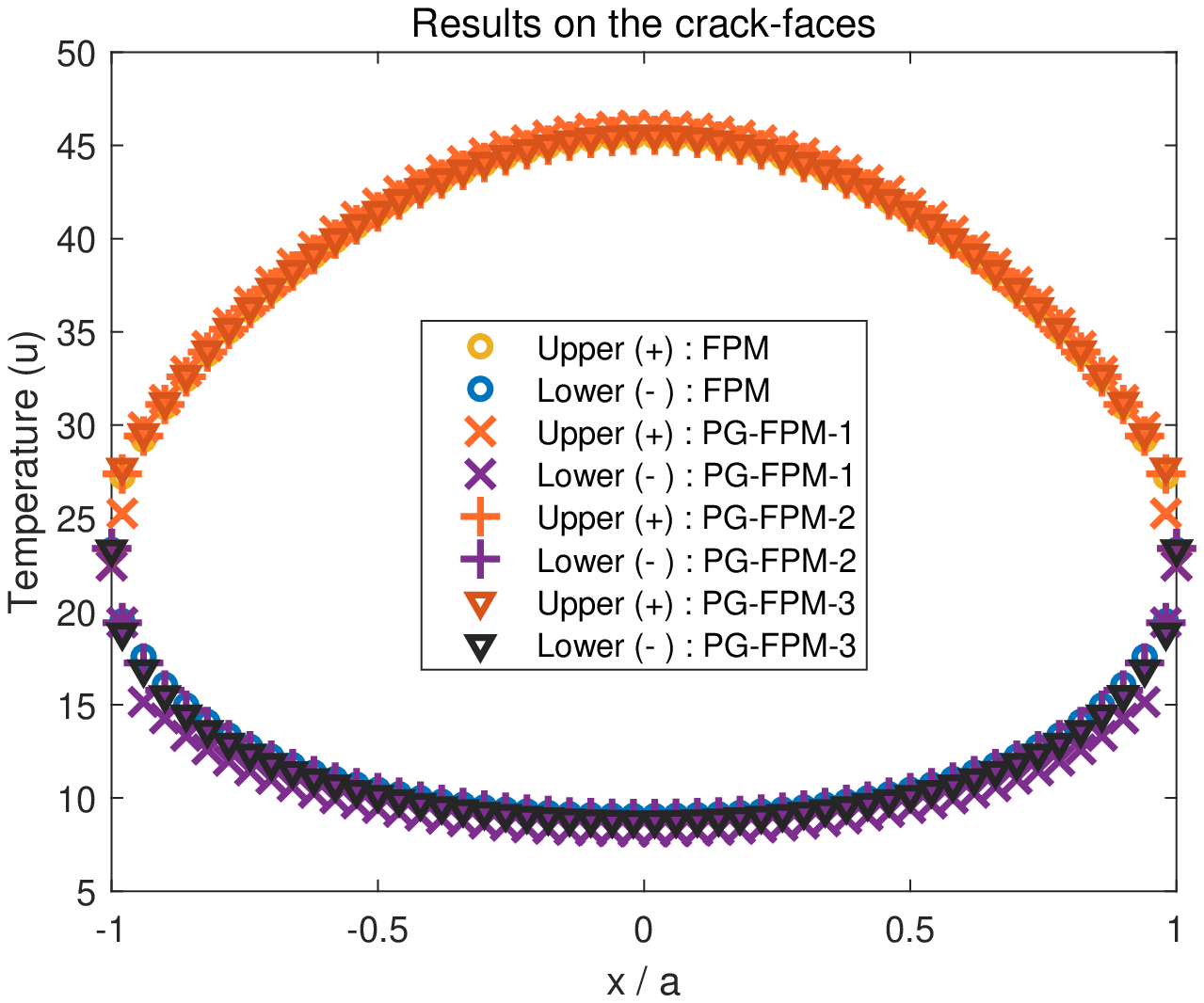}}  
        \subfigure[]{ 
    \label{fig:Ex17_SS_01} 
    \includegraphics[width=0.48\textwidth]{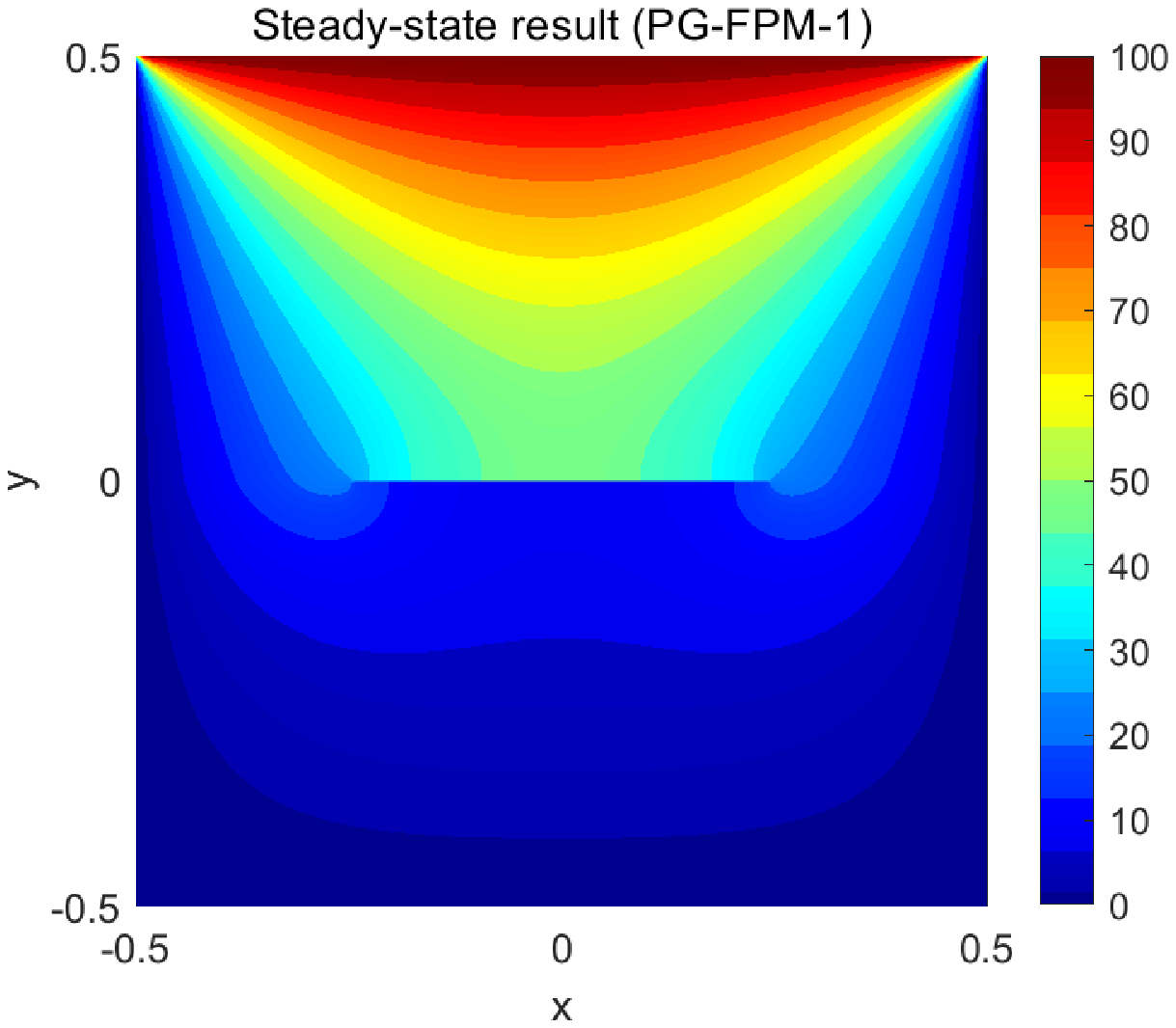}}  
        \subfigure[]{ 
    \label{fig:Ex17_SS_02} 
    \includegraphics[width=0.48\textwidth]{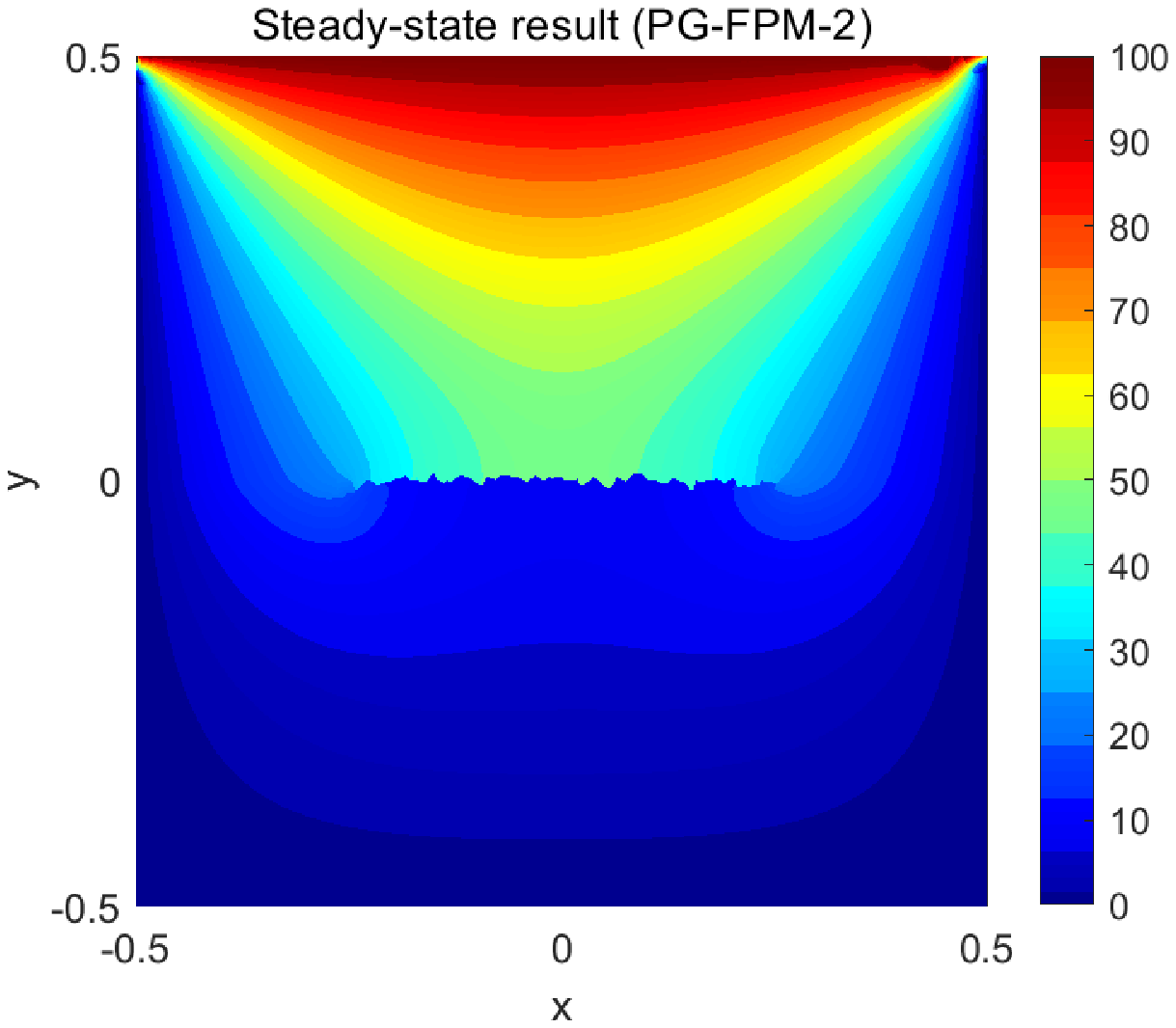}}  
  \caption{Ex. (1.7) – The adiabatic crack and computed solutions. (a) the points, partition and adiabatic crack. (b) the distribution of temperature on the upper (+) and lower (-) crack-faces achieved by the FPM and PG-FPMs. (c) PG-FPM-1 solution. (d) PG-FPM-2 solution.} 
  \label{fig:Ex17} 
\end{figure}

In this section, some practical 2D problems are considered. Ex.~(1.7) is a steady-state heat conduction problem in discontinuous materials with an adiabatic crack. As shown in Fig.~\ref{fig:Ex17_Partition}, the problem domain is a $1~\mathrm{m} \times 1~\mathrm{m}$ square. In the top half of the domain ($ y > 0.5~\mathrm{m}$), the medium is isotropic with thermal conductivity $ k_1 = 2 \mathrm{W / (m ^\circ C)}$, while in the bottom half ($ y < 0.5~\mathrm{m}$), the isotropic thermal conductivity is $ k_2 = 1 \mathrm{W / (m ^\circ C)}$. The boundary conditions on the sides and the adiabatic crack are given as:
\begin{align}
\begin{split}
& \widetilde{u}_D (x, L) = 100~^\circ \mathrm{C}, \qquad  \widetilde{u}_D (x, 0) = \widetilde{u}_D (0, y) = \widetilde{u}_D (L, y) = 0~^\circ \mathrm{C}, \\
& \widetilde{q}_N (x, y) = 0, \qquad \qquad \; \; \text{on} \; -a < x < a, \; y = L/2,
\end{split}
\end{align}
where $L = 1~\mathrm{m}$, $a = 0.25~\mathrm{m}$. The body source density $Q = 0$.

A total of 10000 Points are distributed uniformly or randomly in the domain. In the FPM and PG-FPMs, internal subdomain boundaries shared by two neighboring Points on different sides of the stationary crack are ‘broken’. Thus, we cut off the interaction between the two neighboring Points, i.e., removing the Points from the supporting points set of each other, and convert the cracked internal boundary into two external Neumann boundaries. This algorithm can be further extended to crack development analysis in thermally shocked brittle materials. In the current example, the temperature distribution on the upper (+) and lower (-) crack-faces achieved by the FPM and PG-FPMs are presented in Fig.~\ref{fig:Ex17_Crack}. The methods exhibit great consistency with each other under different point distributions and partitions. Figure~\ref{fig:Ex17_SS_01} and \ref{fig:Ex17_SS_02} display the spatial temperature distribution acquired by the collocation method (PG-FPM-1) with uniform points and by the finite volume method (PG-FPM-2) with random points, respectively. As can be seen, the FPM and PG-FPMs can get a good approximation of the temperature distribution in the global domain even if the subdomain boundaries do not coincide with the crack strictly. These results are also qualitatively consistent with ABAQUS and numerical solutions shown in the previous studies \cite{Liu2019, Guan2020}.

The corresponding computational parameters and times of the different methods in solving Ex.~(1.7) are listed in Table~\ref{table:Ex17}. It is clear that all the proposed PG-FPM approaches can improve the efficiency of the original Galerkin FPM. Among them, the finite volume method (PG-FPM-2) saves more than a half of the computational time. This example validates that the superiority of the finite volume and singular solution methods are especially significant for problems with large degrees-of-freedoms (DoFs).

\begin{table}[htbp]
\caption{Computational times of the FPM and PG-FPMs in solving Ex.~(1.7).}
\centering
{
\begin{tabular*}{500pt}{@{\extracolsep\fill}lcccc@{\extracolsep\fill}}
\toprule
\textbf{Method} & \textbf{Point distribution} & \tabincell{c}{\textbf{Computational} \\ \textbf{parameters}} & $N_{band} (\mathbf{K})$  & \tabincell{c}{\textbf{Computational} \\ \textbf{time (s)}} \\
\midrule
\multirow{2}*{FPM} & uniform& $\eta_1 = 1$, $\eta_2 = 1 \times 10^{5}$ & 13 & 3.9 \\
~ & random &$\eta_1 = 1$, $\eta_2 = 1 \times 10^{5}$ & 42 & 4.5 \\
\midrule
\multirow{2}*{PG-FPM-1} & uniform& $\eta_1 = 1$, $\eta_2 = 1 \times 10^{5}$ & 13 & 3.0 \\
~ & random &$\eta_1 = 1$, $\eta_2 = 1 \times 10^{5}$ & 20 & 3.7 \\
\midrule
\multirow{2}*{PG-FPM-2} & uniform& $\eta_1 = 1$, $\eta_2 = 1 \times 10^{5}$ & 9 & 1.7 \\
~ & random &$\eta_1 = 1$, $\eta_2 = 1 \times 10^{5}$ & 20  & 2.3 \\
\midrule
\multirow{2}*{PG-FPM-3} & uniform& $\eta_1 = 1$, $\eta_2 = 1 \times 10^{5}$ & 9 & 1.9 \\
~ & random &$\eta_1 = 1$, $\eta_2 = 1 \times 10^{5}$ & 20 & 2.7 \\
\bottomrule
\end{tabular*}}
\label{table:Ex17}
\end{table}

\begin{figure}[htbp] 
  \centering 
    \subfigure[]{ 
    \label{fig:Ex18_Partition} 
    \includegraphics[width=0.48\textwidth]{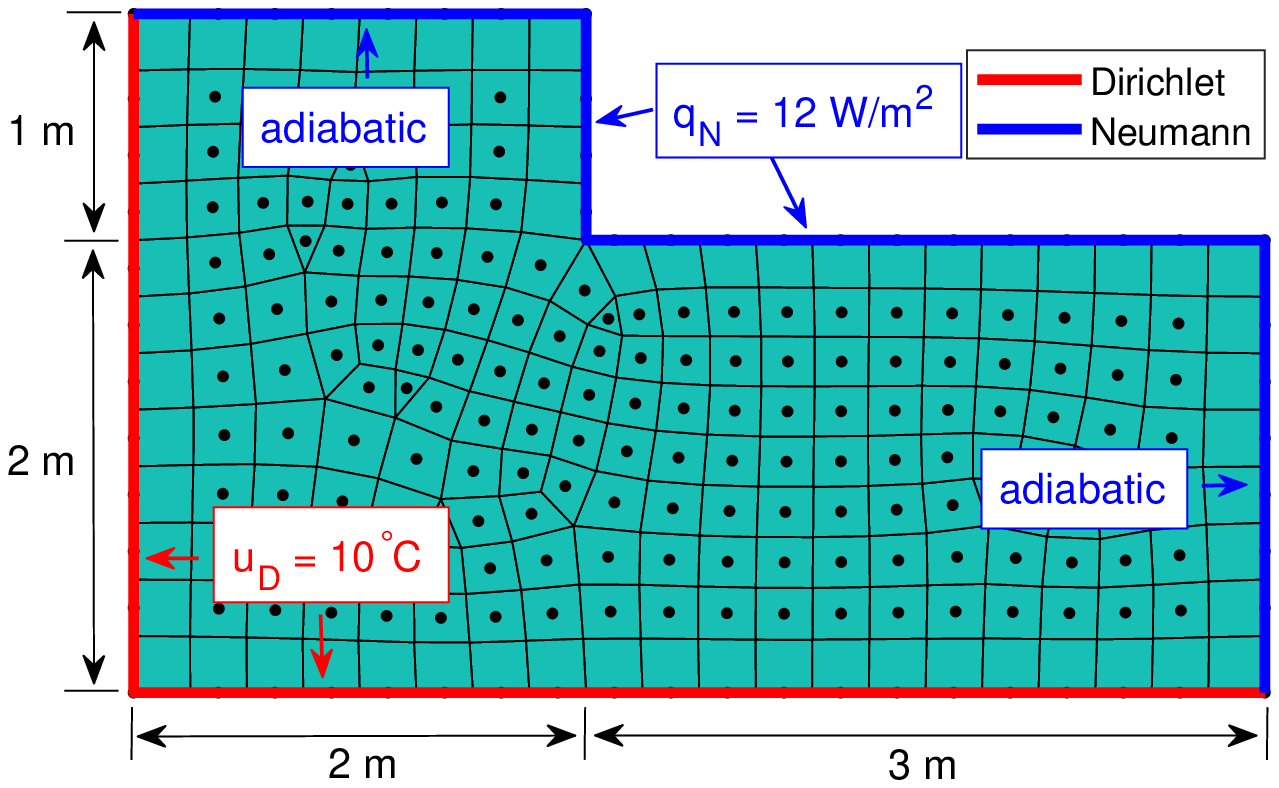}}  
    \subfigure[]{ 
    \label{fig:Ex18_FEA} 
    \includegraphics[width=0.48\textwidth]{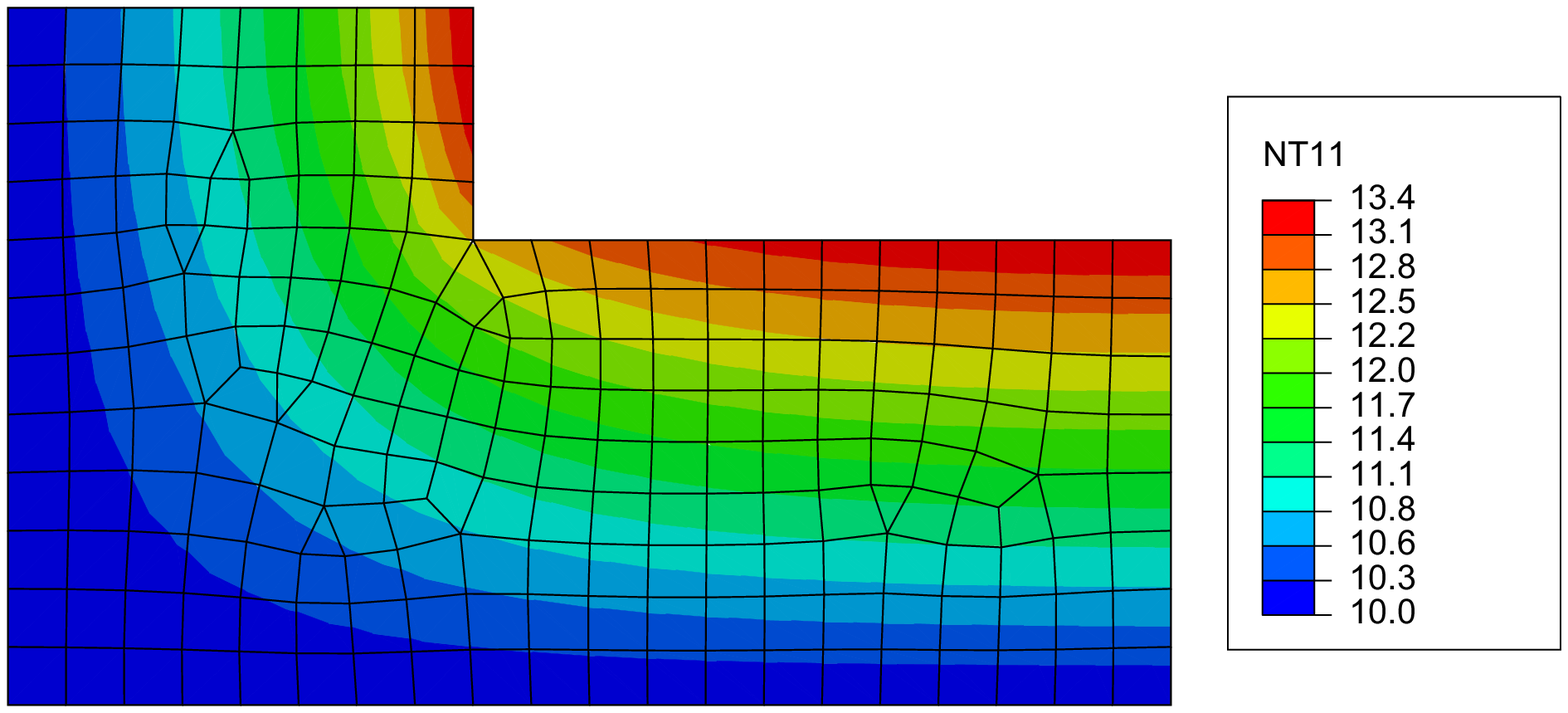}}  
    \subfigure[]{ 
    \label{fig:Ex18_SS_03_abaqus} 
    \includegraphics[width=0.48\textwidth]{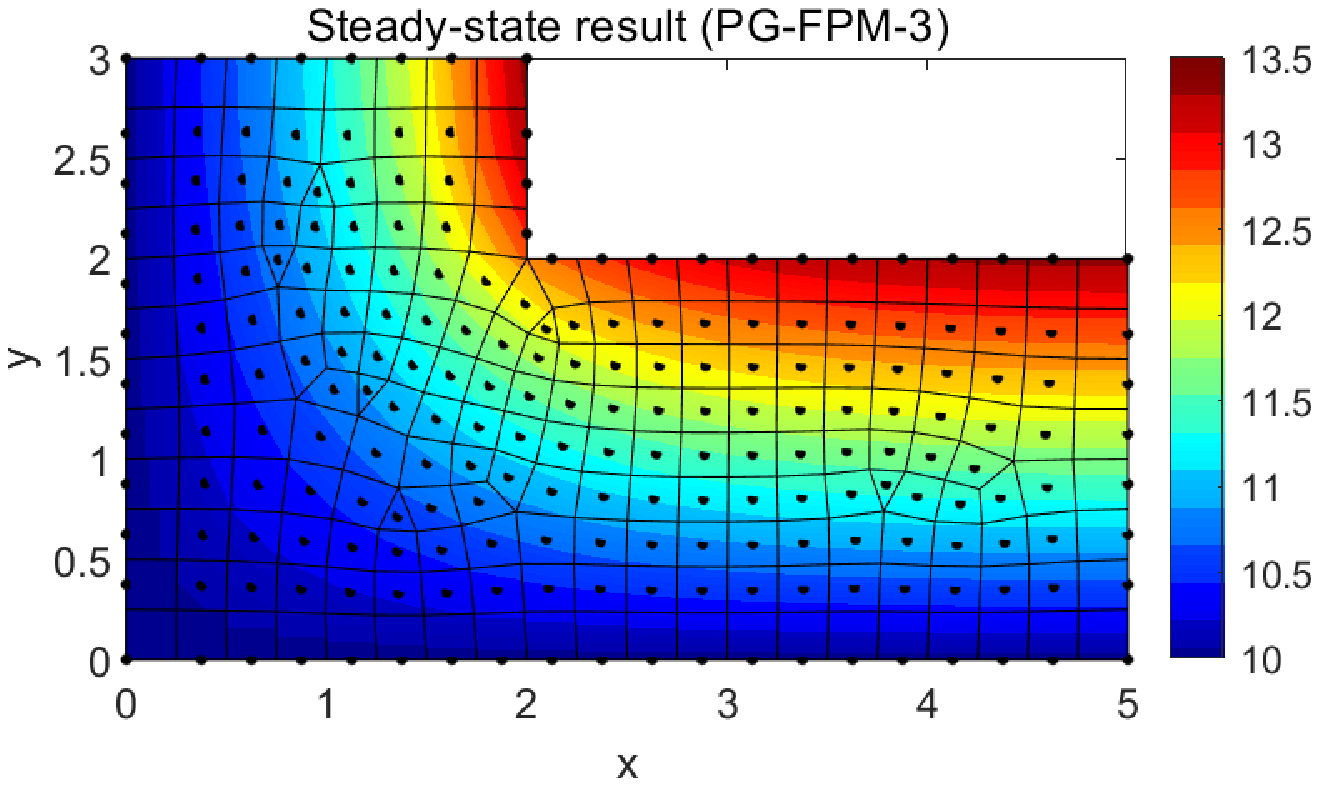}}  
        \subfigure[]{ 
    \label{fig:Ex18_SS_03_uniform} 
    \includegraphics[width=0.48\textwidth]{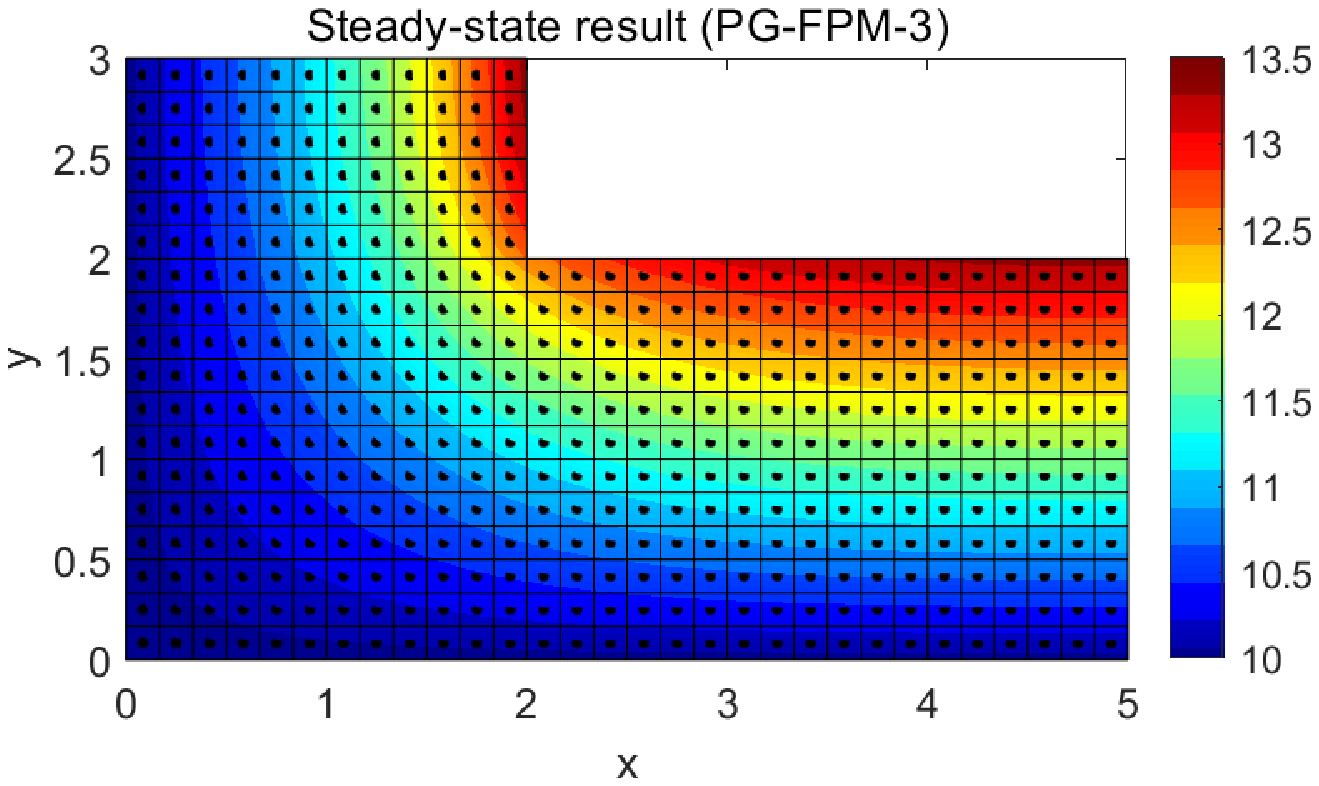}}  
  \caption{Ex. (1.8) – The boundary conditions and computed solutions.  (a) the problem domain and its partition converted from ABAQUS mesh. (b) ABAQUS solution with DC2D3 and DC2D4 elements. (c) PG-FPM-3 solution with ABAQUS partition. (d) PG-FPM-3 solution with uniform point distribution.} 
  \label{fig:Ex18} 
\end{figure}

Ex.~(1.8) is another a steady-state problem in a L-shaped domain. The material is orthotropic with thermal conductivity coefficients $k_{11} = 4~\mathrm{W / m ^\circ C}$, $k_{22} = 7~\mathrm{W / m ^\circ C}$, and $k_{12} = k_{21} = 0$. Dirichlet and Neumann boundary conditions are applied on the external sides. The corresponding values are shown in Fig.~\ref{fig:Ex18_Partition}. In this example, two kinds of domain partitions are considered: a mixed quadrilateral and triangular partition converted from ABAQUS meshes (shown in Fig.~\ref{fig:Ex17_Partition}), and a uniform quadrilateral partition (shown in Fig.~\ref{fig:Ex17_SS_02}). All the Internal Points are placed at the centroid of each subdomain. The computed solution achieved by the FPM and PG-FPMs are approximately the same. Hence, here only the PG-FPM-3 solutions under the ABAQUS and uniform partitions are presented in Fig.~\ref{fig:Ex18_SS_03_abaqus} and \ref{fig:Ex18_SS_03_uniform}. The results also agree well with the FEM solution obtained by ABAQUS using the same element mesh (see Fig.~\ref{fig:Ex18_FEA}). Note that the nodes in the FEM are different from the Fragile Points used in the FPM or PG-FPMs. Table~\ref{table:Ex18} shows the computational times of the FPM and proposed PG-FPMs. In this example, the PG-FPM-1 / 2 /3 approaches elevate the computing efficiency by 14\%, 45\%, and 27\% respectively as compared to the original Galerkin FPM.

\begin{table}[htbp]
\caption{Computational time of the FPM and PG-FPMs in solving Ex.~(1.8).}
\centering
{
\begin{tabular*}{500pt}{@{\extracolsep\fill}lcccc@{\extracolsep\fill}}
\toprule
\textbf{Method} & \tabincell{c}{\textbf{Point distribution} \\ \textbf{(Number of points)}} & \tabincell{c}{\textbf{Computational} \\ \textbf{parameters}} & $N_{band} (\mathbf{K})$  & \tabincell{c}{\textbf{Computational} \\ \textbf{time (s)}} \\
\midrule
\multirow{2}*{FPM} & ABAQUS (218) & $\eta_1 = 1$, $\eta_2 = 1 \times 10^{5}$ & 21 & 0.08 \\
~ & uniform (432) &$\eta_1 = 1$, $\eta_2 = 1 \times 10^{5}$ & 12 & 0.14 \\
\midrule
\multirow{2}*{PG-FPM-1} & ABAQUS (218) & $\eta_1 =1$, $\eta_2 = 1 \times 10^{5}$ & 14 & 0.07 \\
~ & uniform (432) &$\eta_1 = 1$, $\eta_2 = 1 \times 10^{5}$ & 13 & 0.12 \\
\midrule
\multirow{2}*{PG-FPM-2} & ABAQUS (218) & $\eta_1 = 1$, $\eta_2 = 1 \times 10^{5}$ & 12 & 0.04 \\
~ & uniform (432) &$\eta_1 = 1$, $\eta_2 = 1 \times 10^{5}$ & 8  & 0.08 \\
\midrule
\multirow{2}*{PG-FPM-3} & ABAQUS (218) & $\eta_1 = 1$, $\eta_2 = 1 \times 10^{5}$ & 12 & 0.06 \\
~ & uniform (432) &$\eta_1 = 1$, $\eta_2 = 1 \times 10^{5}$ & 8 & 0.10 \\
\bottomrule
\end{tabular*}}
\label{table:Ex18}
\end{table}

\begin{figure}[htbp] 
  \centering 
    \subfigure[]{ 
    \label{fig:Ex19_Partition} 
    \includegraphics[width=0.34\textwidth]{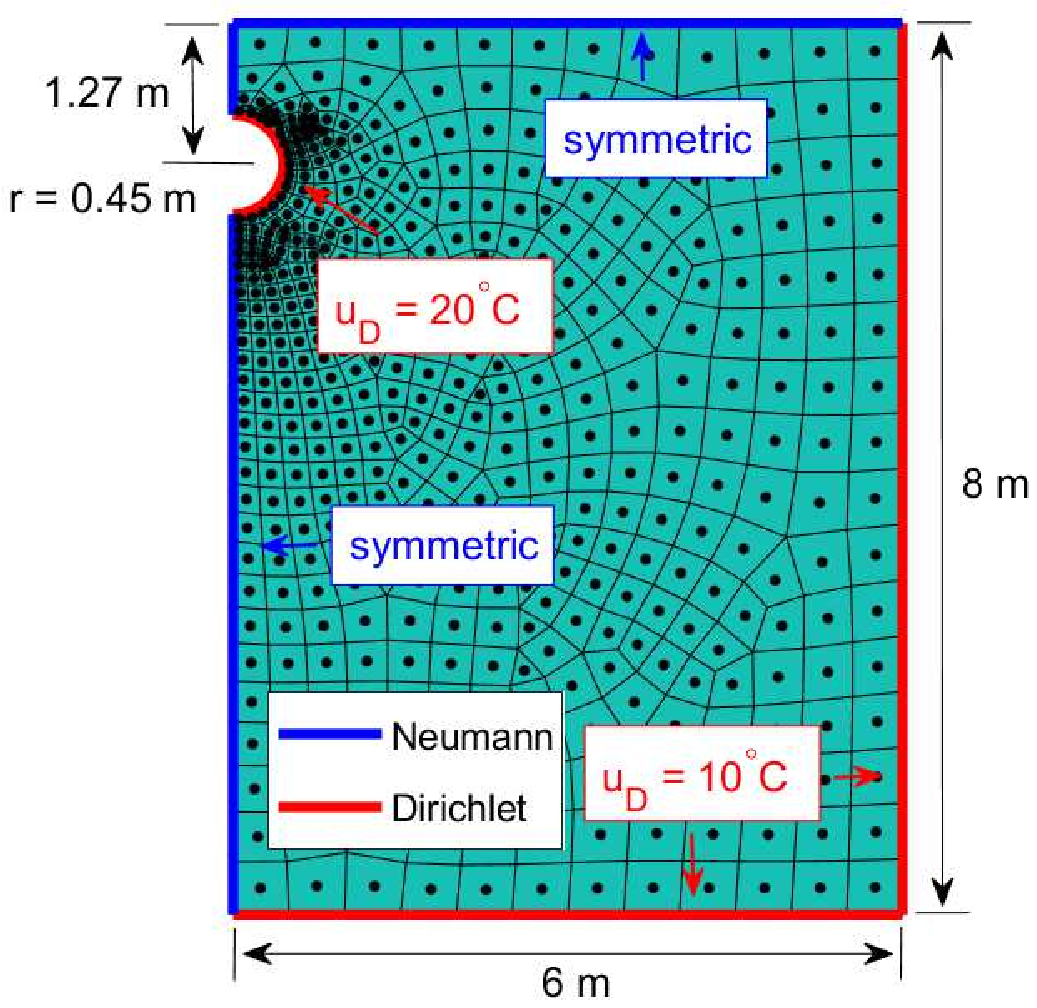}}  
    \subfigure[]{ 
    \label{fig:Ex19_Trans} 
    \includegraphics[width=0.62\textwidth]{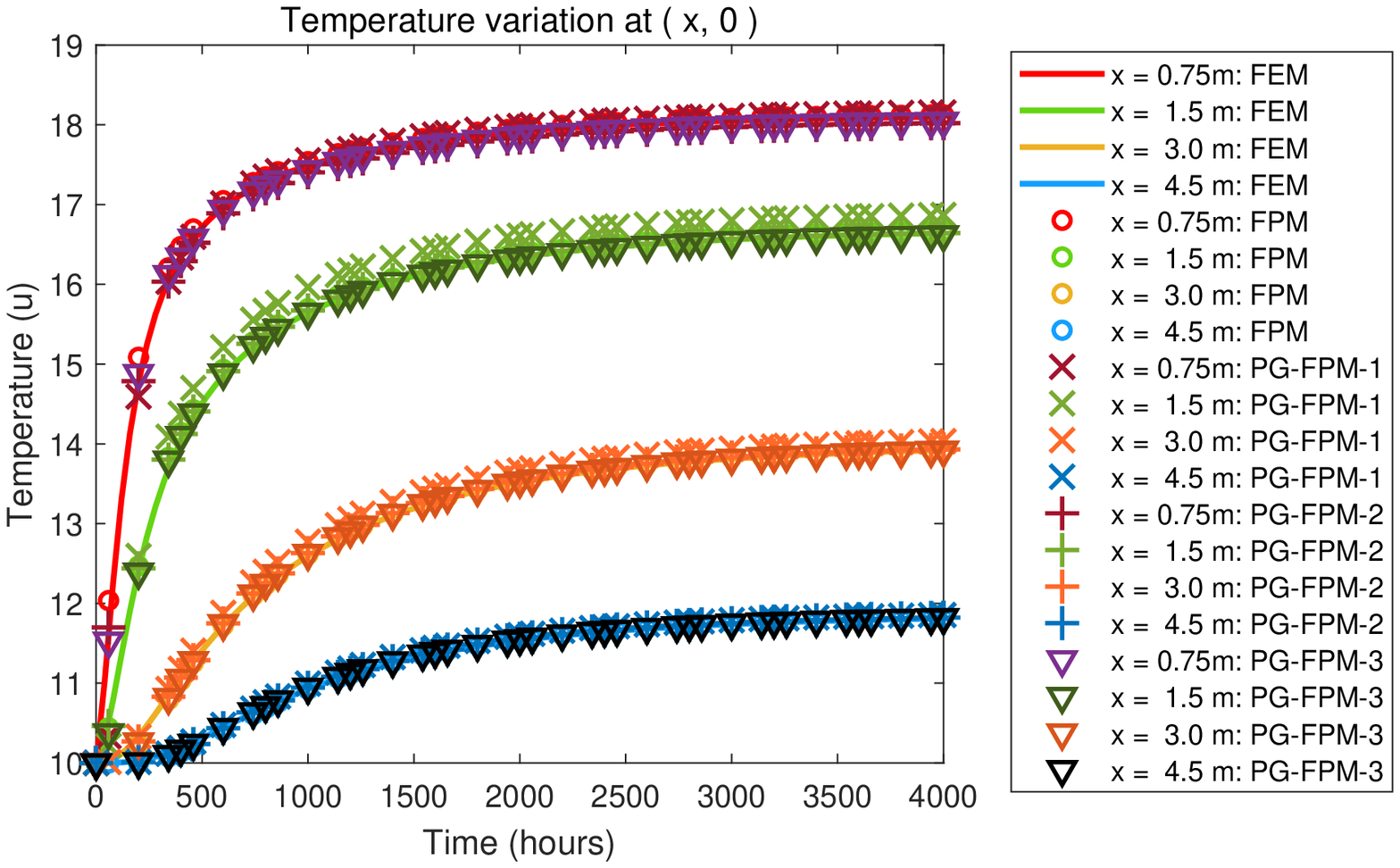}}  
    \subfigure{
    \label{fig:Ex19_01_400} 
    \includegraphics[width=0.3\textwidth]{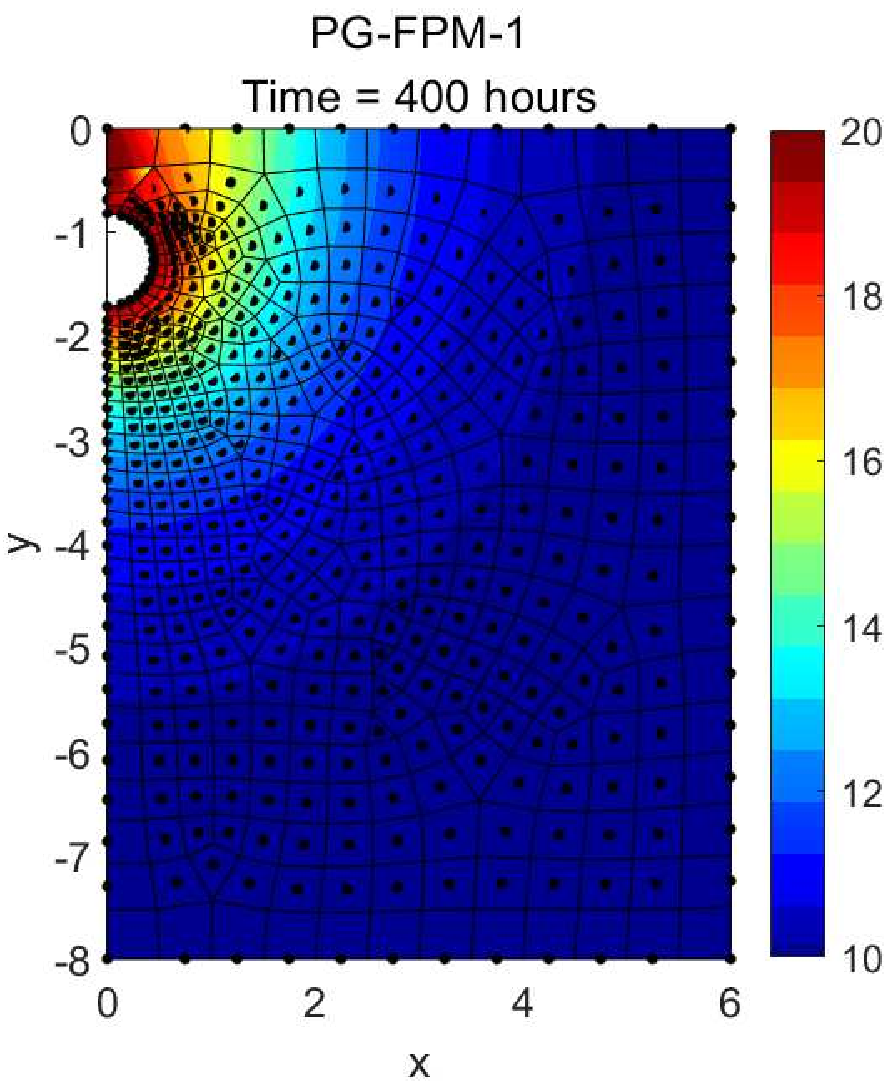}} 
    \addtocounter{subfigure}{-1} 
    \subfigure[]{ 
    \label{fig:Ex19_02_400} 
    \includegraphics[width=0.3\textwidth]{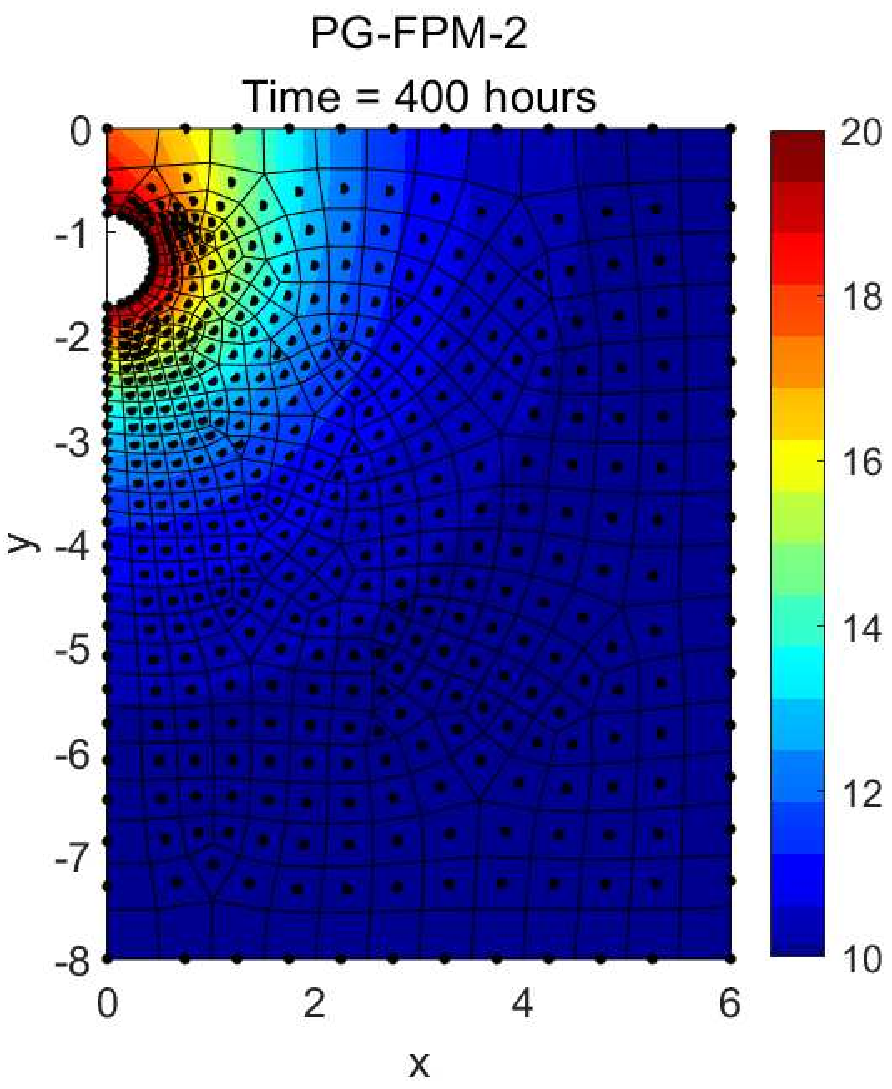}}  
    \subfigure{ 
    \label{fig:Ex19_03_400} 
    \includegraphics[width=0.3\textwidth]{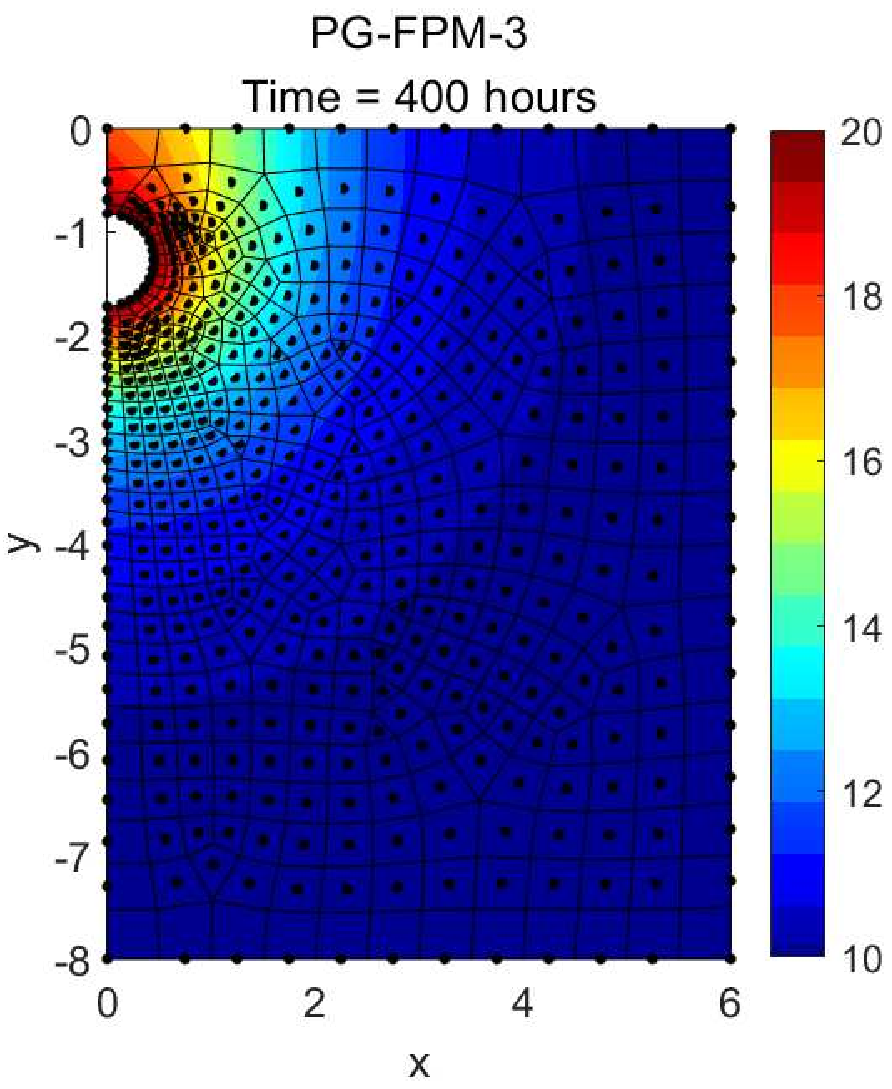}}  
\subfigure{
    \label{fig:Ex19_01_4000} 
    \includegraphics[width=0.3\textwidth]{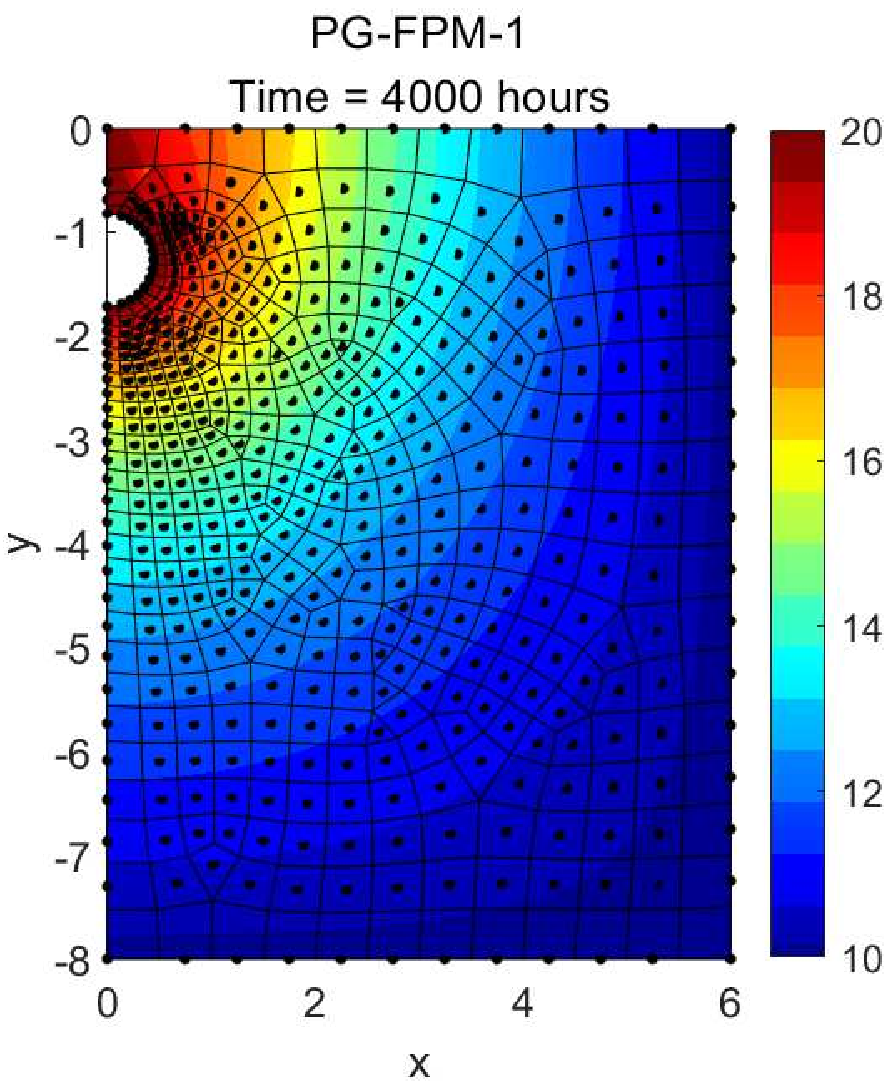}} 
    \addtocounter{subfigure}{-2} 
    \subfigure[]{ 
    \label{fig:Ex19_02_4000} 
    \includegraphics[width=0.3\textwidth]{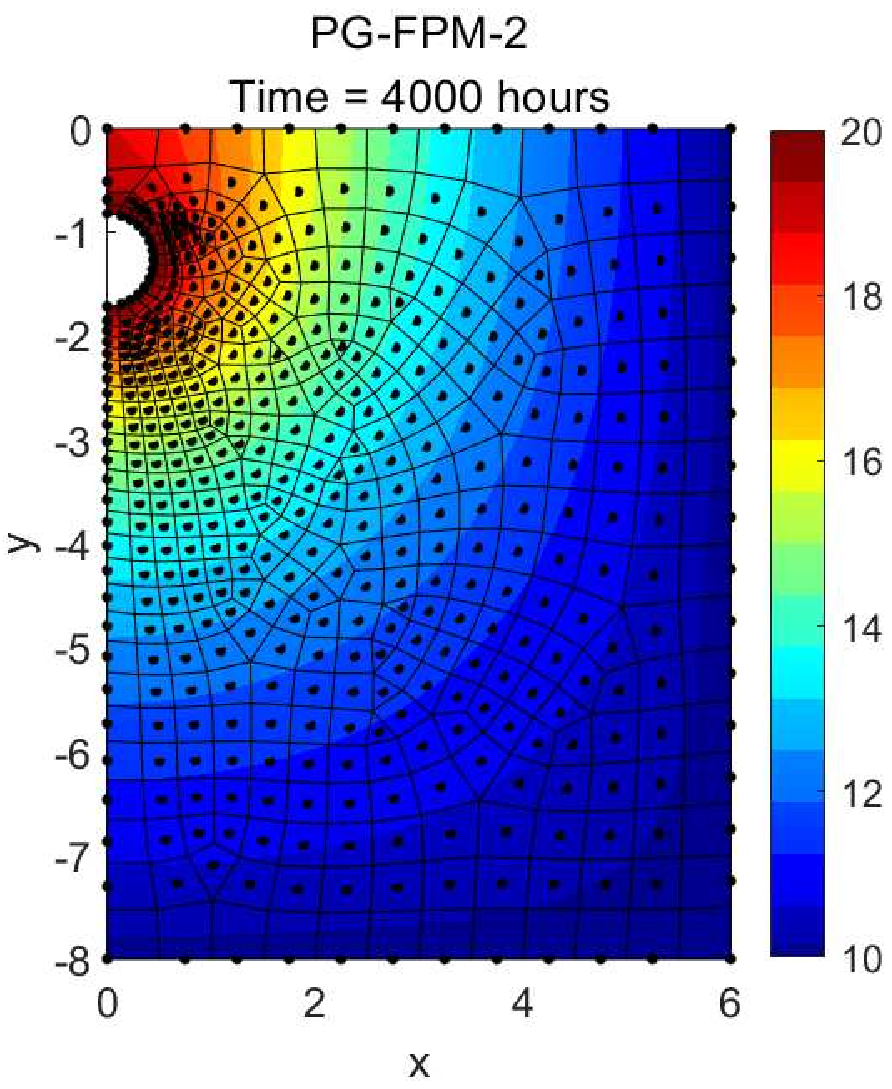}}  
    \subfigure{ 
    \label{fig:Ex19_03_4000} 
    \includegraphics[width=0.3\textwidth]{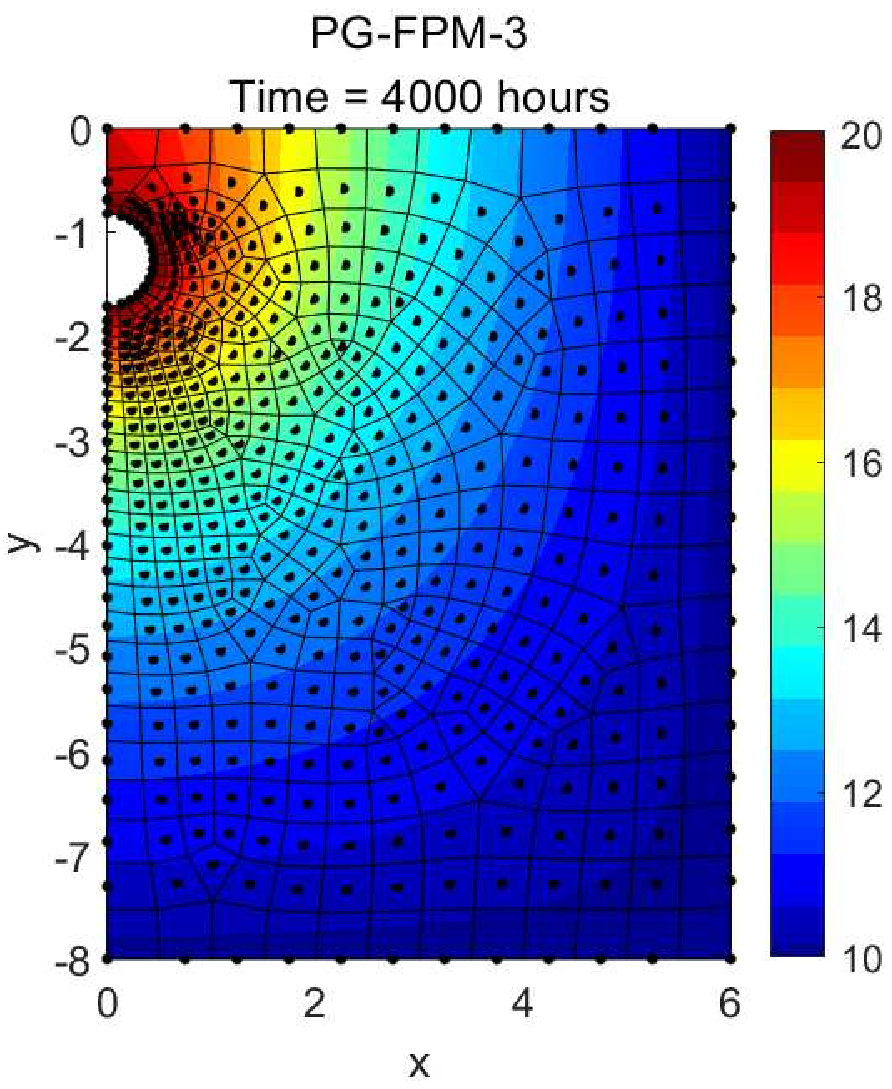}}  
  \caption{Ex. (1.9) – The boundary conditions and computed solutions.  (a) the problem domain and its partition converted from ABAQUS mesh. (b) transient temperature solution. (c) temperature distribution when $t=400$~hours achieved by the three PG-FPMs. (d) temperature distribution when $t=4000$~hours achieved by the three PG-FPMs.} 
  \label{fig:Ex19} 
\end{figure}

The last 2D example (Ex.~(1.9)) is a transient heat conduction problem in a semi-infinite isotropic soil medium caused by an oil pipe. As shown in Fig.~\ref{fig:Ex19_Partition}, according to the symmetry, only one half of the domain is under study. The pipe wall is modeled as a Dirichlet boundary with $\widetilde{u}_D = 20 \mathrm{^\circ C}$. The infinite boundaries are applied as $\widetilde{u}_D = 10 \mathrm{^\circ C}$. The left side is symmetric, and the top side is adiabatic. The two boundary conditions are equivalent in this example. The material properties: $\rho = 2620~\mathrm{kg / m^2}$, $c = 900~\mathrm{J / kg ^\circ C}$, $k = 2.92~\mathrm{W / m ^\circ C}$. The initial condition is $u (x, y, 0) = 10~\mathrm{^\circ C}$.

As can be seen in Fig.~\ref{fig:Ex19_Partition}, a quadrilateral partition is converted from ABAQUS mesh and exploited in the FPM and PG-FPMs. The partition includes 468 subdomains, and the corresponding Fragile Points are placed at the subdomain centroids. The point distribution is uneven. As a more violent variation of temperature is anticipated, more points are scattered in the vicinity of the pipe wall. Figure~\ref{fig:Ex19_Trans} shows the computed time variation of temperature at four representative points on the soil top side using the FPM and PG-FPMs. The results agree well with each other, as well as with a FEM solution achieved by ABAQUS with the same domain partition and DC2D4 elements. In the time domain, the LVIM approach is cooperated with the FPM and PG-FPMs with a time step $\Delta t = 500$~hours, while in ABAQUS, an explicit solver is exploited with $\Delta t = 40$~hours. The computed spatial temperature distribution at $t = 400$~hours and $t = 4000$~hours are displayed in Fig.~\ref{fig:Ex19_02_400} and \ref{fig:Ex19_02_4000} respectively. Here a comparison of the three proposed PG-FPM approaches is presented. A great consistency can be observed for the three PG-FPM approaches in solving complex 2D heat conductivity problems with unevenly distributed points. The accuracy of the results are also validated by ABAQUS solution and numerical solutions in previous studies \cite{Xu2010, Yu2010}.

According to the computational times of the different approaches listed in Table~\ref{table:Ex19}, the finite volume method (PG-FPM-2) still shows the highest efficiency. However, in this example, to assure a good accuracy and stability, more than one integration points are required in each subdomain in the singular solution method (PG-FPM-3). As a result, the PG-FPM-3 shows an unsatisfactory efficiency.

\begin{table}[htbp]
\caption{Computational times of the FPM and PG-FPMs in solving Ex.~(1.9).}
\centering
{
\begin{tabular*}{500pt}{@{\extracolsep\fill}lcccc@{\extracolsep\fill}}
\toprule
\textbf{Method} & \tabincell{c}{\textbf{Point distribution} \\ \textbf{(Number of points)}} & \tabincell{c}{\textbf{Computational} \\ \textbf{parameters}} & \tabincell{c}{$N_{band} (\mathbf{K})$ \& $N_{band} (\mathbf{C})$}  & \tabincell{c}{\textbf{Computational} \\ \textbf{time (s)}} \\
\midrule
\multirow{2}*{FPM} & ABAQUS (468) & $\eta_1 =2$, $\eta_2 = 1 \times 10^{5}$ & 23 \& 1 & 3.4 \\
~ & uniform (463) &$\eta_1 = 2$, $\eta_2 = 1 \times 10^{5}$ & 22 \& 1 & 2.8 \\
\midrule
\multirow{2}*{PG-FPM-1} & ABAQUS (468) & $\eta_1 =2$, $\eta_2 = 1 \times 10^{5}$ & 14 \& 1 & 2.8 \\
~ & uniform (463) &$\eta_1 = 2$, $\eta_2 = 1 \times 10^{5}$ & 14 \& 1 & 2.9 \\
\midrule
\multirow{2}*{PG-FPM-2} & ABAQUS (468) & $\eta_1 = 2$, $\eta_2 = 1 \times 10^{5}$ & 12 \& 1 & 2.1 \\
~ & uniform (463) &$\eta_1 = 2$, $\eta_2 = 1 \times 10^{5}$ & 12 \& 1  & 1.9 \\
\midrule
\multirow{2}*{PG-FPM-3} & ABAQUS (468) & $\eta_1 = 2$, $\eta_2 = 1 \times 10^{5}$ & 12 \& 5 & 9.3 \\
~ & uniform (463) &$\eta_1 = 2$, $\eta_2 = 1 \times 10^{5}$ & 12 \& 5 & 8.0 \\
\bottomrule
\end{tabular*}}
\label{table:Ex19}
\end{table}

\section{3D examples} \label{sec:3D}

\subsection{Anisotropic nonhomogeneous examples in a cubic domain}

\begin{figure}[htbp]
    \centering
    \begin{minipage}[t]{0.48\textwidth}
        \centering
        \includegraphics[width=1\textwidth]{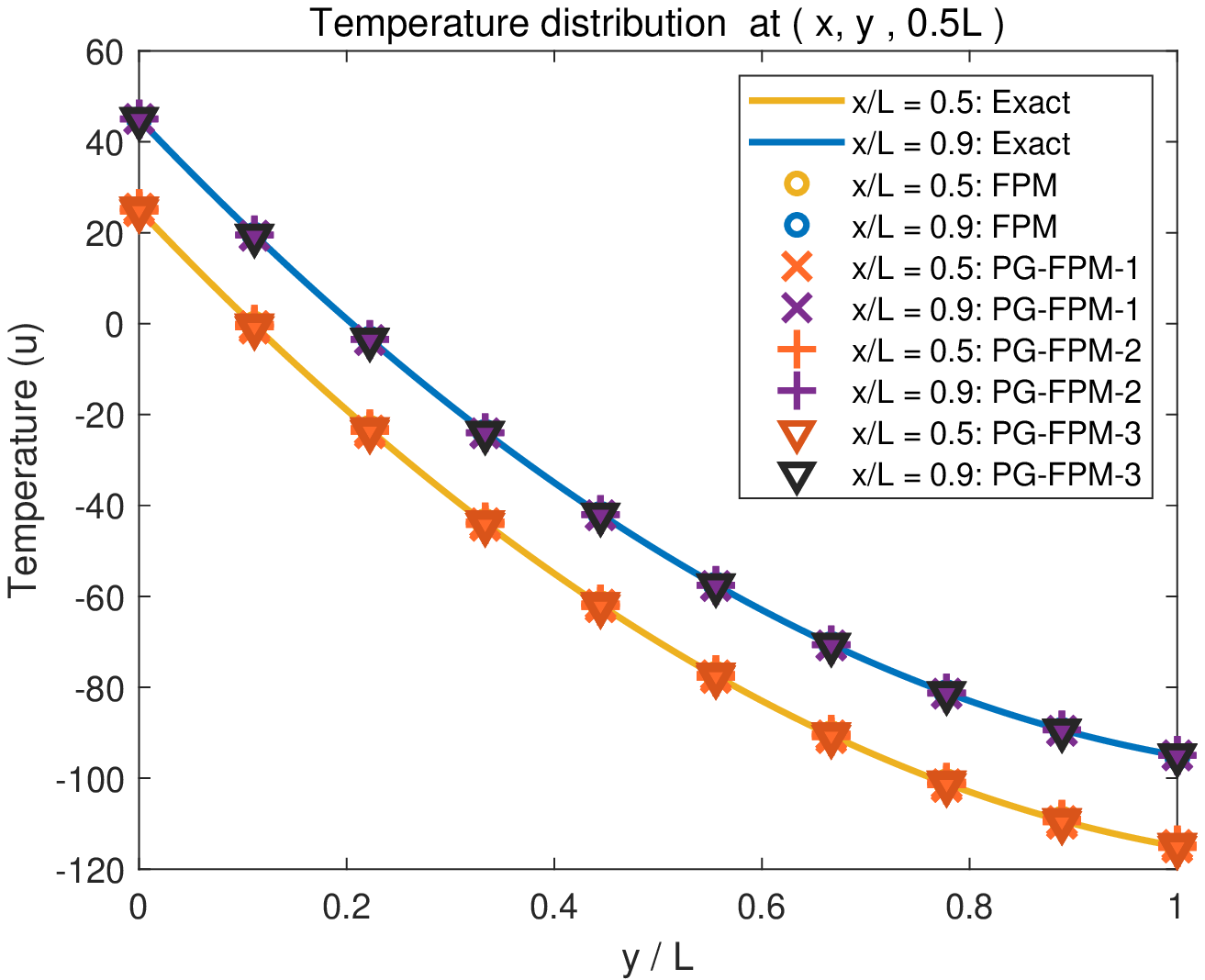}
        \caption{Ex. (2.1) - The computed solutions at $z = 0.5L$.}
        \label{fig:Ex21}
    \end{minipage}
    \begin{minipage}[t]{0.48\textwidth}
        \centering
        \includegraphics[width=1\textwidth]{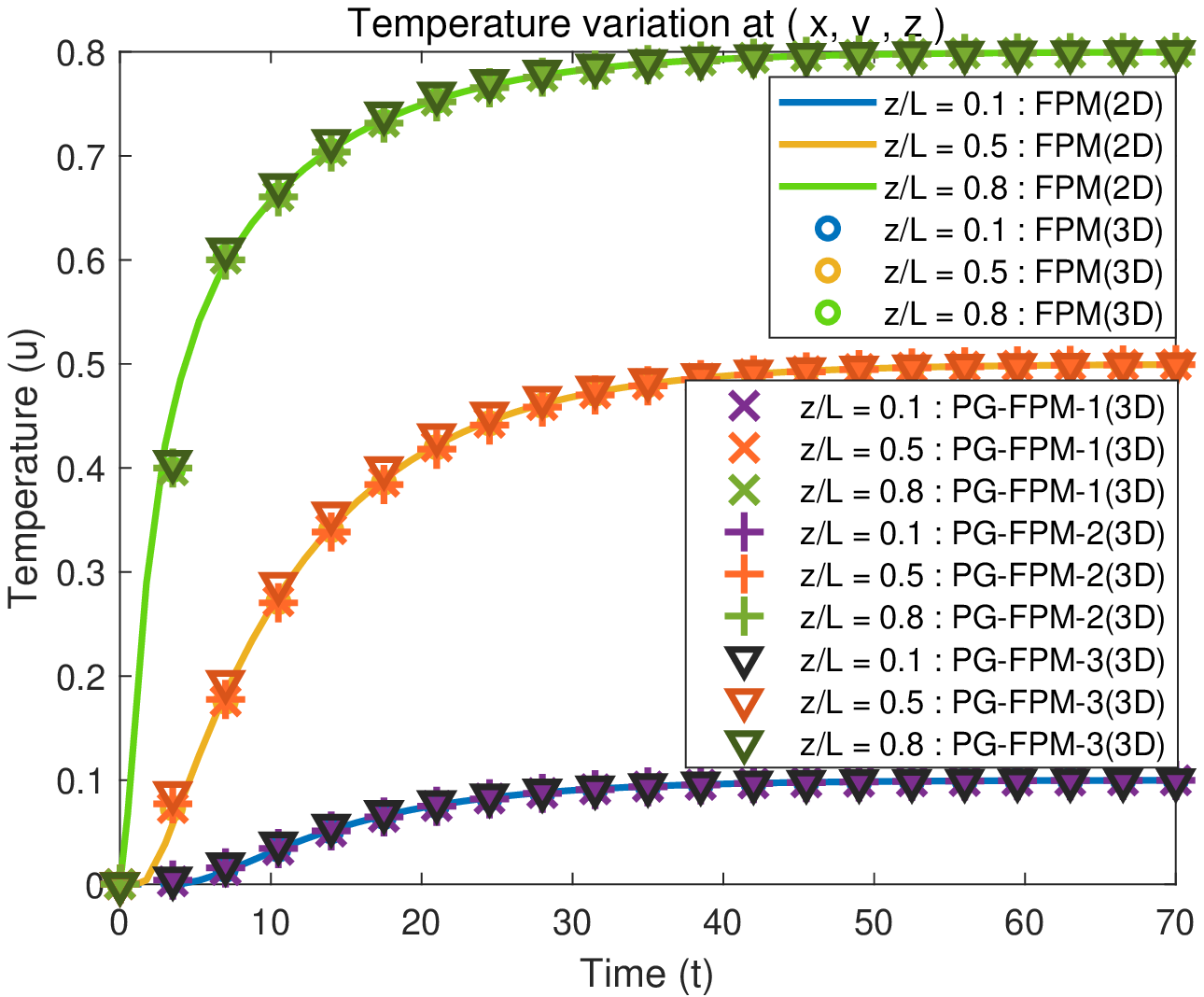}
        \caption{Ex. (2.2) - The computed transient temperature solution.}
        \label{fig:Ex22}
    \end{minipage}
\end{figure}

In the 3D examples, first, we consider a number of heat conduction problems in a $L \times L \times L$ cubic domain. The heat source density $Q$ vanishes in all the following examples. Ex.~(2.1) is a steady-state problem in homogenous anisotropic medium. The thermal conductivity tensor components $k_{11} = k_{22} = k_{33} = 1 \times 10^{-4}$, $k_{23} = 0.2 \times 10^{-4}$, $k_{12} = k_{13} = 0$. A postulated analytical solution is given \cite{Sladek2008, Sladek2012}:
\begin{align}
\begin{split}
u (x, y, z) = y^2+ y -5yz +xz. 
\end{split}
\end{align}
The side length $L = 10$. Dirichlet boundary conditions are applied on all the external faces. The problem is analyzed by the FPM and PG-FPMs with 1000 uniformly distributed Fragile Points. In this example, some Fragile Points are scattered on the external boundaries, hence, the Dirichlet boundary conditions can be enforced either strongly or weakly by the IP Numerical Flux Corrections. In the collocation method (PG-FPM-1), the constant parameter $c = 10$ and remains the same in the following 3D examples. The computed temperature distribution on $z = 0.5 L$ is shown in Fig.~\ref{fig:Ex21}. The corresponding relative errors and computational times are listed in Table~\ref{table:Ex21}. As can be seen, while all the approaches achieve excellent accuracy in this anisotropic example with a relative error less than 1\%, the finite volume method (PG-FPM-2) takes only 60\% of the computational time of the conventional Galerkin FPM and shows the best accuracy among all the approaches. The efficiency of the collocation method (PG-FPM-1) is unsatisfactory, since the superiority of this approach mainly lies in transient analysis.

Next, we consider a simple transient heat conduction problem. The material is homogenous and isotropic, with material properties: $\rho = 1$, $c = 1$, and $k = 1$. The boundary condition on the top surface ($z = L$) is prescribed as a thermal shock $\widetilde{u}_D = H (t-0)$, where $H$ is the Heaviside time step function. The bottom boundary condition on $z = 0$ is $\widetilde{u}_D = 0$. Heat fluxes vanish on all the lateral surfaces. The initial condition $u(x, y, z, 0) = 0$. Clearly, the temperature field is independent of $x$ and $y$ coordinates. Hence the example can also be analyzed in 2D. In Fig.~\ref{fig:Ex22}, the computed transient temperature on $z = 0.1 L$, $z = 0.5 L$ and $z = 0.8 L$ solved by the 3D FPM and PG-FPMs are presented, in comparison with a 2D FPM solution. With 1000 uniformly distributed points and $\Delta t = 7$, $M = 3$ in the time domain, the time costs of the 3D approaches are shown in Table~\ref{table:Ex22}. Achieving approximately identical results, the finite volume method (PG-FPM-2) cuts the computational time by a half as compared to the original Galerkin FPM. The other PG-FPM approaches also improve the performance of the method more or less.

\begin{table}[htbp]
\caption{Relative errors and computational times of the FPM and PG-FPMs in solving Ex.~(2.1).}
\centering
{
\begin{tabular*}{500pt}{@{\extracolsep\fill}lcccc@{\extracolsep\fill}}
\toprule
\textbf{Method} & \tabincell{c}{\textbf{Computational} \\ \textbf{parameters}} &  \textbf{Relative errors} & $N_{band} (\mathbf{K})$ & \tabincell{c}{\textbf{Computational} \\ \textbf{time (s)}} \\
\midrule
FPM & $\eta_1 = 1$, $\eta_2 = 1 \times 10^{5}$ & \tabincell{c}{$e_0 = 5.1 \times 10^{-3}$\\$e_1 = 1.1 \times 10^{-1}$} & 48 & 2.6 \\
PG-FPM-1 & $\eta_1 = 1$, $\eta_2 = 1 \times 10^{5}$ & \tabincell{c}{$e_0 = 2.8 \times 10^{-3}$\\$e_1 = 1.3 \times 10^{-1}$} & 37 & 3.5 \\
PG-FPM-2 & $\eta_1 = 1$, $\eta_2 = 1 \times 10^{5}$ & \tabincell{c}{$e_0 = 5.8 \times 10^{-4}$\\$e_1 = 7.2 \times 10^{-2}$} & 20 & 1.6 \\
PG-FPM-3 & $\eta_1 = 1$, $\eta_2 = 1 \times 10^{5}$ & \tabincell{c}{$e_0 = 9.9 \times 10^{-4}$\\$e_1 = 7.2 \times 10^{-2}$} & 21 & 1.9 \\
\bottomrule
\end{tabular*}}
\label{table:Ex21}
\end{table}

\begin{table}[htbp]
\caption{Computational times of the FPM and PG-FPMs in solving Ex.~(2.2).}
\centering
{
\begin{tabular*}{500pt}{@{\extracolsep\fill}lccc@{\extracolsep\fill}}
\toprule
\textbf{Method} & \tabincell{c}{\textbf{Computational} \textbf{parameters}} & $N_{band} (\mathbf{K})$ \& $N_{band} (\mathbf{C})$ & \tabincell{c}{\textbf{Computational} \textbf{time (s)}} \\
\midrule
FPM & $\eta_1 = 1$, $\eta_2 = 20$ & 49 \& 1 & 6.6 \\
PG-FPM-1 & $\eta_1 = 0$, $\eta_2 = 20$ & 14 \& 1 & 5.2 \\
PG-FPM-2 & $\eta_1 = 1$, $\eta_2 = 20$ & 11 \& 1 & 3.3 \\
PG-FPM-3 & $\eta_1 = 1$, $\eta_2 = 20$ & 21 \& 1 & 4.1 \\
\bottomrule
\end{tabular*}}
\label{table:Ex22}
\end{table}

In Ex.~(2.3) – (2.6), we consider the same initial boundary condition problem as Ex.~(2.2) with different material properties. First, a homogenous anisotropic material is assumed. The thermal conductivity tensor components are given as: $k_{11} = k_{33} = 1$, $k_{22} = 1.5$, $k_{23} = 0.5$, $k_{12} = k_{13} = 0$. Symmetric boundary conditions are given on the left and right surfaces ($x = 0, L$) instead of the free boundary conditions. The temperature field is still independent of $x$ and can be equivalent to a 2D problem. Figure~\ref{table:Ex22} shows the representative transient and spatial temperature distribution solutions acquired by the FPM and PG-FPMs with 1331 Fragile Points, in comparison with a 2D FPM solution. While most of the computed results agree well with each other, the transient temperature on $z = 0.2 L$ achieved by the collocation method (PG-FPM-1) is slightly off the other solutions. This implies that in 3D collocation method, the accuracy in the vicinity of external boundaries may be less than the other Points, since the approximation of high order derivatives using the local RBF-DQ method may have larger errors when the target point is on marginal of the supporting points. This problem can be remedied by distributing more Points close to the external boundaries.

\begin{figure}[htbp] 
  \centering 
    \subfigure[]{ 
    \label{fig:Ex23_Trans} 
    \includegraphics[width=0.48\textwidth]{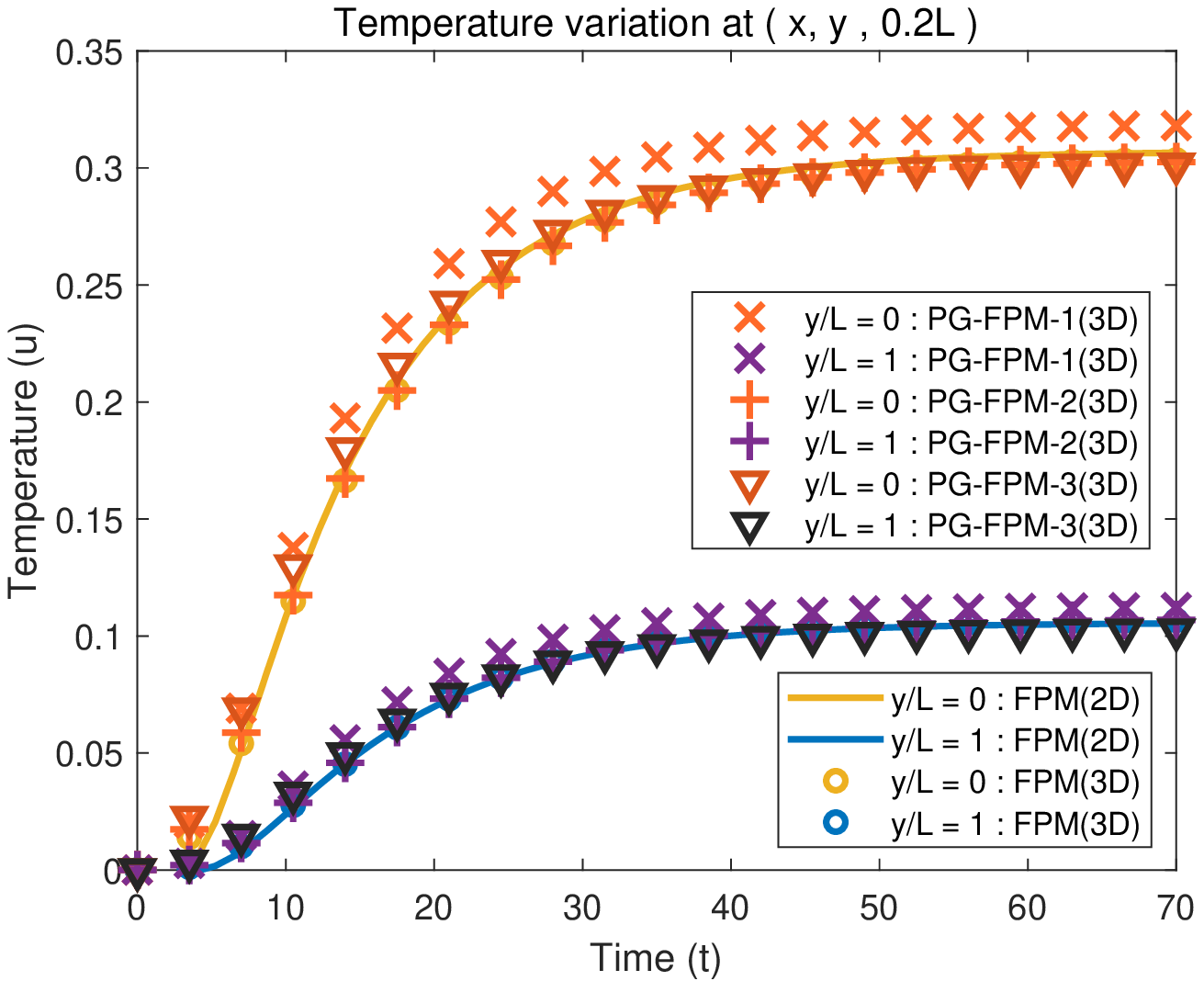}}  
    \subfigure[]{ 
    \label{fig:Ex23_SS} 
    \includegraphics[width=0.48\textwidth]{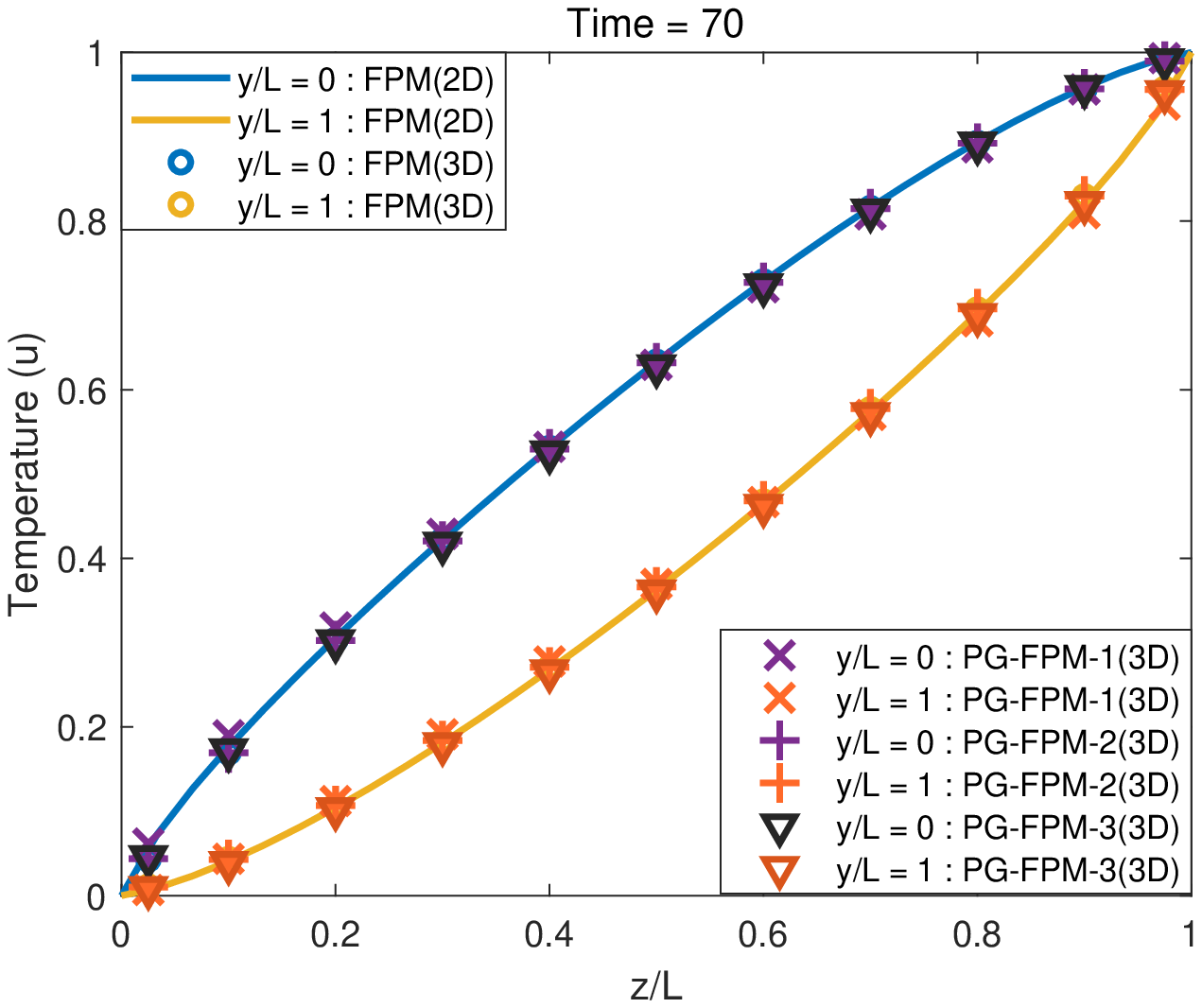}}  
  \caption{Ex. (2.3) - The computed solutions. (a) transient temperature solution at $z = 0.2L$. (b) temperature distribution when $t  = 70$.} 
  \label{fig:Ex23} 
\end{figure}

\begin{figure}[htbp] 
  \centering 
    \subfigure[]{ 
    \label{fig:Ex24_Trans} 
    \includegraphics[width=0.48\textwidth]{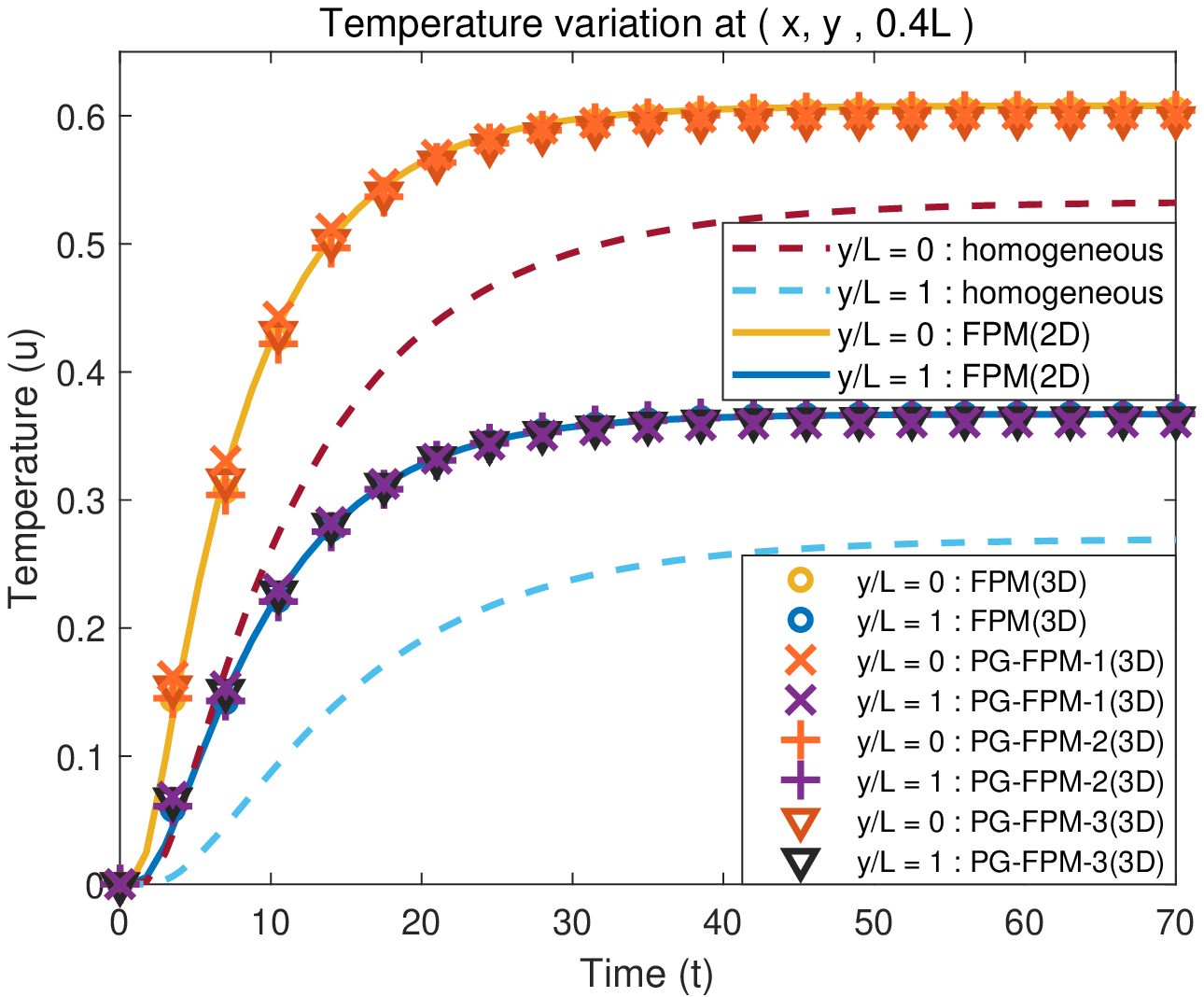}}  
    \subfigure[]{ 
    \label{fig:Ex24_SS} 
    \includegraphics[width=0.48\textwidth]{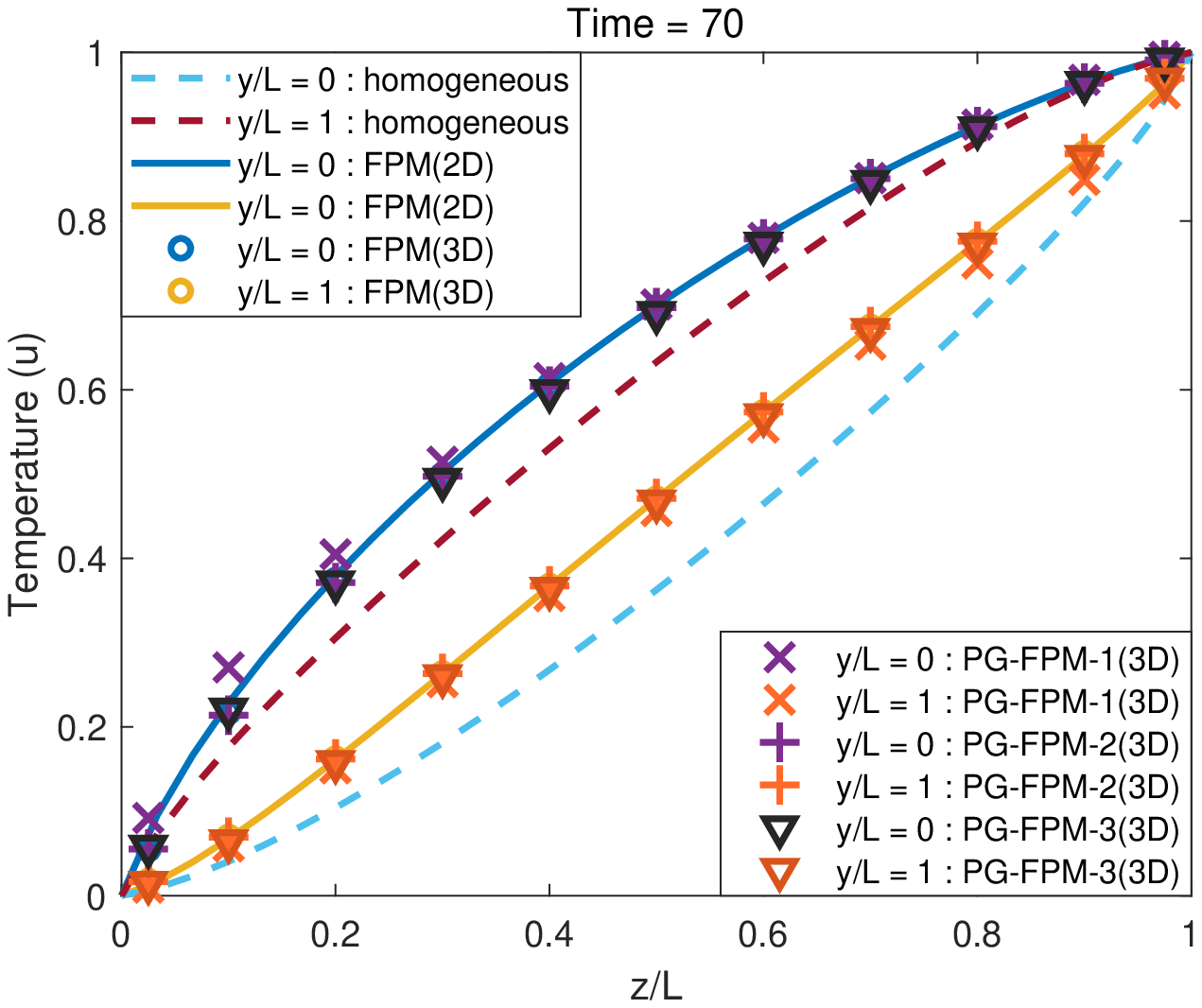}}  
  \caption{Ex. (2.4) - The computed solutions. (a) transient temperature solution at $z = 0.4L$. (b) temperature distribution when $t  = 70$.} 
  \label{fig:Ex24} 
\end{figure}

Next, in Ex.~(2.4), a nonhomogeneous material is considered. The thermal conductivity component $k_{33} (z) = 1 + z/L$, while all the other material properties remain the same as Ex.~(2.3). This example is still equivalent to 2D. The computed solutions obtained by 2D / 3D FPM and 3D PG-FPMs are shown in Fig.~\ref{fig:Ex24}. The homogenous result is also given for comparison. As can be seen, the nonhomogeneity has a considerable influence on the temperature distribution. Yet all the 3D PG-FPM approaches exhibit good performance in analyzing the nonhomogeneous problem. 

\begin{figure}[htbp] 
  \centering 
    \subfigure[]{ 
    \label{fig:Ex25_Trans} 
    \includegraphics[width=0.48\textwidth]{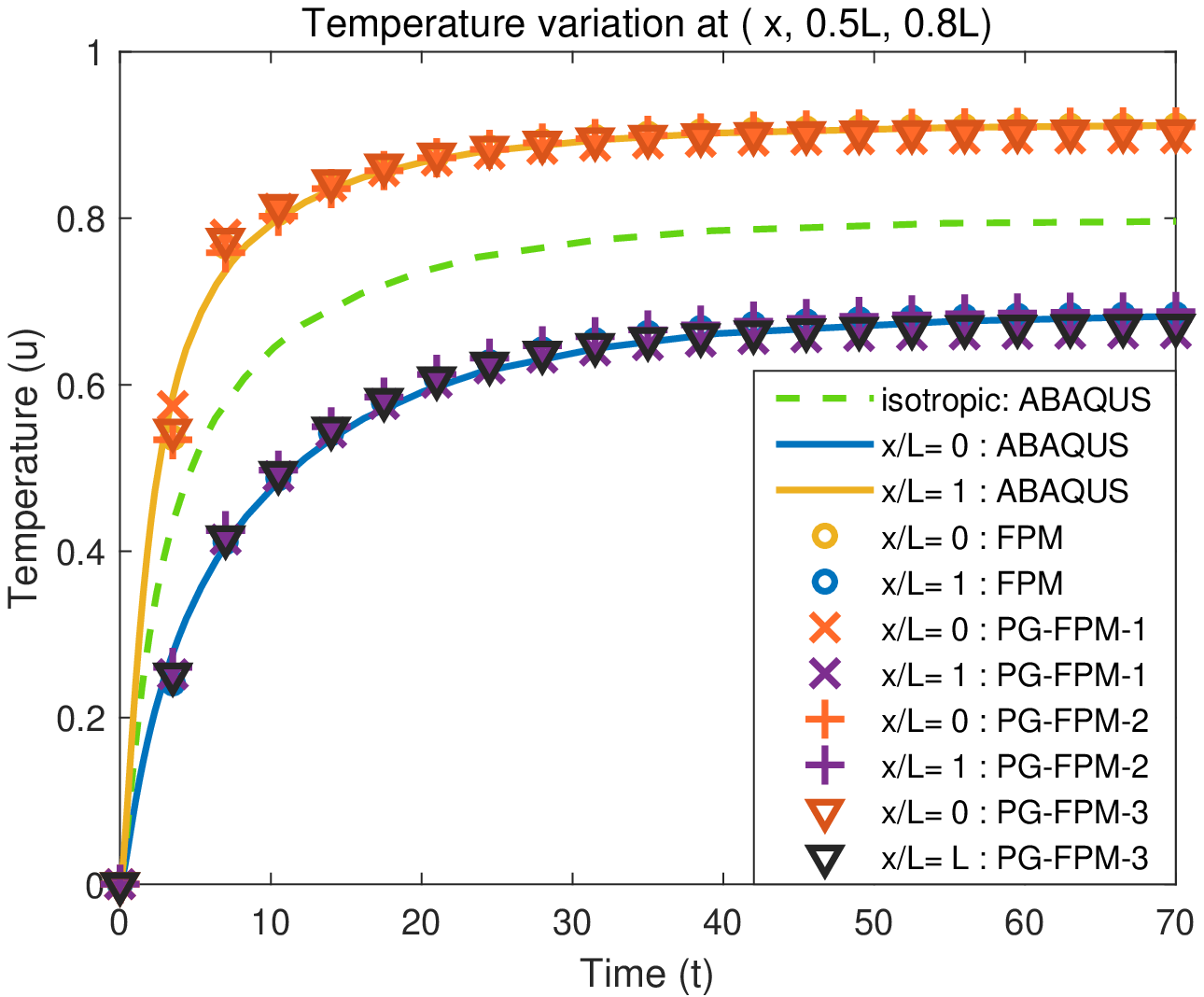}}  
    \subfigure[]{ 
    \label{fig:Ex25_SS} 
    \includegraphics[width=0.48\textwidth]{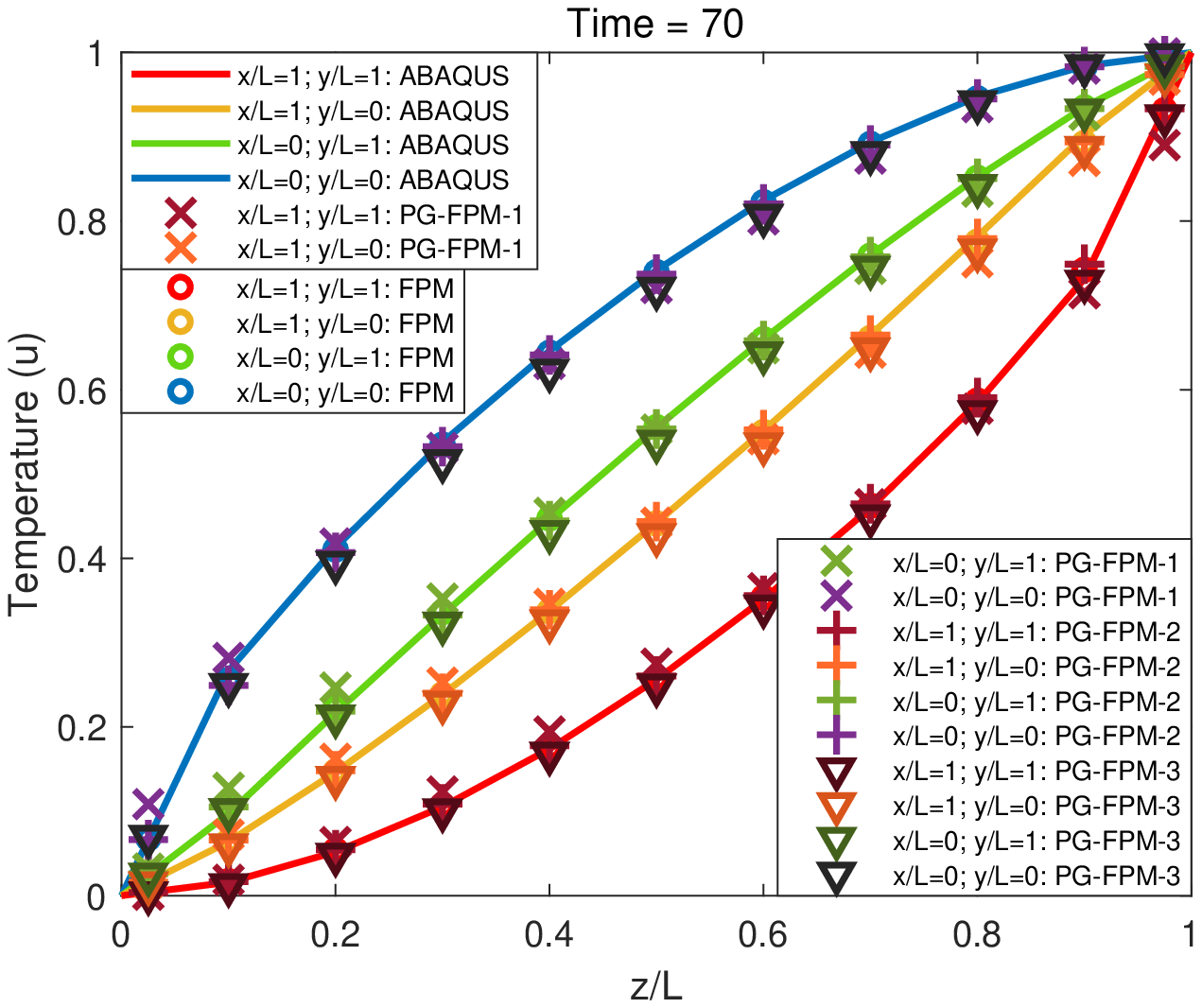}}  
  \caption{Ex. (2.5) - The computed solutions. (a) transient temperature solution at $y = 0.5L$, $z = 0.8L$. (b) temperature distribution when $t  = 70$.} 
  \label{fig:Ex25} 
\end{figure}

In Ex.~(2.5), we consider an example that can no longer be analyzed in 2D. The homogenous anisotropic thermal conductivity coefficients are: $k_{11} = k_{33} = 1$, $k_{22} = 1.5$, $k_{12}=k_{13} = k_{23} = 0.5$. All the lateral surfaces have vanishing heat fluxes. The computed solutions are presented in Fig.~\ref{fig:Ex25}, comparing with the ABAQUS solution using 1000 DC3D8 elements. A good consistency is observed between the PG-FPMs and ABAQUS results. The anisotropy of the material has no influence on the accuracy or efficiency of all the proposed PG-FPM approaches. 

\begin{figure}[htbp]
    \centering
    \begin{minipage}[t]{0.48\textwidth}
        \centering
        \includegraphics[width=1\textwidth]{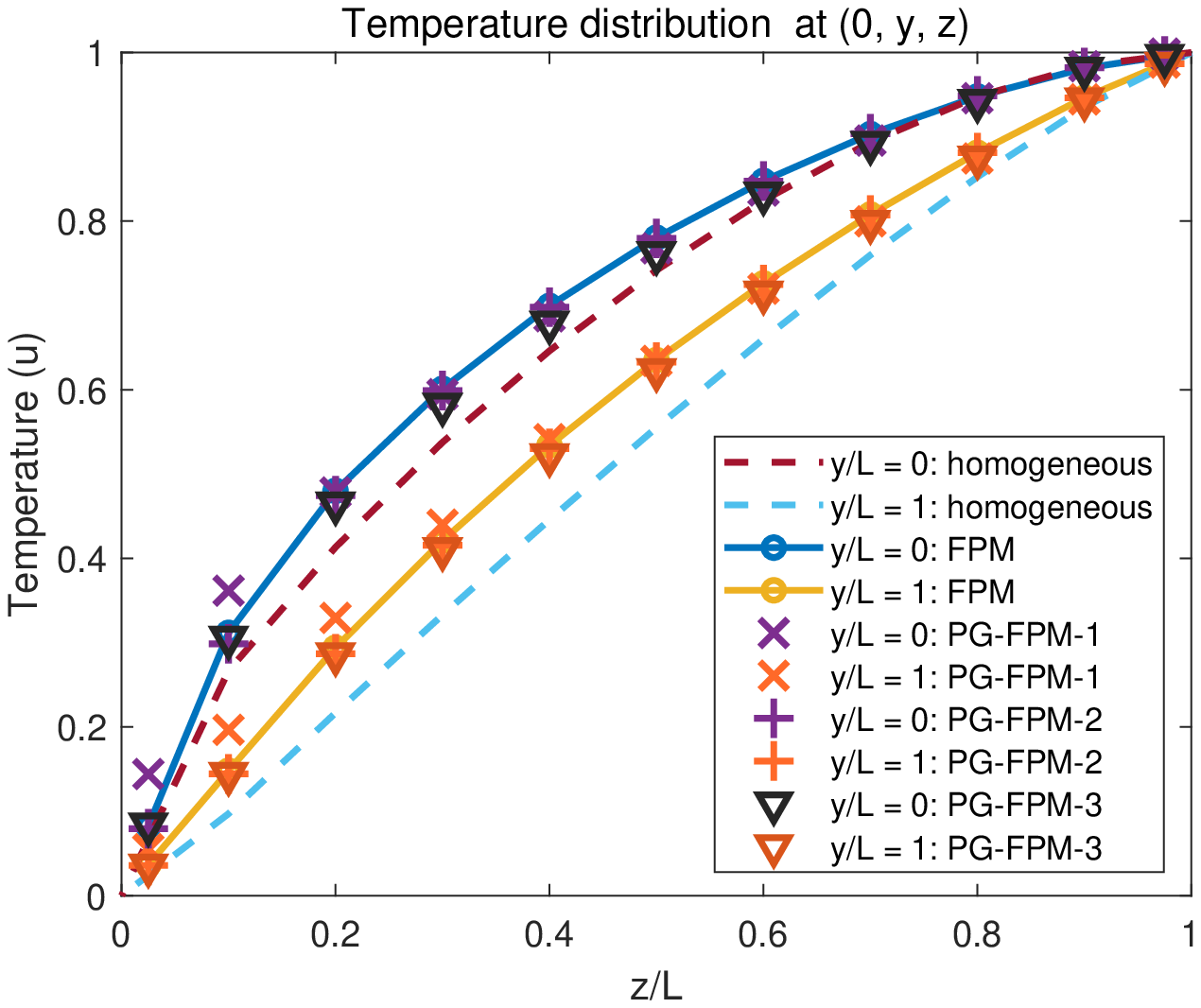}
        \caption{Ex. (2.6) - The computed steady-state result.}
        \label{fig:Ex26}
    \end{minipage}
    \begin{minipage}[t]{0.48\textwidth}
        \centering
        \includegraphics[width=1\textwidth]{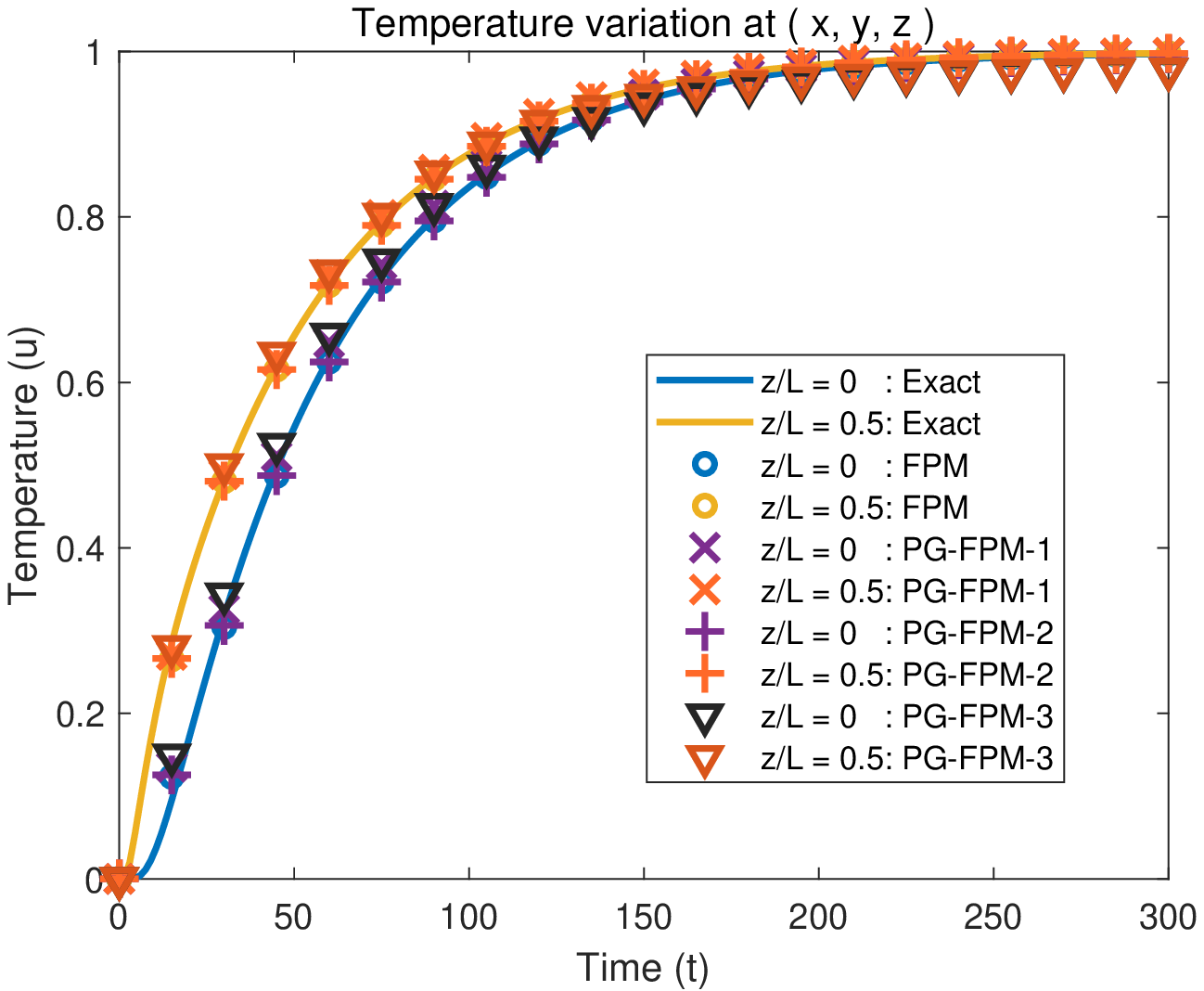}
        \caption{Ex. (2.7) - The computed transient temperature solution.}
        \label{fig:Ex27}
    \end{minipage}
\end{figure}

At last, a nonhomogeneous anisotropic material is considered. All the material properties keep the same as Ex.~(2.5) except $k_{33} (z) = 1 + z/L$. The computed steady-state temperature distribution is shown in Fig.~\ref{fig:Ex26}. All the PG-FPM solutions agree well with the FPM result. These solutions, as well as all the previous solutions in Ex.~(2.2) – Ex.~(2.5), are also consistent with numerical results achieved using the MLPG method \cite{Sladek2008}.

The computational times for Ex.~(2.3) – (2.6) in different materials with nonhomogeneity and / or anisotropy are roughly the same. As listed in Table~\ref{table:Ex23}, all the three PG-FPMs have higher efficiency than the original Galerkin FPM. The finite volume method (PG-FPM-2) shows the best performance, followed by the singular solution method (PG-FPM-3) and collocation method (PG-FPM-1).

In Ex.~(2.7), we come to a transient heat conduction problem with Robin boundary condition. Homogenous and isotropic material is assumed, with $\rho = 1$, $c = 1$, $k = 1$. The top surface is under Robin boundary condition. The heat transfer coefficient $h = 1.0$, and the temperature outside the top surface is $\widetilde{u}_R = H(t-0)$, where $H$ is the Heaviside time step function. All the lateral and bottom surfaces are free ($\widetilde{q}_N = 0$). The initial condition is constant: $u (x,y,z,0)=0$. Clearly, the temperature solution is not dependent on $x$ and $y$ and can be equivalent to a 1D problem. The side length $L = 10$. The computed time variation of temperature on $z = 0$ and $z = 0.5 L$ are presented in Fig.~\ref{fig:Ex27}, comparing with the exact analytical solution \cite{Jaeger1959}. As can be seen, all the FPM and PG-FPMs achieve excellent accuracy in this example. The solutions also agree well with numerical results studied by Sladek et al. \cite{Sladek2008, Sladek2012}. With the same point distribution and time step, the computational time of Ex.~(2.7) is approximately the same as the previous examples and hence is omitted here.

\begin{table}[htbp]
\caption{Computational times of the FPM and PG-FPMs in solving Ex.~(2.3) – (2.6).}
\centering
{
\begin{tabular*}{500pt}{@{\extracolsep\fill}lccc@{\extracolsep\fill}}
\toprule
\textbf{Method} & \tabincell{c}{\textbf{Computational} \textbf{parameters}} & $N_{band} (\mathbf{K})$ \& $N_{band} (\mathbf{C})$ & \tabincell{c}{\textbf{Computational} \textbf{time (s)}} \\
\midrule
FPM & $\eta_1 = 1$, $\eta_2 = 20$ & 49 \& 1 & 7.8 \\
PG-FPM-1 & $\eta_1 = 0$, $\eta_2 = 20$ & 14 \& 1 & 5.4 \\
PG-FPM-2 & $\eta_1 = 1$, $\eta_2 = 20$ & 21 \& 1 & 4.2 \\
PG-FPM-3 & $\eta_1 = 1$, $\eta_2 = 20$ & 21 \& 1 & 4.6 \\
\bottomrule
\end{tabular*}}
\label{table:Ex23}
\end{table}

\subsection{Some practical examples}

\begin{figure}[htbp] 
  \centering 
   \subfigure[]{ 
    \label{fig:Ex28_BC} 
    \includegraphics[width=0.48\textwidth]{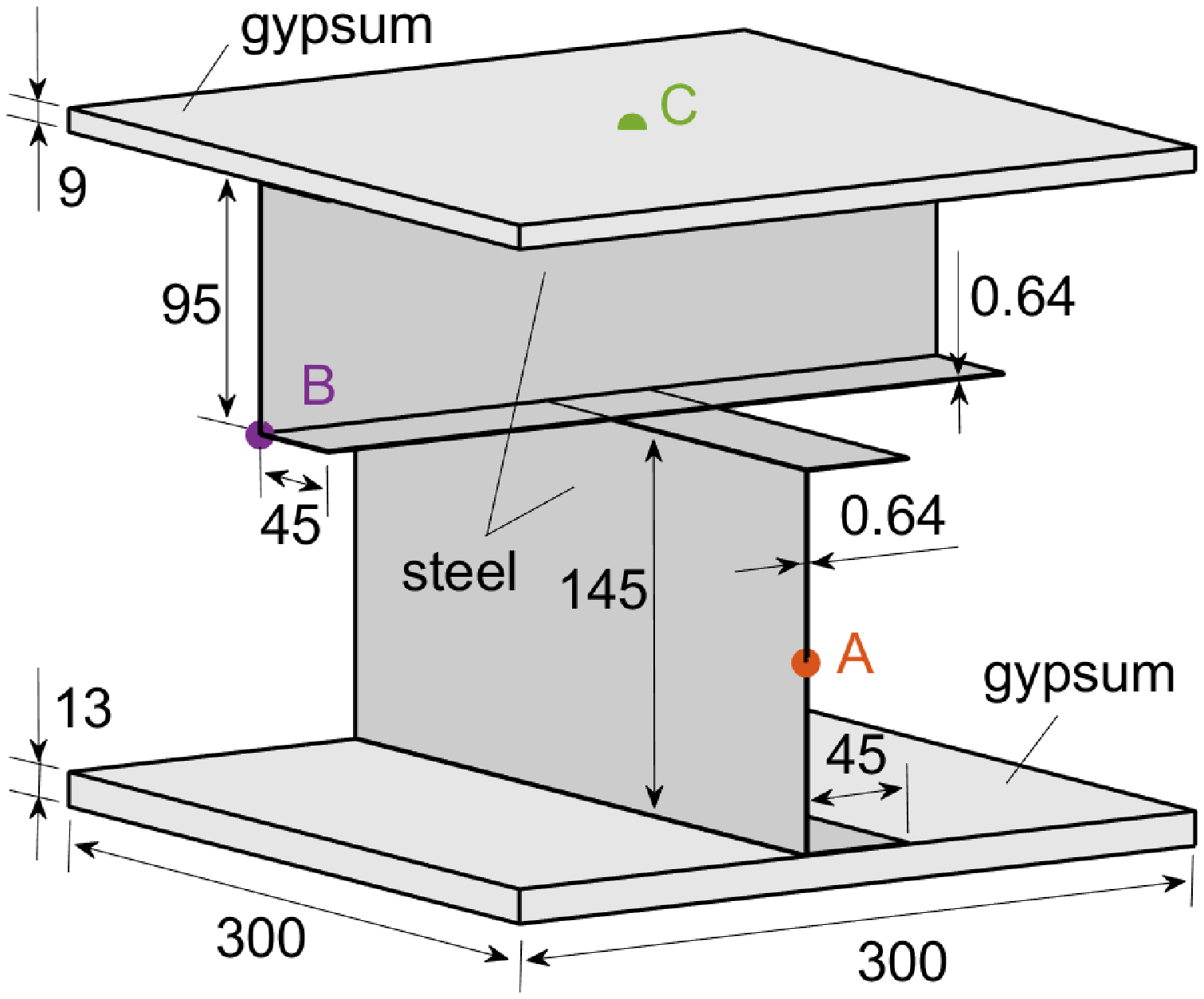}}  
   \subfigure[]{ 
    \label{fig:Ex28_Trans} 
    \includegraphics[width=0.48\textwidth]{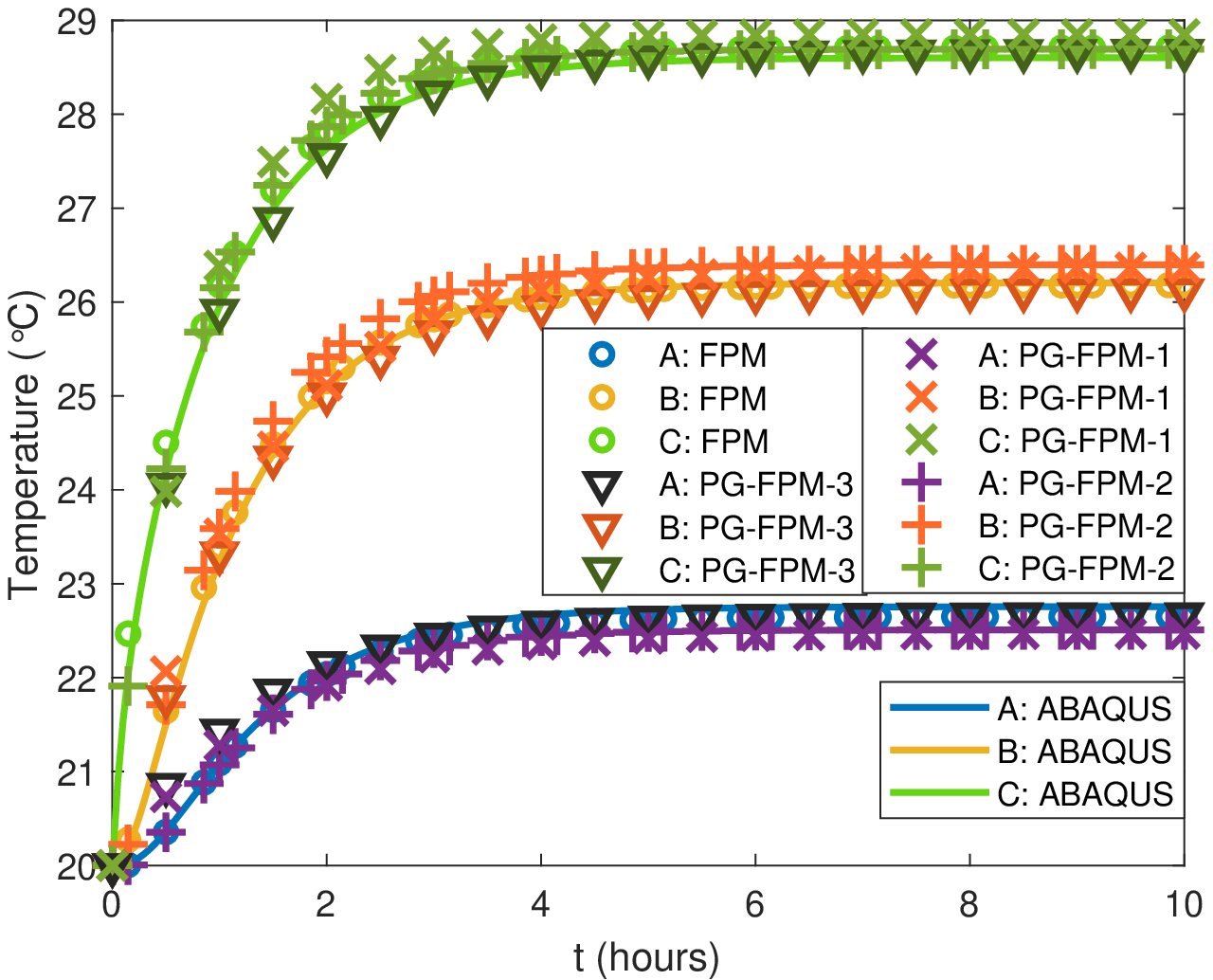}}  
    \subfigure[]{ 
    \label{fig:Ex28_S005_01} 
    \includegraphics[width=0.48\textwidth]{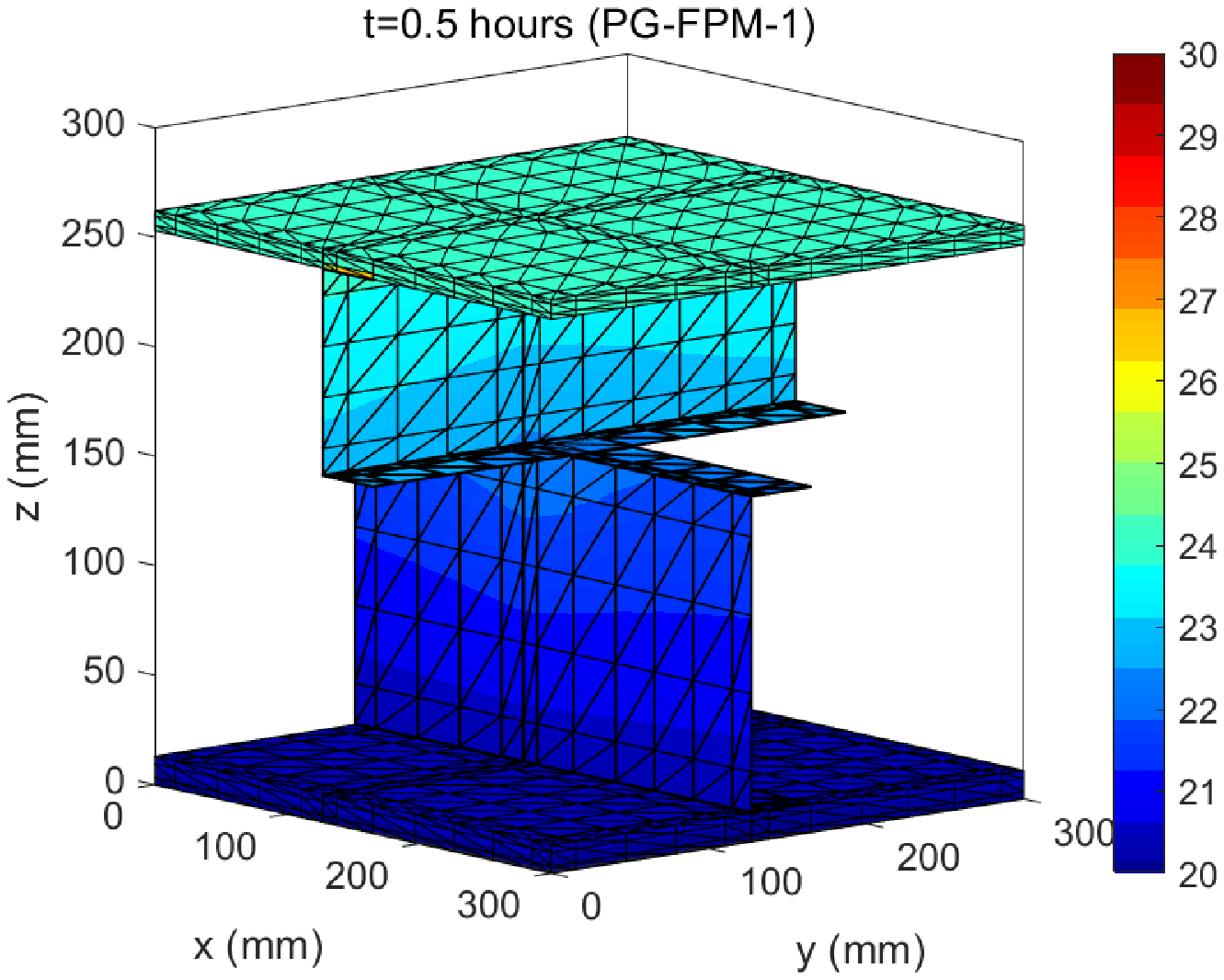}}  
    \subfigure[]{ 
    \label{fig:Ex28_S01_02} 
    \includegraphics[width=0.48\textwidth]{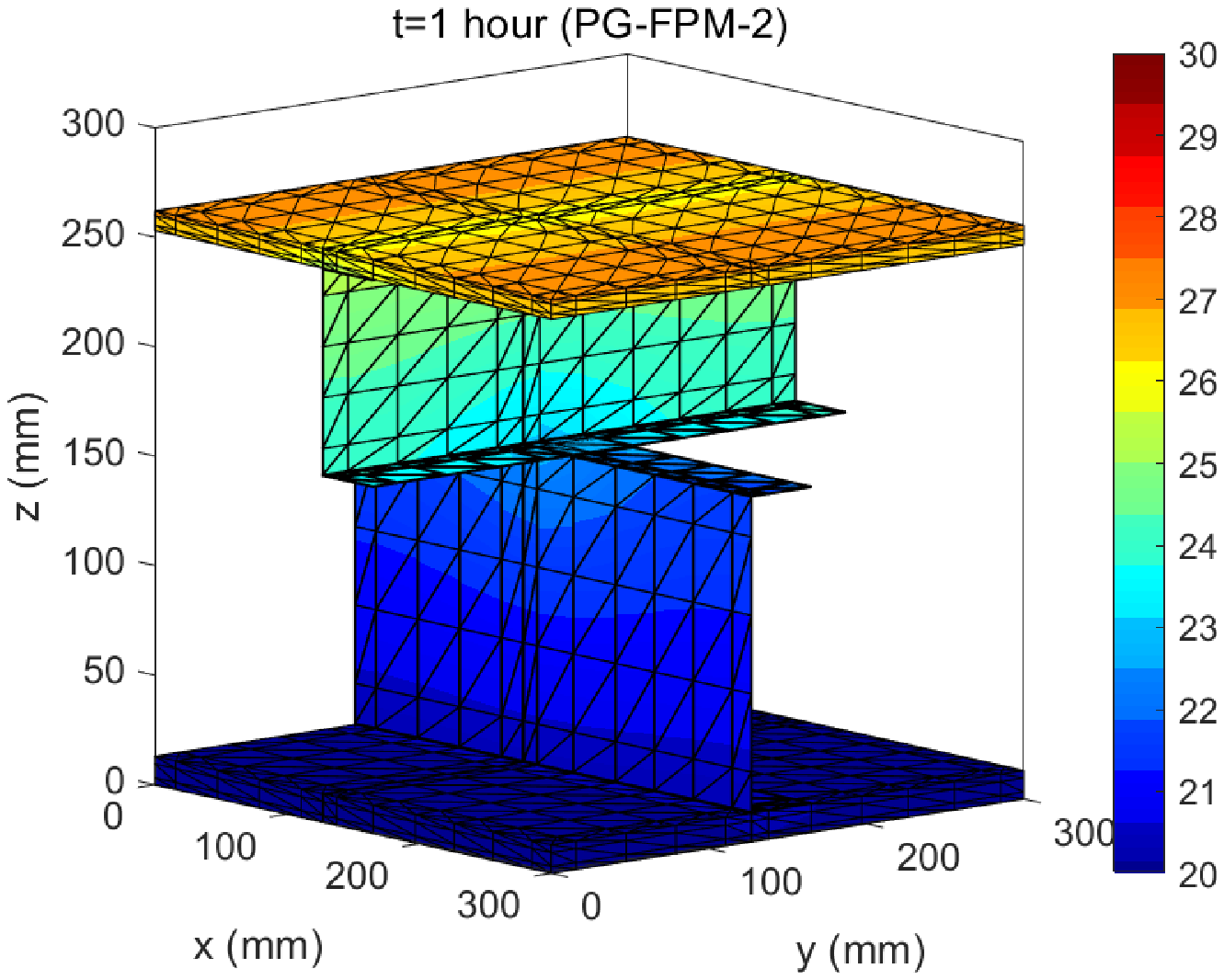}} 
    \subfigure[]{ 
    \label{fig:Ex28_S02_03} 
    \includegraphics[width=0.48\textwidth]{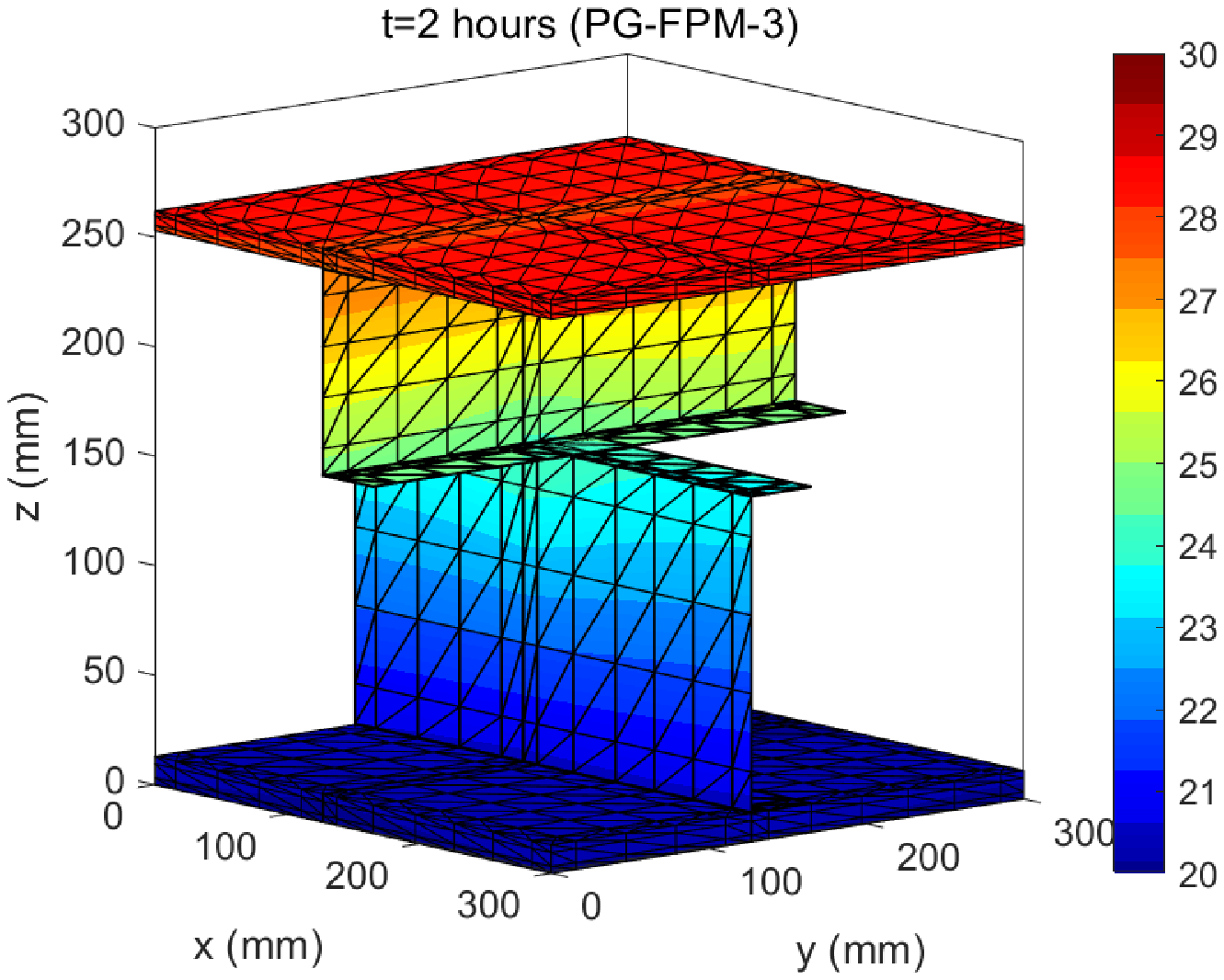}} 
    \subfigure[]{ 
    \label{fig:Ex28_SS} 
    \includegraphics[width=0.48\textwidth]{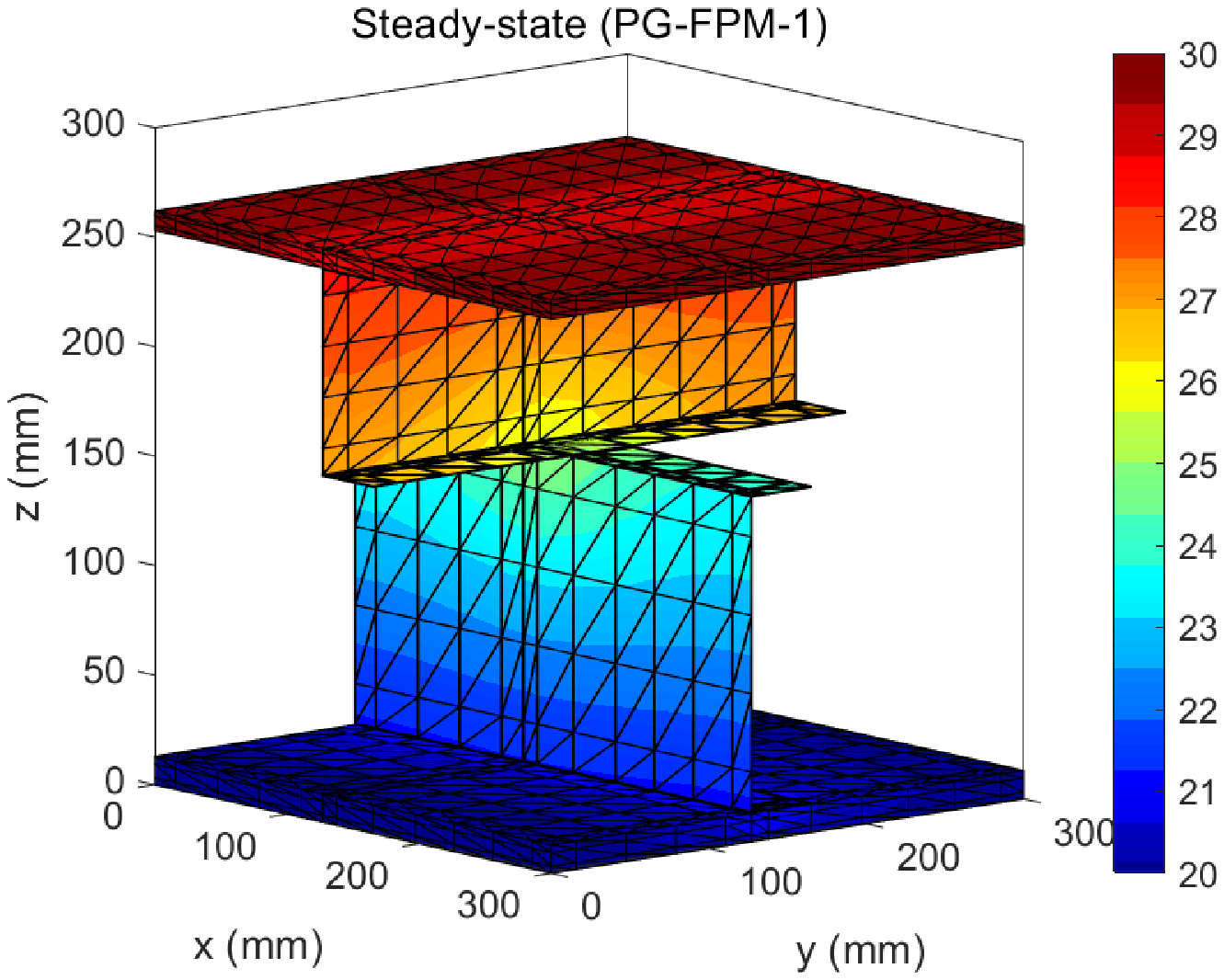}} 
  \caption{Ex. (2.8) – The problem and computed solutions. (a) the problem domain and boundary conditions. (b) transient temperature solution. (c) temperature distribution when $t=0.5$~hours (PG-FPM-1). (d) temperature distribution when $t=1.0$~hour (PG-FPM-2).  (e) temperature distribution when $t=2.0$~hours (PG-FPM-3). (f) steady-state result (PG-FPM-1).} 
  \label{fig:Ex28} 
\end{figure}

In the last two examples, practical problems are under study. Multiple materials and complicated geometries are considered. Ex.~(2.8) is in a wall with crossed U-girders. As shown in Fig.~\ref{fig:Ex28_BC}, the wall is consisted of two gypsum wallboards, two steel crossed U-girders and insulation materials (not presented in the sketch). The U-girders are separated by $300~\mathrm{mm}$. Thus, we focus on the heat conduction in a $300~\mathrm{mm} \times 300~\mathrm{mm} \times 262~\mathrm{mm}$ cell of the wall. The material properties and boundary conditions are listed in Table~\ref{table:Ex28-M} and Table~\ref{table:Ex28_BC}. All the lateral surfaces are symmetric. The initial condition is $u (x, y, z, t) = 20\mathrm{^\circ C}$.

\begin{table}[htbp]
\caption{Material properties in Ex.~(2.8).}
\centering
{
\begin{tabular*}{500pt}{@{\extracolsep\fill}lccc@{\extracolsep\fill}}
\toprule
\textbf{Material} & $\rho$ $(\mathrm{kg/m^3})$ &  $c$ $(\times 10^3~ \mathrm{J/(kg ^\circ C)})$ &  $k$ $(\mathrm{W/ (m^{2 \circ} C)})$ \\
\midrule
gypsum & 2300 & 1.09 & 0.22 \\
steel & 7800 & 0.50 & 60 \\
insulation & 1.29 & 1.01 & 0.036 \\
\bottomrule
\end{tabular*}}
\label{table:Ex28-M}
\end{table}

\begin{table}[htbp]
\caption{Robin boundary conditions in Ex.~(2.8).}
\centering
{
\begin{tabular*}{500pt}{@{\extracolsep\fill}lcc@{\extracolsep\fill}}
\toprule
\textbf{boundary condition} & $\widetilde{u}_R$ $(\mathrm{^\circ C})$ & $h$ $(\mathrm{W/(m^{2 \circ} C)})$ \\
\midrule
outdoor ($z = 0~\mathrm{mm}$) & 20 & 25 \\
indoor ($z = 262~\mathrm{mm}$) & 30 & 7.7 \\
\bottomrule
\end{tabular*}}
\label{table:Ex28_BC}
\end{table}

A total of 3380 Points with hexahedron partition are exploited. Due to the uneven variation of material properties, the Points are distributed unevenly, with a higher density in the gypsum and steel. The temperature at three representative points A, B and C (seen in Fig.~\ref{fig:Ex28_BC}) are monitored. The computed time variation of temperature on these points are presented in Fig.~\ref{fig:Ex28_Trans}. FEM solution achieved by ABAQUS with 9702 DC3D8 elements is also given as a comparison. As can be seen, a good agreement is observed between the ABAQUS result and all the FPM and PG-FPM solutions. Figure~\ref{fig:Ex28_S005_01} – \ref{fig:Ex28_SS} exhibit the temperature distribution in the gypsum wallboards and U-girders at $t = $ 0.5, 1, 2 and 10~hours acquired by multiple PG-FPM approaches. The transient temperature solution approaches the steady-state result gradually. The corresponding computational parameters and time costs of the different approaches are shown in Table~\ref{table:Ex28}. Note that the computational time for the collocation method (PG-FPM-1) decreases remarkably when the penalty parameter $\eta_1 = 0$. Nevertheless, a small positive $\eta_1$ helps to improve the accuracy of the method. In practice, the collocation method with $\eta_1 = 0$ can be employed for a rough estimate in the heat conduction analysis, which costs only one third of the computational time of the original Galerkin FPM. Whereas when higher accuracy is required, the finite volume method is the best choice. Note that the original Galerkin FPM and PG-FPM-2 / 3 based on linear trial functions are not stable with zero penalty parameter $\eta_1$. In the singular solution method (PG-FPM-3), 12 integration points are adopted in each subdomain in this example, which lead to a full Jacobian matrix and lower efficiency compared to all the other approaches.

\begin{table}[htbp]
\caption{Computational times of the FPM and PG-FPMs in solving Ex.~(2.8).}
\centering
{
\begin{tabular*}{500pt}{@{\extracolsep\fill}lccc@{\extracolsep\fill}}
\toprule
\textbf{Method} & \tabincell{c}{\textbf{Computational} \textbf{parameters}} & $N_{band} (\mathbf{K})$ \& $N_{band} (\mathbf{C})$ & \tabincell{c}{\textbf{Computational} \textbf{time (s)}} \\
\midrule
FPM & $\eta_1 = 10$, $\eta_2 = 20$ & 53 \& 1 & 28.7 \\
\midrule
\multirow{2}*{PG-FPM-1} & $\eta_1 = 10$, $\eta_2 = 20$ & 21 \& 1 & 25.1 \\
~ & $\eta_1 = 0$ , $\eta_2 = 20$ & 7 \& 1 & 8.5 \\
\midrule
PG-FPM-2 & $\eta_1 = 10$, $\eta_2 = 20$ & 21 \& 1 & 15.4 \\
\midrule
PG-FPM-3 & $\eta_1 = 10$, $\eta_2 = 20$ & 22 \& 7 & 119.1 \\
\bottomrule
\end{tabular*}}
\label{table:Ex28}
\end{table}

\begin{figure}[htbp] 
  \centering 
   \subfigure[]{ 
    \label{fig:Ex29_BC} 
    \includegraphics[width=0.48\textwidth]{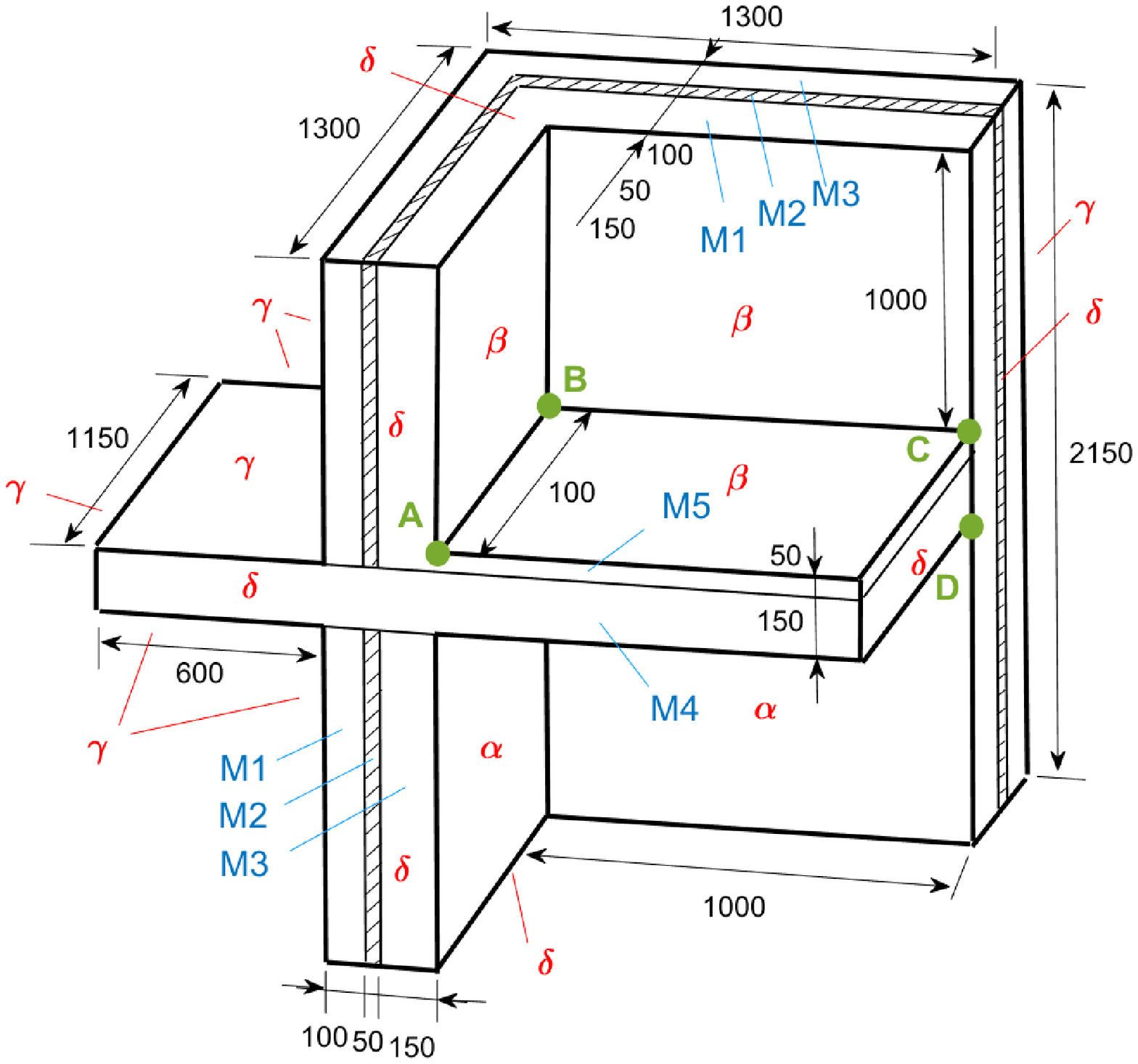}}  
   \subfigure[]{ 
    \label{fig:Ex29_Trans} 
    \includegraphics[width=0.48\textwidth]{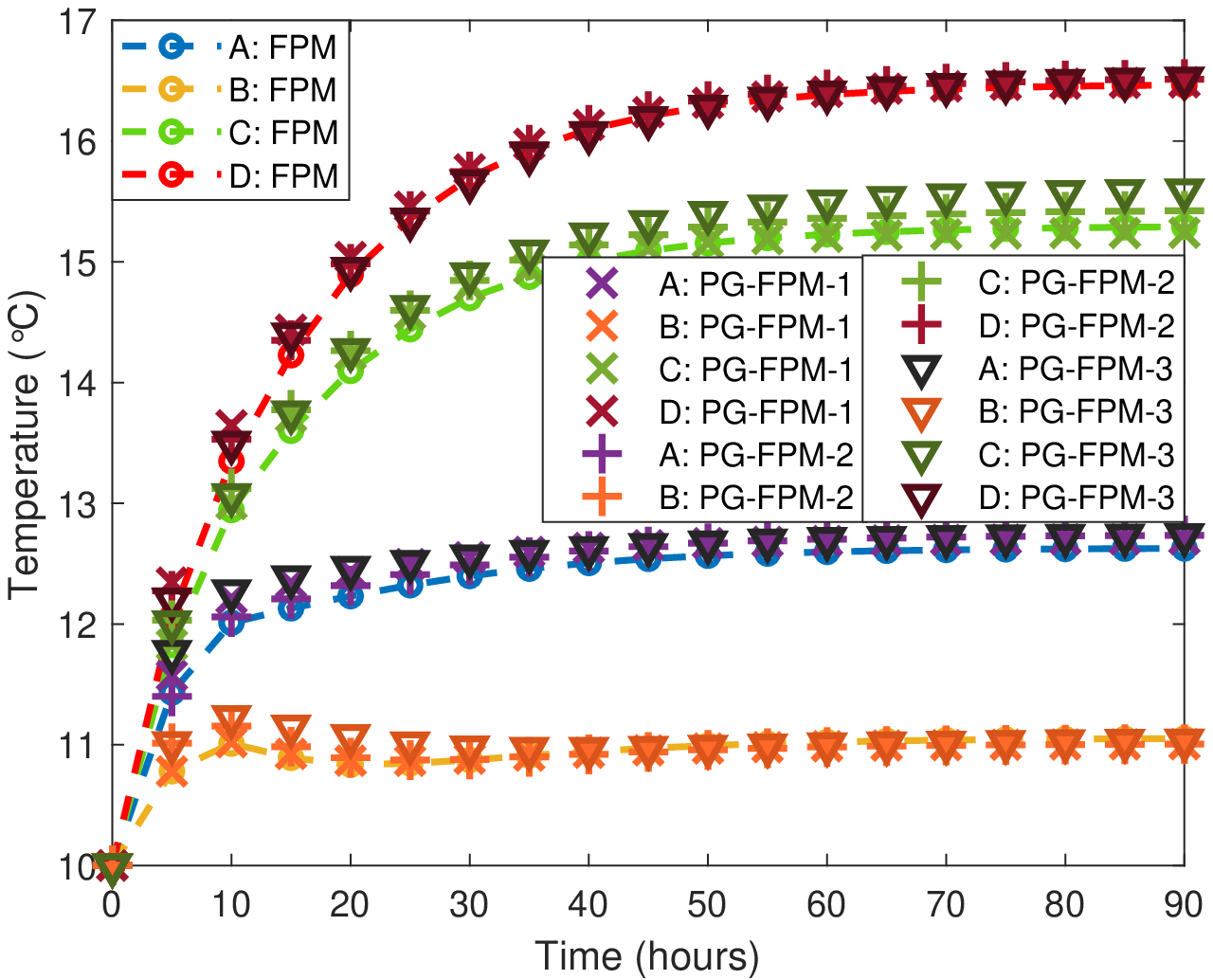}}  
    \subfigure[]{ 
    \label{fig:Ex29_S05_00} 
    \includegraphics[width=0.48\textwidth]{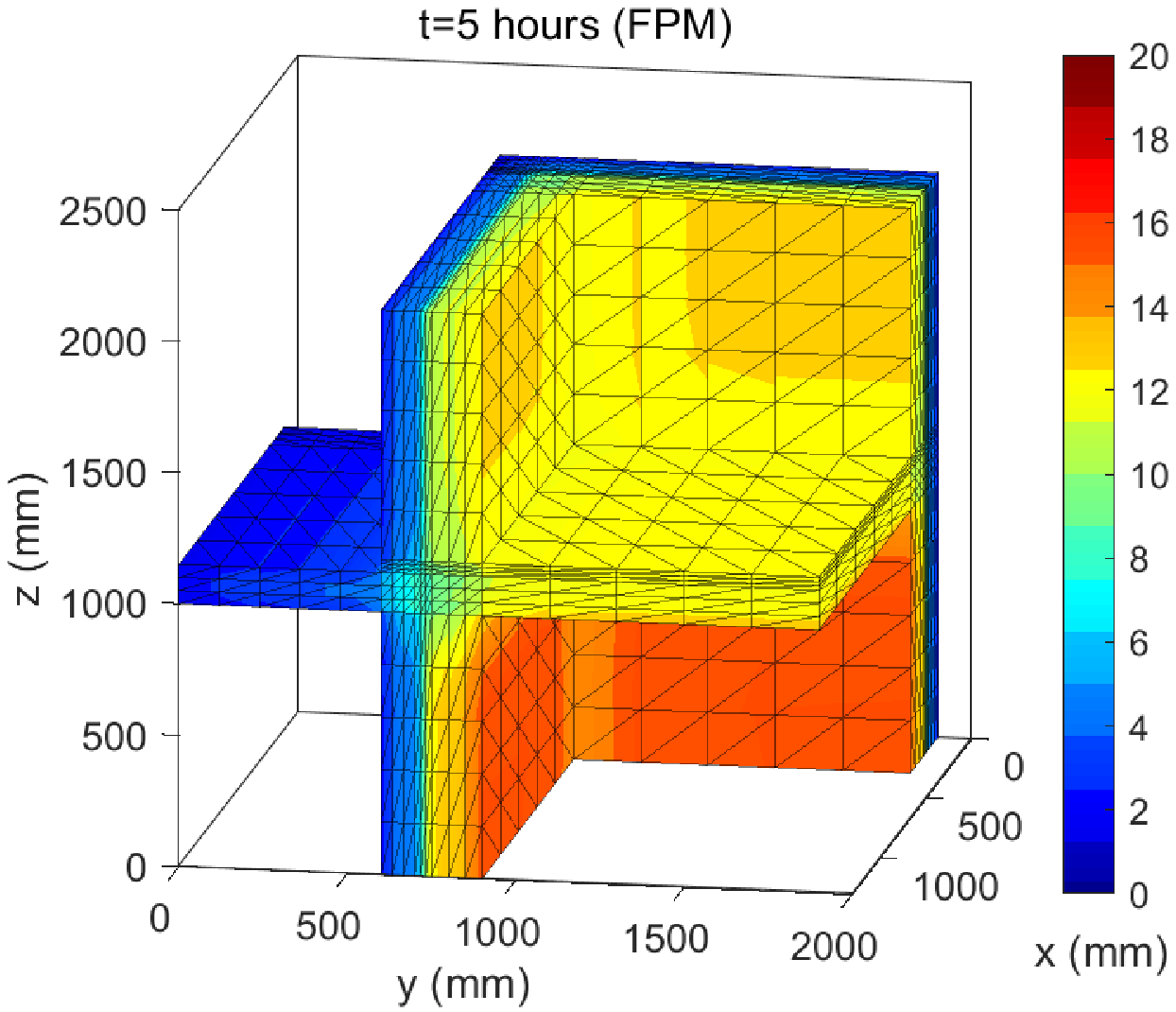}}  
    \subfigure[]{ 
    \label{fig:Ex29_S10_01} 
    \includegraphics[width=0.48\textwidth]{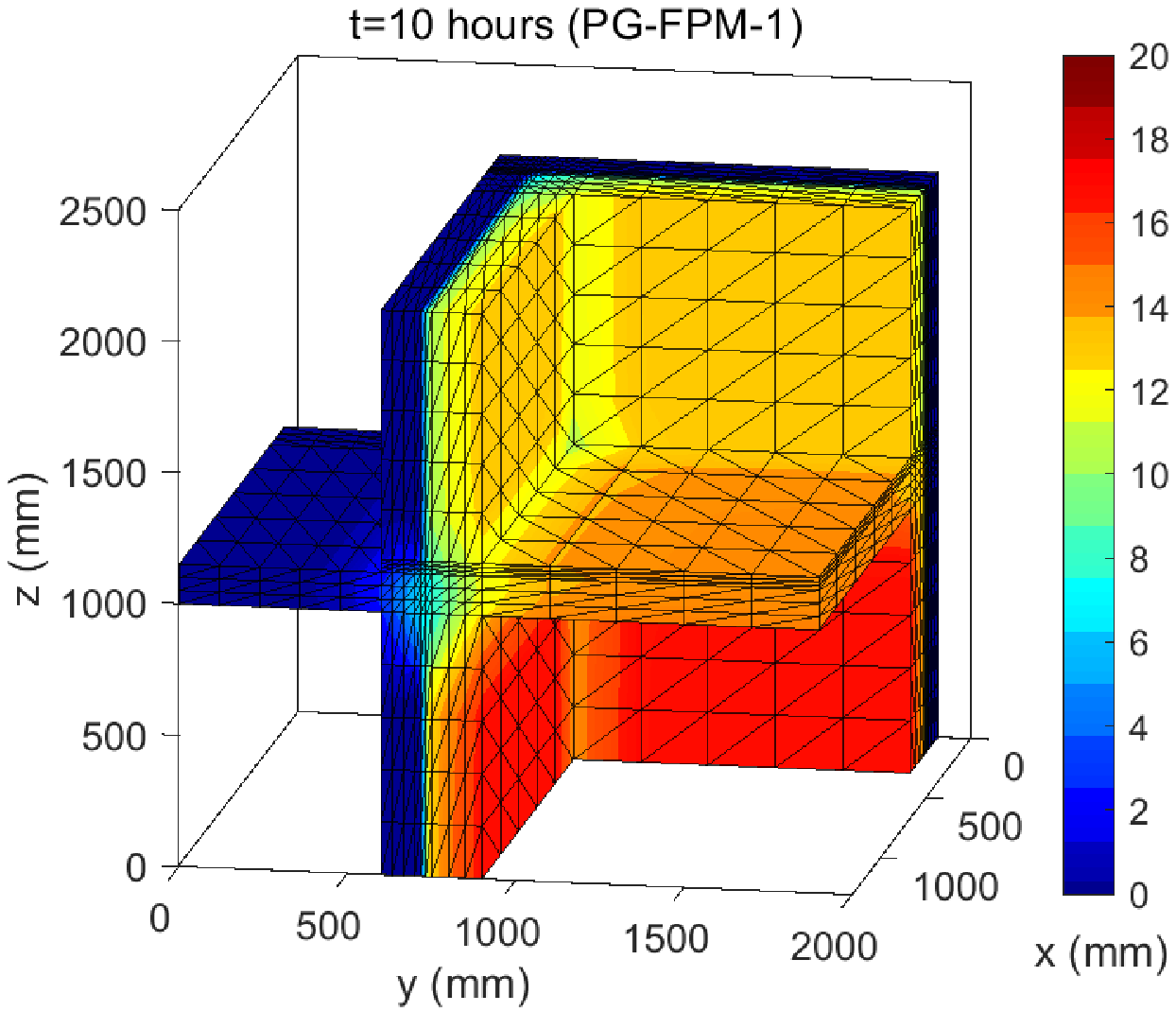}} 
    \subfigure[]{ 
    \label{fig:Ex29_S15_02} 
    \includegraphics[width=0.48\textwidth]{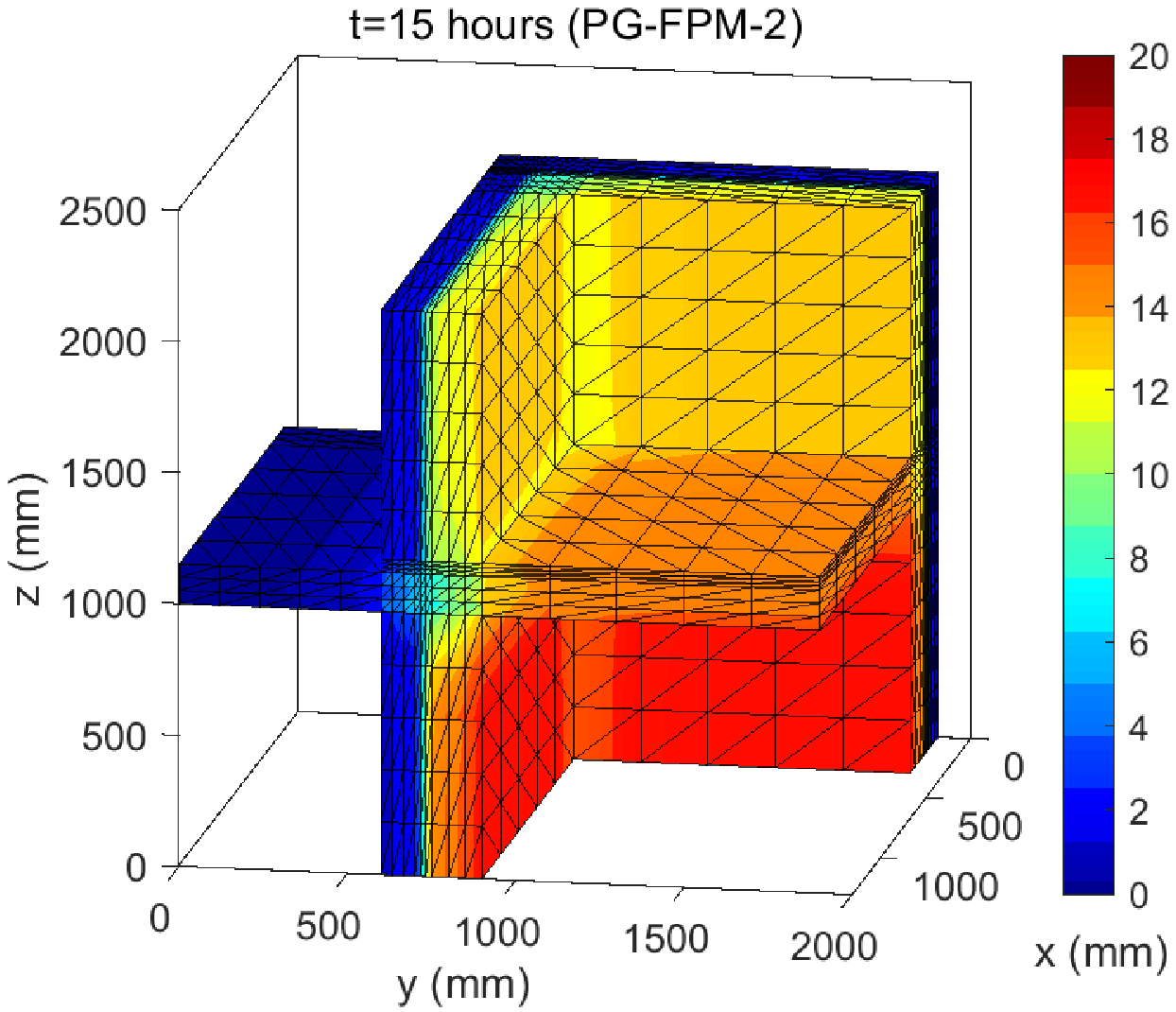}} 
    \subfigure[]{ 
    \label{fig:Ex29_SS_03} 
    \includegraphics[width=0.48\textwidth]{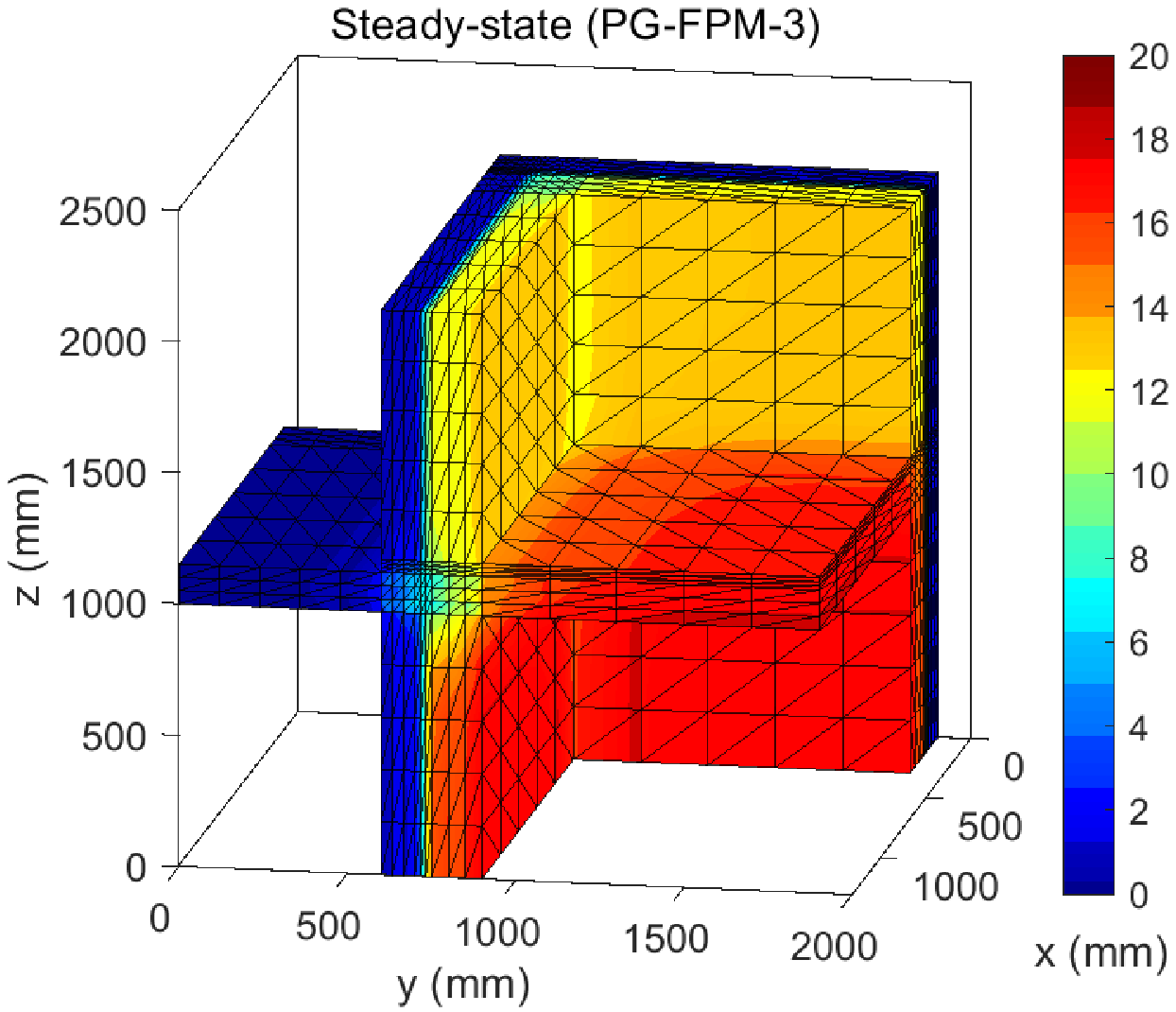}} 
  \caption{Ex. (2.9) – The problem and computed solutions. (a) the problem domain and boundary conditions. (b) transient temperature solution. (c) temperature distribution when $t=5$~hours (FPM). (d) temperature distribution when $t=10$~hour (PG-FPM-1).  (e) temperature distribution when $t=15$~hours (PG-FPM-2). (f) steady-state result (PG-FPM-3).} 
  \label{fig:Ex29} 
\end{figure}

In the last example, we study the heat conduction through a wall corner. The geometry and material distribution in the corner are shown in Fig.~\ref{fig:Ex29_BC}. The properties of the five different materials are listed in Table~\ref{table:EX29-M}. And the four kinds of boundary conditions are presented in Table~\ref{table:EX29-BC}, in which $\alpha$, $\beta$ and $\gamma$ are under Robin boundary conditions, and $\delta$ is adiabatic. The initial condition is $u (x, y, z, t) = 20\mathrm{^\circ C}$.

\begin{table}[htbp]
\caption{Material properties in Ex.~(2.9).}
\centering
{
\begin{tabular*}{500pt}{@{\extracolsep\fill}lccc@{\extracolsep\fill}}
\toprule
\textbf{Material} & $\rho$ $(\mathrm{kg/m^3})$ & $c$ $(\times 10^3~ \mathrm{J/(kg ^\circ C)})$ &  $k$ $(\mathrm{W/ (m^{2 \circ} C)})$ \\
\midrule
M1 & 849 & 0.9 & 0.7 \\
M2 &  80 & 0.84 & 0.04 \\
M3 & 2000 &  0.8& 1.0 \\
M4 & 2711 & 0.88 & 2.5 \\
M5 & 2400 & 0.96 & 1.0 \\
\bottomrule
\end{tabular*}}
\label{table:EX29-M}
\end{table}

\begin{table}[htbp]
\caption{Boundary conditions in Ex.~(2.9).}
\centering
{
\begin{tabular*}{500pt}{@{\extracolsep\fill}lcc@{\extracolsep\fill}}
\toprule
\textbf{Boundary condition} & $\widetilde{u}_R$ $(\mathrm{^\circ C})$ & $h$ $(\mathrm{W/(m^{2 \circ} C)})$ \\
\midrule
$\alpha$ & 20 & 5 \\
$\beta$ & 15 & 5 \\
$\gamma$ & 0 & 20 \\
$\delta$ & -- & 0 \text{(adiabatic)} \\
\bottomrule
\end{tabular*}}
\label{table:EX29-BC}
\end{table}

3006 Points are employed in the analysis. Figure~\ref{fig:Ex29_Trans} shows the time variation of temperatures on four representative points achieved by FPM and PG-FPMs. The result approaches steady state as time increases, and is finally consistent with the data shown in the European standard (CEN, 1995) \cite{ISO10211, Blomberg1996}. The spatial temperature distribution in the corner when $t = $ 5, 10, 15 and 90~hours are presented in Fig.~\ref{fig:Ex29_S05_00} – \ref{fig:Ex29_SS_03}. The solutions achieved by multiple Galerkin FPM and PG-FPM approaches are approximately identical. Hence the result of only one approach is presented at each time frame. In Table~\ref{table:Ex29}, the computational times of these PG-FPM approaches are listed. The same as Ex.~(2.8), the collocation method (PG-FPM-1) with $\eta_1 = 0$ is the most efficient approach, while the finite volume method (PG-FPM-2) has the highest efficiency when a more accurate solution is required.

\begin{table}[htbp]
\caption{Computational times of the FPM and PG-FPMs in solving Ex.~(2.9).}
\centering
{
\begin{tabular*}{500pt}{@{\extracolsep\fill}lccc@{\extracolsep\fill}}
\toprule
\textbf{Method} & \tabincell{c}{\textbf{Computational} \textbf{parameters}} & $N_{band} (\mathbf{K})$ \& $N_{band} (\mathbf{C})$ & \tabincell{c}{\textbf{Computational} \textbf{time (s)}} \\
\midrule
FPM & $\eta_1 = 0.1$, $\eta_2 = 20$ & 52 \& 1 & 26.0 \\
\midrule
\multirow{2}*{PG-FPM-1} & $\eta_1 = 0.1$, $\eta_2 = 20$ & 24 \& 1 & 15.7 \\
~  & $\eta_1 = 0$  , $\eta_2 = 20$ & 7 \& 1 & 10.6 \\
\midrule
PG-FPM-2 & $\eta_1 = 0.1$, $\eta_2 = 20$ & 21 \& 1 & 14.5 \\
\midrule
PG-FPM-3 & $\eta_1 = 0.1$, $\eta_2 = 20$ & 22 \& 7 & 97.0 \\
\bottomrule
\end{tabular*}}
\label{table:Ex29}
\end{table}

\section{Discussion of the computational parameters} \label{sec:PS}

\begin{figure}[htbp] 
  \centering 
   \subfigure{ 
    \label{fig:Ex11_e01} 
    \includegraphics[width=0.48\textwidth]{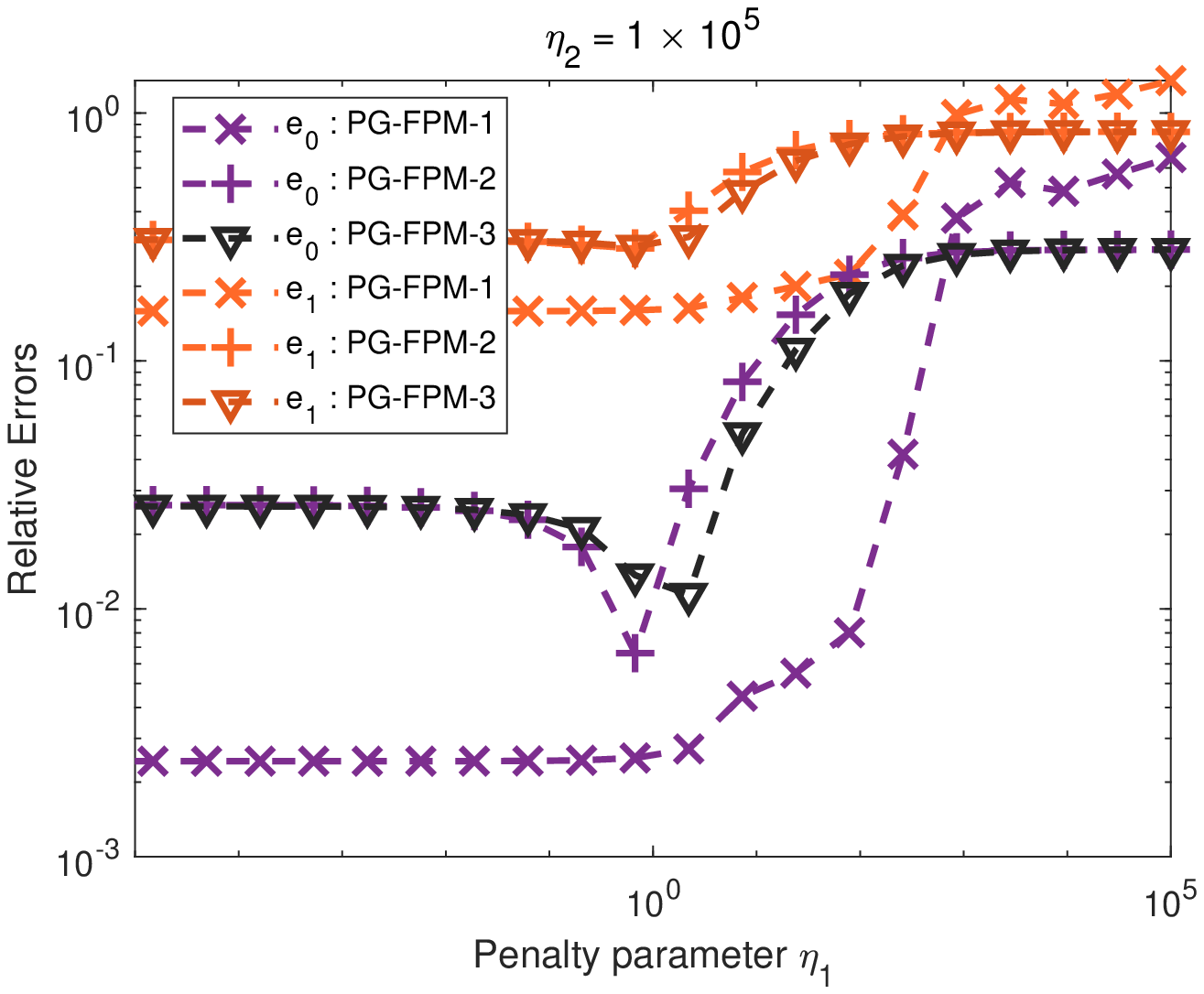}}  
    \subfigure{ 
    \label{fig:Ex21_e01} 
    \includegraphics[width=0.48\textwidth]{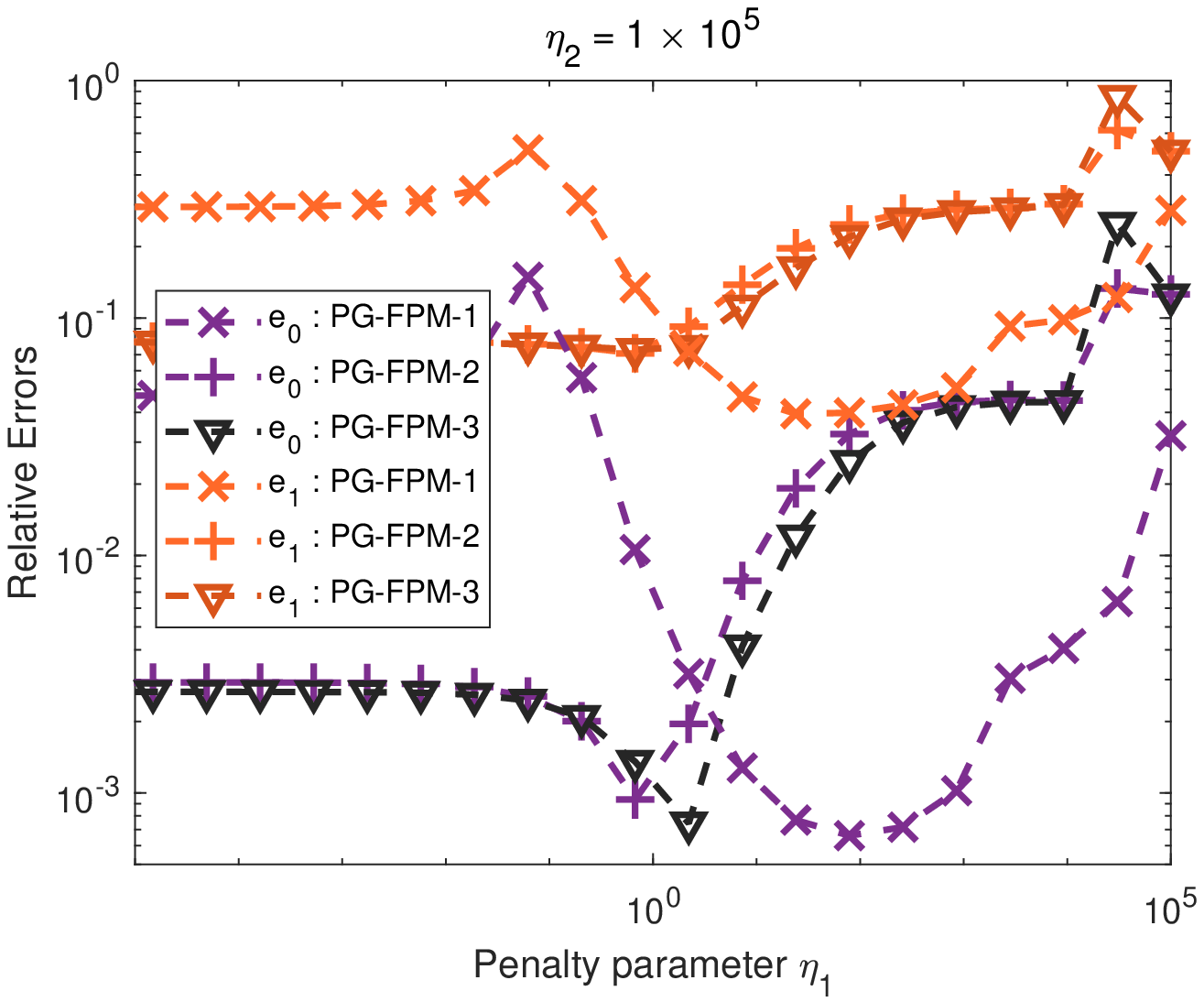}}   
    \addtocounter{subfigure}{-2}
   \subfigure[]{ 
    \label{fig:Ex11_e02} 
    \includegraphics[width=0.48\textwidth]{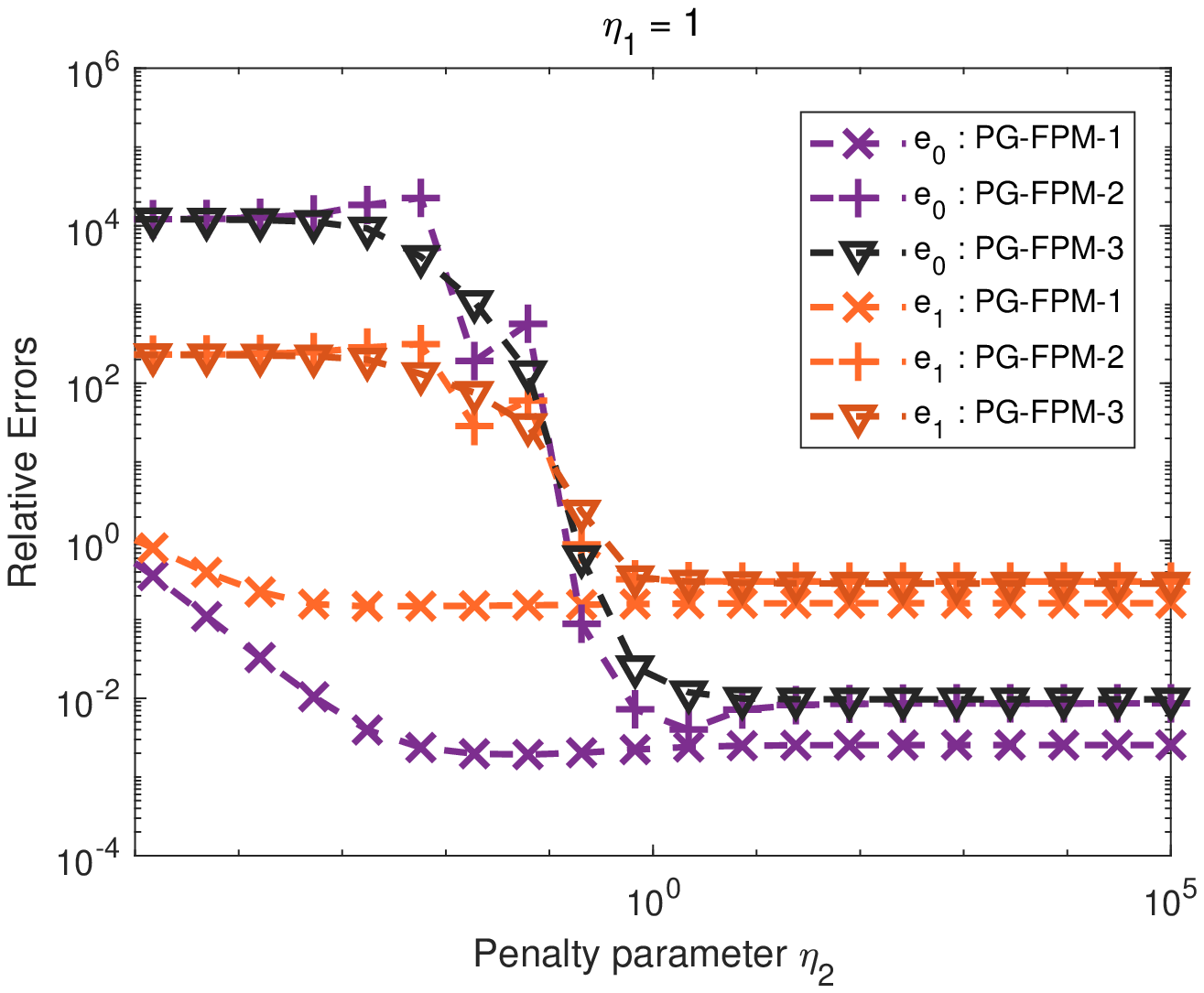}} 
    \subfigure[]{ 
    \label{fig:Ex21_e02} 
    \includegraphics[width=0.48\textwidth]{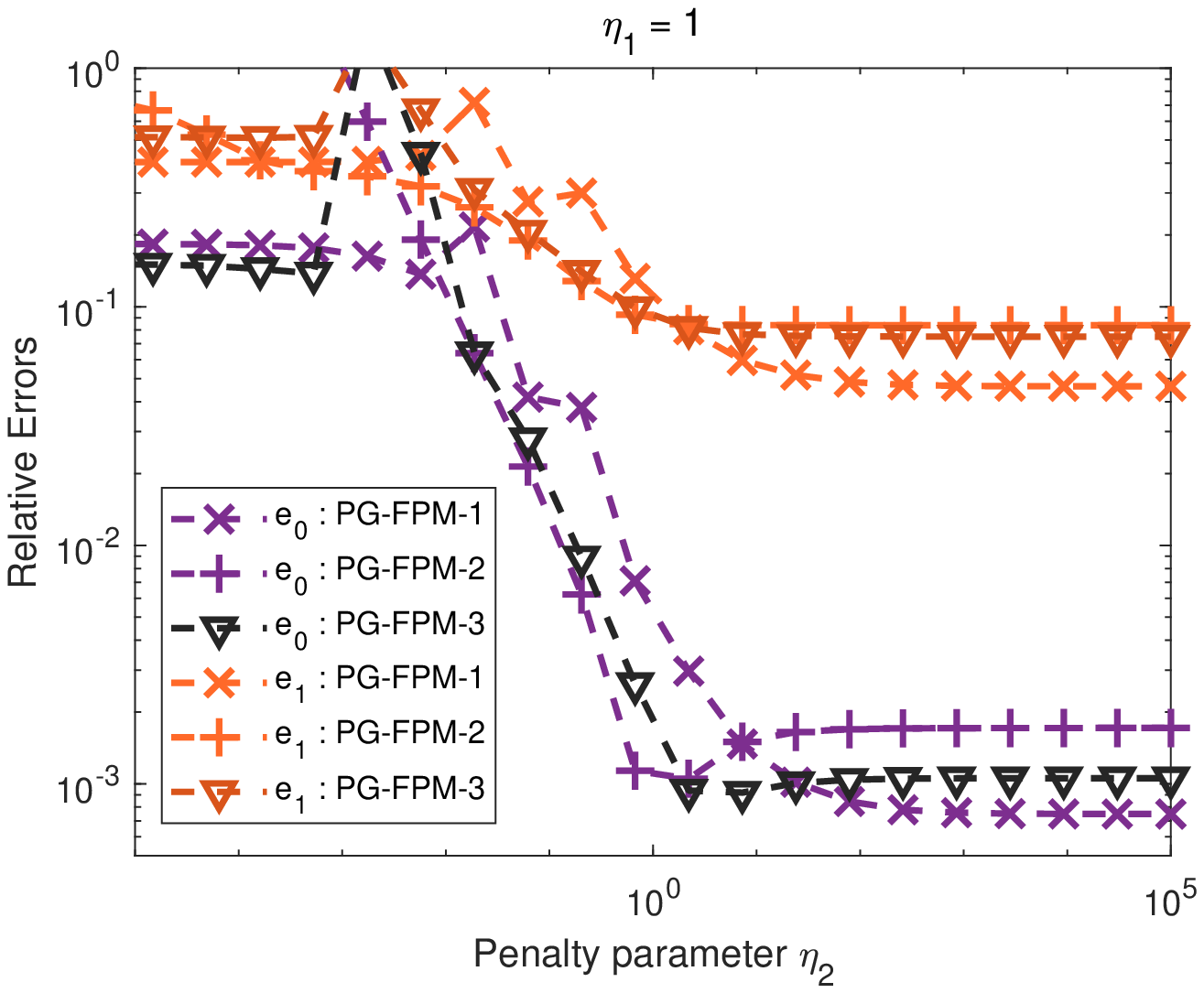}}  
  \caption{Parametric studies on $\eta_1$ and $\eta_2$. (a) 2D case: Ex.~(1.1). (b) 3D case: Ex.~(2.1).} 
  \label{fig:PS} 
\end{figure}

In the last part of this paper, we give a parametric study on the penalty parameters $\eta_1$ and $\eta_2$ for the proposed PG-FPM approaches. In 2D case, Ex.~(1.1) is reconsidered. With 600 Fragile Points, when the penalty parameters varies from $10^{-5}$ to $10^5$, the corresponding relative errors are exhibited in Fig.~\ref{fig:Ex11_e02}. As can be seen, for the finite volume method (PG-FPM-2) and singular solution method (PG-FPM-3), the solution is only stable and accurate when $0.2 < \eta_1 < 2$. Whereas for the collocation method, $\eta_1$ can be as small as zero. There is no upper limit for $\eta_2$ which enforces the boundary conditions, while the lower limit is approximate 1. Similarly, for 3D problems, the parametric study is carried out for Ex.~(2.1). With 1000 Fragile Points, the corresponding results are shown in Fig.~\ref{fig:Ex21_e02}. After taking the other examples into consideration, The recommended ranges for $\eta_1$ and $\eta_2$ in 2D and 3D analysis are given in Table~\ref{table:Range}. Note that the recommended values for different PG-FPM approaches can be different.

\begin{table}[htbp]
\caption{Recommended ranges of the penalty parameters.}
\centering
{
\begin{tabular*}{500pt}{@{\extracolsep\fill}lccc@{\extracolsep\fill}}
\toprule
\textbf{Method} & \textbf{Problem dimension} & $\eta_1$ & $\eta_2$ \\
\midrule
\multirow{2}*{FPM} & 2D & $1$ –  $50$ & $> 50$ \\
~  & 3D & $0.1$ –  $50$ & $50$ – $10^4$\\
\midrule
\multirow{2}*{PG-FPM-1} & 2D & $0$ –  $10$ & $> 1$ \\
~  & 3D & $0$ –  $10^3$ & $> 10$\\
\midrule
\multirow{2}*{PG-FPM-2} & 2D & $0.1$ –  $10$ & $> 10$ \\
~  & 3D & $0.1$ –  $10$ & $> 10$\\
\midrule
\multirow{2}*{PG-FPM-3} & 2D & $0.1$ –  $10$ & $> 10$ \\
~  & 3D & $0.1$ –  $10$ & $> 10$\\
\bottomrule
\end{tabular*}}
\label{table:Range}
\end{table}

The recommended range for the constant parameter $c$ in the collocation method (PG-FPM-1) is $1$ – $10$, with a recommended value 4 for 2D problems and 10 for 3D problems. All the parameters discussed in this section are nondimensional.

\section{Conclusion}

In the current work, three Fragile Points Methods based on Petrov-Galerkin weak-forms (PG-FPMs) are developed. As improved versions of the original Galerkin FPM, the trial functions in the PG-FPMs are still chosen to be local, polynomial, and piecewise-continuous. A modified local Radial Basis Function-based Differential Quadrature (RBF-DQ) method is introduced to approximate the first and higher derivatives at each Fragile Point. With different test functions, the three proposed PG-FPMs are also named as the collocation method (PG-FPM-1), the finite volume method (PG-FPM-2) and the singular solution method (PG-FPM-3) according to their features. Table~\ref{table:Comp_2} illustrates a comparison between the original Galerkin FPM and the three proposed PG-FPMs.

The PG-FPM approaches keep the advantages of the Galerkin FPM including simple numerical integration, trial functions with Delta function property, free of locking and mesh distortion problems, etc., which are superior to the EFG, MLPG and other methods in previous literatures. Meanwhile, all the three proposed approaches simplify the computing process of the original Galerkin FPM. A number of numerical results are given to validate the accuracy and efficiency of the three proposed PG-FPM approaches. 2D and 3D examples with Dirichlet, Neumann, Robin and symmetric boundary conditions are carried out. Both functionally graded and composite materials are considered. The computed solutions achieved by the three PG-FPM approaches are compared with each other, as well as solutions of the Galerkin FPM, FEM, and analytical solution (if applicable).

The PG-FPM-1 with Dirac delta function as the test function leads to collocation equations at each Point. The approach has advantages in transient analysis. As the heat capacity matrix is always diagonal, the PG-FPM-1 has the same computing efficiency with arbitrary point distributions and partitions. Moreover, for problems with high DoFs, the collocation method with penalty parameter $\eta_1 = 0$ gives a rough estimate of the solution costing as low as one-third to one-half of the computational time compared to the original FPM.

The PG-FPM-2 with the Heaviside step function as the test function is analogous to the conventional Finite Volume Method (FVM) yet the conservation law is not strictly satisfied but enforced by Interior Penalty Flux Corrections. The PG-FPM-2 approach is the best choice for a balance between accuracy and efficiency. With Fragile Points placed at the centroid of each subdomain, the finite volume method achieves similar or better accuracy as the original Galerkin FPM while saving 25\% – 50\% of the computational time. The method is also the most efficient in steady-state analysis.

The PG-FPM-3 with the local fundamental solution as the test function has vanishing integrals in subdomains in the weak-form formulation. The approach is recommended for steady-state analysis only. As more than one integration points may be required in each subdomain, the PG-FPM-3 has an unsatisfactory efficiency in analyzing transient problems. Yet its performance is as good as the finite volume method (PG-FPM-2) in steady-state analysis.

The recommended ranges of the nondimensional computational parameters utilized in the PG-FPMs are given at last. We conclude that, with suitable computational parameters, all the three proposed PG-FPM approaches, especially the finite volume method (PG-FPM-2) can improve the performance of the original Galerkin FPM, and are considerably superior as compared to heat conduction analyzing approaches in earlier literatures.

\begin{table}[htbp]
\caption{A comparison of the original FPM and the three presented PG-FPMs.}
\centering
{
\begin{tabular}{ m{83pt}  m{83pt}<{\centering}  m{83pt}<{\centering}  m{83pt}<{\centering}   m{83pt}<{\centering}  }
\toprule
\textbf{Method} & \textbf{FPM} & \textbf{PG-FPM - 1 (collocation method)} & \textbf{PG-FPM - 2 (finite volume method)} & \textbf{PG-FPM - 3 (singular solution method)} \\
\midrule
Trial functions & Discontinuous linear & Discontinuous quadratic & Discontinuous linear & Discontinuous linear \\
\midrule
Test functions & As above & Dirac delta function & Heaviside step function & Local fundamental solution \\
\midrule
Local approximation method & GFD, RBF-DQ & RBF-DQ & GFD, RBF-DQ & GFD, RBF-DQ \\
\midrule
Support of each Fragile points & Nearest neighboring & Nearest and second neighboring & Nearest neighboring  & Nearest neighboring \\
\midrule
Weak-form & Galerkin & Petrov-Galerkin & Petrov-Galerkin & Petrov-Galerkin  \\
\midrule
Weak-form integration & Very simple (one-point quadrature) & Very simple (one-point quadrature) & Very simple (one-point quadrature) & Simple (Gaussian quadrature)  \\
\midrule
Heat capacity matrix ($\mathbf{C}$) & Symmetric, sparse & Diagonal & Asymmetric, sparse & Asymmetric, sparse \\
\midrule
Thermal conductivity matrix ($\mathbf{K}$) & Symmetric, sparse & Asymmetric, sparse & Asymmetric, highly sparse & Asymmetric, highly sparse \\
\midrule
Jacobian matrix ($\mathbf{J} = - \mathbf{C}^\mathrm{-1} \mathbf{K}$) & Symmetric, sparse or full & Asymmetric, sparse & Asymmetric, sparse or full & Asymmetric, sparse or full \\
\bottomrule
\end{tabular}}
\label{table:Comp_2}
\end{table}

\section*{Acknowledgment}

We thank Prof. Leiting Dong (Beihang University) for several stimulating conversations on these topics. We thankfully acknowledge the financial support for Dr. Guan’s work, provided through the funding for Professor Atluri’s Presidential Chair at TTU.

\bibliography{Part_1}

\end{document}